\documentclass[12pt]{article}
\usepackage{amssymb,amsfonts}
\usepackage{amsmath}
\newtheorem{theorem}{Theorem}[part]
\newtheorem{lemma}{Lemma}[part]
\newtheorem{prop}{Proposition}[part]
\newtheorem{corollary}{Corollary}[part]
\begin{document}
\title{Geometry of Optimal Control Problems and Hamiltonian
Systems}
\author{A.~A.~Agrachev}
\date{}
\maketitle

\section*{Preface}

These notes are based on the mini-course given in June 2004 in
Cetraro, Italy, in the frame of a C.I.M.E. school. Of course, they
contain much more material that I could present in the 6 hours
course. The goal was to give an idea of the general variational and
dynamical nature of nice and powerful concepts and results mainly
known in the narrow framework of Riemannian Geometry. This
concerns Jacobi fields, Morse's index formula, Levi
Civita connection, Riemannian curvature and related topics.

I tried to make the presentation as light as possible: gave more
details in smooth regular situations and referred to the
literature in more complicated cases. There is an evidence that the
results described in the notes and treated in technical papers
we refer to are just parts of a united beautiful subject to be
discovered on the crossroads of Differential Geometry, Dynamical
Systems, and Optimal Control Theory. I will be happy if the
course and the notes encourage some young ambitious researchers to
take part in the discovery and exploration of this subject.

\smallskip\noindent {\sl Acknowledgments.} I would like to express
my gratitude to
Professor Gamkrelidze for his permanent interest to this topic and
many inspiring discussions and to thank participants of the school for
their surprising and encouraging will to work in the relaxing
atmosphere of the Mediterranean resort.

\tableofcontents

\part{Lagrange multipliers' geometry}

\section{Smooth optimal control problems}
In these lectures we discuss some geometric constructions
and results emerged from the investigation of smooth optimal
control problems. We'll consider problems with integral costs and
fixed endpoints. A standard formulation of such a problem is as
follows:\quad Minimize a functional
$$
J^{t_1}_{t_0}(u(\cdot))=\int\limits_{t_0}^{t_1}\varphi(q(t),u(t))\,dt, \eqno
(1)
$$
where
$$
\dot q(t)=f(q(t),u(t)),\quad u(t)\in U, \quad \forall t\in[t_0,t_1],
\eqno (2)
$$
$q(t_0)=q_0, \ q(t_1)=q_1$. Here $q(t)\in\mathbb R^n,\
U\subset\mathbb R^k$, a {\it control function} $u(\cdot)$ is
supposed to be measurable bounded while $q(\cdot)$ is
Lipschitzian; scalar function $\varphi$ and vector function $f$
are smooth. A pair $(u(\cdot),q(\cdot))$ is called an {\it
admissible pair} if it satisfies differential equation (2) but may
violate the boundary conditions.

We usually assume that Optimal Control Theory generalizes
classical Calculus of Variations. Unfortunately, even the most
classical geometric variational problem, the length minimization
on a Riemannian manifold, cannot be presented in the just
described way. First of all, even simplest manifolds,
like spheres, are not domains in $\mathbb R^n$. This does not look
as a serious difficulty: we slightly generalize original
formulation of the optimal control problem assuming that $q(t)$
belongs to a smooth manifold $M$ instead of $\mathbb R^n$. Then
$\dot q(t)$ is a tangent vector to $M$ i.e. $\dot q(t)\in T_{q(t)}M$
and we assume that $f(q,u)\in T_qM,\ \forall q,u.$ Manifold $M$ is
called the {\it state space} of the optimal control problem.

Now we'll try to give a natural formulation of the length
minimization problem as an optimal control problem on a Riemannian
manifold $M$. Riemannian structure on $M$ is (by definition) a
family of Euclidean scalar products $\langle\cdot,\cdot\rangle_q$ on
$T_qM,\ q\in M$, smoothly depending on $q$. Let
$f_1(q),\ldots,f_n(q)$ be an orthonormal basis of $T_qM$ for the
Euclidean structure $\langle\cdot,\cdot\rangle_q$ selected in such
a way that $f_i(q)$ are smooth with respect to $q$. Then any
Lipschitzian curve on $M$ satisfies a differential equation of the
form:
$$
\dot q=\sum\limits_{i=1}^nu_i(t)f_i(q), \eqno (3)
$$
where $u_i(\cdot)$ are measurable bounded scalar functions. In
other words, any Lipschitzian curve on $M$ is an admissible
trajectory of the control system (3). The Riemannian length of the
tangent vector $\sum\limits_{i=1}^nu_if_i(q)$ is
$\left(\sum\limits_{i=1}^nu_i^2\right)^{1/2}$. Hence the
length of a trajectory of system (3) defined on the segment
$[t_0,t_1]$ is $\ell(u(\cdot))=
\int_{t_0}^{t_1}\left(\sum\limits_{i=1}^nu^2_i(t)\right)^{1/2}\,dt$.
Moreover, it is easy to derive from the Cauchy--Schwarz
inequality that the length minimization is equivalent to the
minimization of the functional
$J^{t_1}_{t_0}(u(\cdot))=\int_{t_0}^{t_1}\sum\limits_{i=1}^nu^2_i(t)\,dt$.
The length minimization problem is thus reduced to a specific
optimal control problem on the manifold of the form (1), (2).

Unfortunately, what I've just written was wrong. It would be
correct if we could select a smooth orthonormal frame
$f_i(q),\ q\in M,\ i=1,\ldots,n$. Of course, we can always do it
locally, in a coordinate neighborhood of $M$ but, in general, we
cannot do it globally. We cannot do it even on the 2-dimensional
sphere: you know very well that any continuous vector field on the
2-dimensional sphere vanishes somewhere. We thus need another more
flexible formulation of a smooth optimal control problem.

Recall that a {\it smooth locally trivial bundle} over $M$ is a
submersion \linebreak
$\pi:V\to M$, where all {\it fibers} $V_q=\pi^{-1}(q)$
are diffeomorphic to each other and, moreover, any $q\in M$
possesses a neighborhood $O_q$ and a diffeomorphism
$\Phi_q:O_q\times V_q\to\pi^{-1}(O_q)$ such that
$\Phi_q(q',V_q)=V_{q'},\ \forall q'\in O_q$. In a less formal
language one can say that a smooth locally trivial bundle is a
smooth family of diffeomorphic manifolds $V_q$ (the fibers)
parametrized by the points of the manifold $M$ (the base). Typical
example is the tangent bundle $TM=\bigcup\limits_{q\in M}T_qM$
with the canonical projection $\pi$ sending $T_qM$ into $q$.

\medskip\noindent{\bf Definition.} A smooth control system with
the state space $M$ is a smooth mapping $f:V\to TM$, where $V$ is
a locally trivial bundle over $M$ and $f(V_q)\subset T_qM$ for any
fiber $V_q,\ q\in M$. An admissible pair is a
bounded\footnote{the term `bounded' means that the closure of the
image of the mapping is compact} measurable mapping
$v(\cdot):[t_0,t_1]\to V$ such that $t\mapsto\pi(v(t))=q(t)$ is a
Lipschitzian curve in $M$ and $\dot q(t)=f(v(t))$ for almost all
$t\in[t_0,t_1]$. Integral cost is a functional
$J^{t_1}_{t_0}(v(\cdot))=\int\limits_{t_0}^{t_1}\varphi(v(t))\,dt$, where
$\varphi$ is a smooth scalar function on $V$.

\medskip\noindent{\bf Remark.} The above more narrow
definition of an optimal control
problem on $M$ was related to the case of a {\it trivial bundle}
$V=M\times U,\ V_q=\{q\}\times U$. For the length minimization
problem we have $V=TM,\ f=\mathrm{Id},\ \varphi(v)=\langle
v,v\rangle_q,\ \forall v\in T_qM,\ q\in M$.

\medskip
Of course, any general smooth control system on the manifold $M$
 is locally equivalent to a standard control system on $\mathbb
R^n$. Indeed, any point $q\in M$ possesses a coordinate neighborhood
$O_q$ diffeomorphic to $\mathbb R^n$ and a mapping
$\Phi_q:O_q\times V_q\to\pi^{-1}(O_q)$ trivializing the
restriction of the bundle $V$ to $O_q$; moreover, the fiber $V_q$
can be embedded in $\mathbb R^k$ and thus serve as a set of
control parameters $U$.

Yes, working locally we do not obtain new systems with respect
of those in
$\mathbb R^n$. Nevertheless, general intrinsic definition is very
useful and instructive even for a purely local geometric analysis.
Indeed, we do not need to fix specific coordinates on $M$ and a
trivialization of $V$ when we study a control system defined in
the intrinsic way. A change of coordinates in $M$ is actually a
smooth transformation of the state space while a change of the
trivialization results in the feedback transformation of the
control system. This means that an intrinsically defined control
system represents actually the whole class of systems that are
equivalent with respect to smooth state and feedback
transformations. All information on the system obtained in the
intrinsic language is automatically invariant with respect to
smooth state and feedback transformations. And this is what any
geometric analysis intends to do: to study properties of the
object under consideration preserved by the natural transformation
group.

We denote by $L_\infty([t_0,t_1];V)$ the space of measurable
bounded mappings from $[t_0,t_1]$ to $V$ equipped with the
$L_\infty$-topology of the uniform convergence on a full measure
subset of $[t_0,t_1]$. If $V$ would an Euclidean space, then
$L_\infty([t_0,t_1];V)$ would have a structure of a Banach space.
Since $V$ is only a smooth manifold, then $L_\infty([t_0,t_1];V)$
possesses a natural structure of a smooth Banach manifold modeled on the
Banach space $L_\infty([t_0,t_1];\mathbb R^{\dim V})$.

Assume that $V\to M$ is a locally trivial bundle with the
$n$-dimensional base and $m$-dimensional fibers; then $V$ is an
$(n+m)$-dimensional manifold.

\begin{prop} Let $f:V\to TM$ be a smooth control system; then the
space $\mathcal V$ of admissible pairs of this system is a smooth
Banach submanifold of $L_\infty([t_0,t_1];V)$ modeled on
$\mathbb R^n\times L_\infty([t_0,t_1];\mathbb R^m)$.
\end{prop}
{\bf Proof.} Let $v(\cdot)$ be an admissible pair and
$q(t)=\pi(v(t)),\ t\in[t_0,t_1]$. There exists a Lipschitzian with
respect to $t$ family of local trivializations \linebreak
$R_t:O_{q(t)}\times U\to\pi^{-1}(O_{q(t)})$, where $U$ is
diffeomorphic to the fibers $V_q$. The construction of such a
family is a boring exercise which we omit.

Consider the system
$$
\dot q=f\circ R_t(q,u),\quad u\in U. \eqno (4)
$$
Let $v(t)=R_t(q(t),u(t))$; then $R_t,\ t_0\le t\le t_1,$
induces a diffeomorphism of an $L_\infty$-neighborhood of
$(q(\cdot),u(\cdot))$ in the space of admissible pairs for (4) on
a neighborhood of $v(\cdot)$ in $\mathcal V$. Now fix $\bar
t\in[t_0,t_1]$. For any $\hat q$ close enough to $q(\bar t)$ and any
$u'(\cdot)$ sufficiently close to $u(\cdot)$ in the
$L_\infty$-topology there exists a unique Lipschitzian path
$q'(\cdot)$ such that $\dot q'(t)=f\circ R_t(q'(t),u'(t))),\
t_0\le t\le t_1,\ q'(\bar t)=\hat q$; moreover the mapping $(\hat
q,u'(\cdot))\mapsto q'(\cdot)$ is smooth. In other words, the
Cartesian product of a neighborhood of $q(\bar t)$ in $M$ and a
neighborhood of $u(\cdot)$ in  $L_\infty([t_0,t_1],U)$ serves as
a coordinate chart for a neighborhood of $v(\cdot)$ in $\mathcal
V$. This finishes the proof since $M$ is an $n$-dimensional
manifold and $L_\infty([t_0,t_1],U)$ is a Banach manifold modeled
on $L_\infty([t_0,t_1],\mathbb R^m). \quad \square $

An important role in our study will be played by the
``evaluation mappings" $F_t:v(\cdot)\mapsto q(t)=\pi(v(t))$. It is
easy to show that $F_t$ is a smooth mapping from $\mathcal V$ to
$M$. Moreover, it follows from the proof of Proposition~I.1 that
$F_t$ is a submersion. Indeed, $q(t)=F_t(v(\cdot))$ is, in fact a
part of the coordinates of $v(\cdot)$ built in the proof (the
remaining part of the coordinates is the control $u(\cdot)$.

\section{Lagrange multipliers}
Smooth optimal control problem is a special case of the general
smooth conditional minimum problem on a Banach manifold $\mathcal
W$. The general problem consists of the minimization of a smooth
functional $J:\mathcal W\to\mathbb R$ on the level sets
$\Phi^{-1}(z)$ of a smooth mapping $\Phi:\mathcal W\to N$, where
$N$ is a finite-dimensional manifold. In the optimal control
problem we have $\mathcal W=\mathcal V,\ N=M\times M,\
\Phi=(F_{t_0},F_{t_1})$.

An efficient classical way to study the conditional minimum
problem is the Lagrange multipliers rule. Let us give a coordinate
free description of this rule.  Consider the mapping
$$
\bar\Phi=(J,\Phi):\mathcal W\to\mathbb R\times N,\quad
\bar\Phi(w)=(J(w),\Phi(w)),\ w\in\mathcal W.
$$
It is easy to see that any point of the local conditional minimum
or maximum (i.e. local minimum or maximum of $J$ on a level set of
$\Phi$) is a critical point of $\bar\Phi$. I recall that $w$ is a
critical point of $\bar\Phi$ if the differential
$D_w\bar\Phi:T_w\mathcal W\to T_{\bar\Phi(w)}\left(\mathbb R\times
N\right)$ is {\it not} a surjective mapping. Indeed, if
$D_w\bar\Phi$ would surjective then, according to the implicit
function theorem, the image $\bar\Phi(O_w)$ of an arbitrary
neighborhood $O_w$ of $w$ would contain a neighborhood of
$\bar\Phi(w)=(J(w),\Phi(w))$; in particular, this image would
contain an interval
$\left((J(w)-\varepsilon,J(w)+\varepsilon),\Phi(w)\right)$ that
contradicts the local conditional minimality or maximality of
$J(w)$.

The linear mapping $D_w\bar\Phi$ is not surjective if and only if
there exists a nonzero linear form $\bar\ell$ on
$T_{\bar\Phi(w)}\left(\mathbb R\times N\right)$ which annihilates
the image of $D_w\bar\Phi$. In other words, $\bar\ell
D_w\bar\Phi=0$, where $\bar\ell D_w\bar\Phi:T_w\mathcal
W\to\mathbb R$ is the composition of $D_w\bar\Phi$ and the linear
form $\bar\ell: T_{\bar\Phi(w)}\left(\mathbb R\times N\right)\to
\mathbb R$.

We have $T_{\bar\Phi(w)}\left(\mathbb R\times N\right)=\mathbb
R\times T_{\Phi(w)}N$. Linear forms on $\left(\mathbb R\times N\right)$
constitute the adjoint space $\left(\mathbb R\times N\right)^*=
\mathbb R\oplus T^*_{\Phi(w)}N$, where $T^*_{\Phi(w)}N$ is the
adjoint space of $T_{\Phi(w)}M$ (the {\it cotangent space} to $M$
at the point $\Phi(w)$). Hence $\ell=\nu\oplus\ell$, where
$\nu\in\mathbb R,\ \ell\in T^*_{\Phi(w)}N$ and
$$
\bar\ell D_w\bar\Phi=(\nu\oplus\ell)\left(d_wJ,D_w\Phi\right)=\nu
d_wJ+\ell D_w\Phi.
$$
We obtain the equation
$$
\nu d_wJ+\ell D_w\Phi=0. \eqno (5)
$$
This is the Lagrange multipliers rule: if $w$ is a local
conditional extremum, then there exists a nontrivial pair
$(\nu,\ell)$ such that equation (5) is satisfied. The pair
$(\nu,\ell)$ is never unique: indeed, if $\alpha$ is a nonzero
real number, then the pair $(\alpha\nu,\alpha\ell)$ is also
nontrivial and satisfies equation (5). So the pair is actually
defined up to a scalar multiplier; it is natural to treat this
pair as an element of the projective space
$\mathbb P\left(\mathbb R\oplus T^*_{\Phi(w)}N\right)$ rather
than an element of the linear space.

The pair $(\nu,\ell)$ which satisfies (5) is called the {\it Lagrange
multiplier} associated to the critical point $w$. The Lagrange
multiplier is called {\it normal} if $\nu\ne 0$ and abnormal if
$\nu=0$. In these lectures we consider only normal Lagrange
multipliers, they belong to a distinguished coordinate chart of
the projective space
$\mathbb P\left(\mathbb R\oplus T^*_{\Phi(w)}N\right)$.

Any normal Lagrange multiplier has a unique representative of the
form $(-1,\ell)$; then (5) is reduced to the equation
$$
\ell D_w\Phi=d_wJ. \eqno (6)
$$
The vector $\ell\in T^*_{\Phi(w)}N$ from equation (6) is also
called a normal Lagrange multiplier (along with $(-1,\ell)$).

\section{Extremals}
Now we apply the Lagrange multipliers rule to the optimal control
problem. We have $\Phi=\left(F_{t_0},F_{t_1}\right):\mathcal V\to
M\times M$. Let an admissible pair $v\in\mathcal V$ be a critical
point of the mapping $\left(J^{t_1}_{t_0},\Phi\right)$, the curve
 $q(t)=\pi(v(t)),\ t_0\le t\le t_1$ be the corresponding trajectory, and
$\ell\in T^*_{(q(t_0),q(t_1))}(M\times M)$ be a normal Lagrange
multiplier associated to $v(\cdot)$. Then
$$
\ell D_v\left(F_{t_0},F_{t_1}\right)=d_vJ^{t_1}_{t_0}. \eqno (7)
$$
We have $T^*_{(q(t_0),q(t_1))}(M\times M)=T^*_{q(t_0)}M\times
T^*_{q(t_1)}M$, hence $\ell$ can be presented in the form
$\ell=(-\lambda_{t_0},\lambda_{t_1})$, where $\lambda_{t_i}\in
T^*_{q(t_i)}M,\ i=0,1$. Equation (7) takes the form
$$
\lambda_{t_1}D_vF_{t_1}-\lambda_{t_0}D_vF_{t_0}=d_vJ^{t_1}_{t_0}.
 \eqno (8)
 $$
 Note that $\lambda_{t_1}$ in (8) is uniquely defined by
 $\lambda_{t_0}$ and $v$. Indeed, assume that
 $\lambda'_{t_1}D_vF_{t_1}-\lambda_{t_0}D_vF_{t_0}=d_vJ_{t_0}^{t_1}$
 for some $\lambda'_{t_1}\in T^*_{q(t_1)}M$. Then
 $(\lambda'_{t_1}-\lambda_{t_1})D_vF_{t_1}=0$. Recall that
 $F_{t_1}$ is a submersion, hence $D_vF_{t_1}$ is a surjective
 linear map and $\lambda'_{t_1}-\lambda_{t_1}=0$.

 \begin{prop} Equality (8) implies that for any $t\in[t_0,t_1]$
 there exists a unique $\lambda_t\in T^*_{q(t)}M$ such that
 $$
 \lambda_{t}D_vF_{t}-\lambda_{t_0}D_vF_{t_0}=d_vJ^t_{t_0} \eqno
 (9)
 $$
 and $\lambda_t$ is Lipschitzian with respect to $t$.
 \end{prop}
 {\bf Proof.} The uniqueness of $\lambda_t$ follows from the fact
 that $F_t$ is a submersion as it was explained few lines above.
 Let us proof the existence. To do that we use the
 coordinatization of $\mathcal V$ introduced in the proof of
 Proposition~I.1, in particular, the family of local trivializations
$R_t:O_{q(t)}\times U\to\pi^{-1}(O_{q(t)})$. Assume that
$v(t)=R_t(q(t),u(t)),\ t_0\le t\le t_1$, where $v(\cdot)$ is the
referenced admissible pair from (8).

Given $\tau\in[t_0,t_1],\ \hat q\in O_{q(\tau)}$ let $t\mapsto
Q_\tau^t(\hat q)$ be the solution of the differential equation
$\dot q=R_t(q,u(t))$ which satisfies the condition
$Q^\tau_\tau(\hat q)=\hat q$. In particular,
$Q_\tau^t(q(\tau))=q(t)$. Then $Q_\tau^t$ is a diffeomorphism of a
neighborhood of $q(\tau)$ on a neighborhood of $q(t)$. We define a
Banach submanifold $\mathcal V_\tau$ of the Banach manifold
$\mathcal V$ in the following way:
$$
\mathcal V_\tau=\{v'\in\mathcal
V:\pi(v'(t))=Q_\tau^t(\pi(v'(\tau))),\ \tau\le t\le t_1\}.
$$
It is easy to see that $F_{t_1}\Bigr|_{\mathcal
V_\tau}=Q^{t_1}_\tau\circ F_\tau\Bigr|_{\mathcal V_\tau}$ and
$J_\tau^{t_1}\Bigr|_{\mathcal V_\tau}=a_\tau\circ F_\tau$, where
$a(\hat q)=\int\limits_\tau^t\varphi\left(\Phi_t(Q_\tau^t(\hat
q),u(t))\right)\,dt$. On the other hand, the set
$\{v'\in\mathcal V : v'|_{[t_0,\tau]}\in\mathcal
V_\tau\bigr|_{[t_0,\tau]}\}$ is a neighborhood of $v$ in $\mathcal
V$. The restriction of (8) to $\mathcal V_\tau$ gives:
$$
\lambda_{t_1}D_v\left(Q^{t_1}_\tau\circ F_\tau\right)-
\lambda_{t_0}D_vF_{t_0}=d_vJ^\tau_{t_0}+d_v\left(a_\tau\circ F_\tau\right).
$$
Now we apply the chain rule for the differentiation and obtain:
$$
\lambda_\tau D_vF_\tau-\lambda_{t_0}D_vF_{t_0}=d_vJ^\tau_{t_0},
$$
where
$\lambda_\tau=\lambda_{t_1}D_{q(\tau)}Q_\tau^{t_1}-d_{q(\tau)}a_\tau$.
$\quad \square$

\medskip\noindent {\bf Definition.} A Lipschitzian curve
$t\mapsto\lambda_t,\ t_0\le t\le t_1,$ is called a {\it normal
extremal} of the given optimal control problem if there exists an
admissible pair $v\in\mathcal V$ such that equality (9) holds. The
projection $q(t)=\pi(\lambda_t)$ of a normal extremal is called a
(normal) {\it extremal path} or a (normal) {\it extremal
trajectory}.

\medskip According to Proposition~I.2, normal Lagrange multipliers
are just points of normal extremals. A good thing about normal
extremals is that they satisfy a nice differential equation which
links optimal control theory with a beautiful and powerful
mathematics and, in many cases, allows to explicitly characterize
all extremal paths.

\section{Hamiltonian system}
Here we derive equations which characterize normal extremals; we
start from coordinate calculations. Given $\tau\in[t_0,t_1]$,
fix a coordinate neighborhood $\mathcal O$ in $M$ centered at
$q(\tau)$, and focus on the piece of the extremal path $q(\cdot)$
which contains $q(\cdot)$ and is completely contained in $\mathcal
O$. Identity (9) can be rewritten in the form
$$
\lambda_{t}D_vF_{t}-\lambda_{\tau}D_vF_{\tau}=d_vJ^t_{\tau}, \eqno
 (10)
$$
where $q(t)$ belongs to the piece of $q(\cdot)$ under
consideration. Fixing coordinates and a local trivialization of
$V$ we (locally) identify our optimal control problem with
a problem (1), (2) in $\mathbb R^n$. We have
$T^*\mathbb R^n\cong\mathbb R^n\times\mathbb
R^n=\{(p,q):p,q\in\mathbb R^n\}$, where $T_q^*\mathbb R^n=\mathbb
R^n\times\{q\}$. Then $\lambda_t=\{p(t),q(t)\}$ and
$\lambda_tD_vF_t\cdot=\langle p(t),D_vF_t\cdot\rangle=
D_v\langle p(t),F_t\rangle$.

Admissible pairs of (2) are parametrized by $\hat q=F_\tau(v'),\
v'\in\mathcal V$, and control functions $u'(\cdot)$; the pairs
have the form: $v'=\left(u'(\cdot),q'(\cdot;\hat
q,u'(\cdot))\right)$, where
$\frac\partial{\partial t}q'(t;\hat q,u'(\cdot))=f\left(q'(t;\hat
q,u'(\cdot)),u'(t)\right)$ for all available $t$
and $q'(\tau;\hat q,u(\cdot))=\hat q$. Then
$F_t(v')=q'(t;\hat q,u'(\cdot))$.

Now we differentiate identity (10) with respect to $t$:
$\frac\partial{\partial t}D_v\langle
p(t),F_t\rangle=\frac\partial{\partial t}d_vJ^t_\tau$ and change
the order of the differentiation
$D_v\frac\partial{\partial t}\langle p(t),F_t\rangle=
d_v\frac\partial{\partial t}J^t_\tau$. We compute the derivatives
with respect to $t$ at $t=\tau$:
$$
\frac\partial{\partial t}\langle p(t),F_t\rangle\bigr|_{t=\tau}=
\langle\dot p(\tau),\hat q\rangle+\langle p(\tau),f(\hat
q,u'(\tau)\rangle, \quad
\frac\partial{\partial t}J^t_\tau\bigr|_{t=\tau}=\varphi(\hat
q,u'(\tau)).
$$
Now we have to differentiate with respect to
$v'(\cdot)=(u'(\cdot),q'(\cdot))$. We however see that the
quantities to differentiate depend only on the values of
$u'(\cdot)$ and $q'(\cdot)$ at $\tau$, i.e. on the
finite-dimensional vector $(u'(\tau),\hat q)$. We derive:
$$
\dot p(\tau)+\frac\partial{\partial q}\langle
p(\tau),f(q(\tau),u(\tau))\rangle=\frac{\partial\varphi}{\partial
q}(q(t),u(t)),
$$
$$
 \frac\partial{\partial u}\langle
p(\tau),f(q(\tau),u(\tau))\rangle=\frac{\partial\varphi}{\partial
u}(q(\tau),u(\tau)),
$$
where $v(\cdot)=(q(\cdot),u(\cdot))$.

Of course, we can change $\tau$ and perform the differentiation at
any available moment $t$. Finally, we obtain that (10) is
equivalent to the identities
$$
\dot p(t)+\frac\partial{\partial q}\left(\langle p(t),f(q(t),u(t))\rangle-
\varphi(q(t),u(t))\right)=0,
$$
$$
\frac\partial{\partial u}\left(\langle p(t),f(q(t),u(t))\rangle-
\varphi(q(t),u(t))\right)=0,
$$
which can be completed by the equation
$\dot q=f(q(t),u(t))$. We introduce a function
$h(p,q,u)=\langle p,f(q,u)\rangle-\varphi(q,u)$ which is called
the {\it Hamiltonian} of the optimal control problem (1), (2).
This function permits us to present the obtained relations in a nice
Hamiltonian form:
$$
\left\{\begin{aligned}\dot p & =-\frac{\partial h}{\partial
q}(p,q,u)\\
\dot q & =\frac{\partial h}{\partial p}(p,q,u)\\
\end{aligned}\right., \quad \frac{\partial h}{\partial
u}(p,q,u)=0. \eqno (11)
$$

A more important fact is that system (11) has an intrinsic
coordinate free interpretation. Recall that in the triple
$(p,q,u)$ neither $p$ nor $u$ has an intrinsic meaning; the pair
$(p,q)$ represents $\lambda\in T^*M$ while the pair $(q,u)$
represents $v\in V$. First we consider an intermediate case
$V=M\times U$ (when $u$ is separated from $q$ but coordinates in
$M$ are not fixed) and then turn to the completely intrinsic
setting.

If $V=M\times U$, then $f:M\times U\to TM$ and $f(q,u)\in T_qM$.
The Hamiltonian of the optimal control problem is a function
$h:T^*M\times U\to\mathbb R$ defined by the formula
$h(\lambda,u)=\lambda(f(q,u))-\varphi(q,u)$, $\forall\lambda\in
T^*_qM,\ q\in M,\ u\in U$. For any $u\in U$ we obtain a function
$h_u\stackrel{def}{=}h(\cdot,u)$ on $T^*M$. The cotangent bundle
$T^*M$ possesses a canonical symplectic structure which provides a
standard way to associate a {\it Hamiltonian vector field} to any
smooth function on $T^*M$. We'll recall this procedure.

Let $\pi: T^*M\to M$ be the projection, $\pi(T^*_qM)=\{q\}$. The
{\it Liouville} (or {\it tautological}) differential 1-form
$\varsigma$ on $T^*M$ is defined as follows. Let \linebreak
$\varsigma_\lambda:T_\lambda(T^*M)\to \mathbb R$ be the value of
$\varsigma$ at $\lambda\in T^*M$, then
$\varsigma_\lambda=\lambda\circ\pi_*$, the composition of
$\pi_*:T_\lambda(T^*M)\to T_{\pi(\lambda)}M$ and the cotangent vector
$\lambda:T_{\pi(\lambda)}M\to\mathbb R$. The coordinate presentation
of the Liouville form is:
$\varsigma_{(p,q)}=\langle p,dq\rangle=\sum\limits_{i=1}^np^idq^i$,
where $p=(p^1,\ldots,p^n)$, $q=(q^1,\ldots,q^n)$. The {\it canonical
symplectic structure} on $T^*M$ is the differential 2-form
$\sigma=d\varsigma$; its coordinate representation is:
$\sigma=\sum\limits_{i=1}^ndp^i\wedge dq^i$. The Hamiltonian
vector field associated to a smooth function $a:T^*M\to\mathbb R$
is a unique vector field $\vec a$ on $T^*M$ which satisfies the
equation $\sigma(\cdot,\vec a)=da$. The coordinate representation of
this field is:
$\vec a=\sum\limits_{i=1}^n\left(\frac{\partial a}{\partial p_i}
\frac\partial{\partial q_i}-\frac{\partial a}{\partial q_i}
\frac\partial{\partial p_i}\right)$. Equations (11) can be rewritten in the form:
$$
\dot\lambda=\vec h_u(\lambda), \quad \frac{\partial h}{\partial u}(\lambda,u)=0.
\eqno (12)
$$
Now let $V$ be an arbitrary locally trivial bundle over $M$. Consider
the Cartesian product of two bundles:
$$
T^*M\times_MV=\{(\lambda,v):v\in V_q,\lambda\in T^*_qM,q\in M\}
$$
that is a bundle over $M$ whose fibers are Cartesian products of the correspondent
fibers of $V$ and $T^*M$. Hamiltonian of the optimal control problem takes the form
$h(\lambda,v)=\lambda(f(v))-\varphi(v)$; this is a well-defined smooth function on
$T^*M\times_MU$. Let $\mathfrak p:T^*M\times_MV\to T^*M$ be the projection on the first factor,
$\mathfrak p:(\lambda,v)\mapsto\lambda$. Equations (11) (or (12)) can be rewritten
in the completely
intrinsic form as follows: $(\mathfrak p^*\sigma)_v(\cdot,\dot\lambda)=dh$.
One may check this fact
in any coordinates; we leave this simple calculation to the reader.

Of course, by fixing a local trivialization of $V$, we turn the last relation back
into a more convinient to study equation (12). A domain $\mathcal D$
in $T^*M$  is called regular for the Hamiltonian $h$ if for any $\lambda\in\mathcal D$
there exists a unique
solution $u=\bar u(\lambda)$ of the equation $\frac{\partial h}{\partial u}(\lambda,u)=0$,
where $\bar u(\lambda)$ is smooth with respect to $\lambda$. In particular, if $U$ is an affine
space and the functions $u\mapsto h(\lambda,u)$ are strongly concave (convex) and
bounded from above (below) for $\lambda\in\mathcal D$,
then $\mathcal D$ is regular
and $\bar u(\lambda)$ is defined by the relation
$$
h(\lambda,\bar u(\lambda))=\max\limits_{u\in U}h(\lambda,u) \quad
\left(h(\lambda,\bar u(\lambda))=\min\limits_{u\in U}h(\lambda,u)\right).
$$
In the regular domain, we set $H(\lambda)=h(\lambda,\bar u(\lambda))$, where
$\frac{\partial h}{\partial u}(\lambda,\bar u(\lambda))=0$. It is easy to see that equations
(12) are equivalent to one Hamiltonian system $\dot\lambda=\vec H(\lambda)$.
Indeed, the equality $d_{(\lambda,\bar u(\lambda))}h=d_\lambda h_{\bar u(\lambda)}+
\frac{\partial h_{\bar u(\lambda)}}{\partial u}du=d_\lambda h_{\bar u(\lambda)}$ immediately
implies that $\vec H(\lambda)=\vec h_{\bar u(\lambda)}(\lambda)$.

\section{Second order information}

We come back to the general setting of Section 2 and try to go beyond the Lagrange
multipliers rule. Take a pair $(\ell,w)$ which satisfies equation (6). We call such
pairs (normal) {\it Lagrangian points}. Let $\Phi(w)=z$. If $w$ is a regular point of $\Phi$,
then $\Phi^{-1}(z)\cap O_w$ is a smooth codimension $\dim N$ submanifold of $\mathcal W$,
for some neighborhood $O_w$ of $w$. In this case $w$ is a critical point of
$J\bigl|_{\Phi^{-1}(z)\cap O_w}$. We are going to compute the Hessian of
$J\bigl|_{\Phi^{-1}(z)}$ at $w$ without resolving the constraints $\Phi(w)=z$.
The formula we obtain makes sense without the regularity assumptions as well.

Let $s\mapsto\gamma(s)$ be a smooth curve in $\Phi^{-1}(z)$ such that $\gamma(0)=w$.
Differentiation of the identity $\Phi(\gamma(s))=z$ gives:
$$
D_w\Phi\dot\gamma=0, \quad D^2_w\Phi(\dot\gamma,\dot\gamma)+D_w\Phi\ddot\gamma=0,
$$
where $\dot\gamma$ and $\ddot\gamma$ are the first and the second derivatives of
$\gamma$ at $s=0$. We also have:
$$
\frac{d^2}{ds^2}J(\gamma(s))|_{s=0}=D^2_wJ(\dot\gamma,\dot\gamma)+D_wJ\ddot\gamma
\stackrel{\mathrm{eq.}(6)}{=}
$$
$$
D^2_wJ(\dot\gamma,\dot\gamma)+\ell D_w\Phi\ddot\gamma=D^2_wJ(\dot\gamma,\dot\gamma)-
\ell D^2_w\Phi(\dot\gamma,\dot\gamma).
$$
Finally,
$$
\mathrm{Hess}_w(J\bigl|_{\Phi^{-1}(z)})=(D^2_wJ-\ell D^2_w\Phi)\bigl|_{\ker D_w\Phi}.
\eqno (13) $$

\begin{prop} If quadratic form (13) is positive (negative) definite, then $w$ is
a strict local minimizer (maximizer) of $J\bigl|_{\Phi^{-1}(z)}$.
\end{prop}

If $w$ is a regular point of $\Phi$, then the proposition is obvious but one can
check that it remains valid without the regularity assumption.
On the other hand, without the regularity assumption, local
minimality does not imply nonnegativity of form (13). What local
minimality (maximality) certainly implies is nonnegativity
(nonpositivity) of form (13) on a finite codimension subspace of
$\ker D_w\Phi$ (see \cite[Ch.~20]{as} and references there).

\medskip\noindent {\bf Definition.} A Lagrangian point $(\ell,w)$
is called {\it sharp} if quadratic form (13) is nonnegative or
nonpositive on a finite codimension subspace of $\ker D_w\Phi$.

\medskip
Only sharp Lagrangian points are counted in the conditional
extremal problems under consideration. Let $Q$ be a real quadratic
form defined on a linear space $E$. Recall that the {\it negative
inertia index} (or the {\it Morse index}) $\mathrm{ind}Q$ is the
maximal possible dimension of a subspace in $E$ such that the
restriction of $Q$ to the subspace is a negative form. The
{\it positive inertia index} of $Q$ is the Morse index of $-Q$.
Each of these indices is a nonnegative integer or $+\infty$. A
Lagrangian point $(\ell,w)$ is sharp if the negative or positive
inertia index of form (13) is finite.

In the optimal
control problems, $\mathcal W$ is a huge infinite dimensional manifold while $N$
usually has a modest dimension. It is much simpler to characterize Lagrange
multipliers in $T^*N$ (see the previous section)
than to work directly with $J\bigl|_{\Phi^{-1}(z)}$. Fortunately, the
information on the sign and, more generally, on the inertia indices of the infinite
dimensional quadratic form (13) can also be
extracted from the Lagrange multipliers or, more precisely, from the so called
{\it $\mathcal L$-derivative} that can be treated as a dual to the form (13) object.

$\mathcal L$-derivative concerns the linearization of equation (6) at a given Lagrangian
point. In order to linearize the equation we have to present its left- and right-hand
sides as smooth mappings of some manifolds. No problem with the right-hand side:
$w\mapsto d_wJ$ is a smooth mapping from $\mathcal W$ to $T^*\mathcal W$. The variables
$(\ell,w)$ of the left-hand side live in the manifold
$$
\Phi^*T^*N=\{(\ell,w): \ell\in T^*_{\Phi(w)},\ w\in\mathcal W\}\subset
T^*N\times\mathcal W.
$$
Note that $\Phi^*T^*N$ is a locally trivial bundle over $\mathcal W$ with the
projector $\pi:(\ell,w)\mapsto w$; this is nothing else but the {\it induced
bundle} from $T^*N$ by the mapping $\Phi$. We treat equation (6) as the equality
of values of two mappings from $\Phi^*T^*N$ to $T^*\mathcal W$. Let us rewrite
this equation in local coordinates.

So let $N=\mathbb R^m$ and $\mathcal W$ be a Banach space. Then
$T^*N=\mathbb R^{m*}\times\mathbb R^m$ (where $T_zN=\mathbb
R^{m*}\times\{z\}$), $T^*\mathcal W=\mathcal W^*\times\mathcal W$,
$\Phi^*T^*N=\mathbb R^{m*}\times\mathbb R^m\times\mathcal W$.
Surely, $\mathbb R^{m*}\cong\mathbb R^m$ but in the forthcoming
calculations it is convenient to treat the first factor in the
product $\mathbb R^{m*}\times\mathbb R^m$ as the space of linear
forms on the second factor. We have:
$\ell=(\zeta,z)\in\mathbb R^{m*}\times\mathbb R^m$ and equation
(6) takes the form
$$
\zeta\frac{d\Phi}{dw}=\frac{dJ}{dw}, \quad \Phi(w)=z. \eqno (14)
$$

Linearization of system (14) at the point $(\zeta,z,w)$ reads:
$$
\zeta'\frac{d\Phi}{dw}+\zeta\frac{d^2\Phi}{dw^2}(w',\cdot)=
\frac{d^2J}{dw^2}(w',\cdot), \quad \frac{d\Phi}{dw}w'=z'. \eqno
(15)
$$
We set
$$
\mathcal L^0_{(\ell,w)}(\bar\Phi)=\{\ell'=(\zeta',z')\in
T_\ell (T^*N) :\exists w'\in\mathcal W\ \mathrm{s.t.}\
(\zeta',z',w')\ \mathrm{satisfies}\ (15)\}.
$$
Note that subspace $\mathcal L^0_{(\ell,w)}(\bar\Phi)\subset
T_{\ell}(T^*N)$ does not depend on the choice of local
coordinates. Indeed, to construct this subspace we take all
$(\ell',w')\in T_{(\ell,w)}(\Phi^*T^*N)$ which satisfy the
linearized equation (6) and then apply the projection
$(\ell',w')\mapsto\ell'$.

Recall that $T_\ell(T^*N)$ is a symplectic space endowed with the
canonical symplectic form $\sigma_\ell$ (cf. Sec.~4). A subspace
$S\subset T_\ell(T^*N)$ is {\it isotropic} if
$\sigma_{\ell}|_S=0$. Isotropic subspaces of maximal possible
dimension $m=\frac 12\dim T_\ell(T^*N)$ are called {\it Lagrangian
subspaces}.

\begin{prop} $\mathcal L^0_{(\ell,w)}(\bar\Phi)$ is an isotropic
subspace of $T_\ell(T^*N)$. If $\dim\mathcal W<\infty$, then
$\mathcal L^0_{(\ell,w)}(\bar\Phi)$ is a Lagrangian subspace.
\end{prop}
{\bf Proof.} First we'll prove the isotropy of $\mathcal
L^0_{(\ell,w)}(\bar\Phi)$. Let $(\zeta',z'),(\zeta'',z'')\in
T_\ell(T^*N)$. We have
$\sigma_\ell((\zeta',z'),(\zeta'',z''))=\zeta'z''-\zeta''z'$;
here the symbol $\zeta z$ denotes the result of the application
of the linear form $\zeta\in\mathbb R^{m*}$ to the vector
$z\in\mathbb R^n$ or, in the matrix terminology, the product of
the row $\zeta$ and the column $z$. Assume that $(\zeta',z',w')$
and $(\zeta'',z'',w'')$ satisfy equations (15); then
$$
\zeta'z''=\zeta'\frac{d\Phi}{dw}w''=\frac{d^2J}{dw^2}(w',w'')-
\zeta\frac{d^2\Phi}{dw^2}(w',w''). \eqno (16)
$$
The right-hand side of (16) is symmetric with respect to $w'$ and
$w''$ due to the symmetry of second derivatives. Hence
$\zeta'z''=\zeta''z'$. In other words,
$\sigma_{\ell}((\zeta',z'),(\zeta'',z''))=0$. So $\mathcal
L^0_{(\ell,w)}(\bar\Phi)$ is isotropic and, in particular,
\linebreak $\dim\left(\mathcal L^0_{(\ell,w)}(\bar\Phi)\right)\le m$.

Now show that the last inequality becomes the equality as soon as
$\mathcal W$ is finite dimensional. Set
$Q=\frac{d^2J}{dw^2}-\zeta\frac{d^2\Phi}{dw^2}$ and consider the
diagram:
$$
\zeta'\frac{d\Phi}{dw}-Q(w',\cdot)\ \stackrel{left}{\longleftarrow}
\ (\zeta',w')\ \stackrel{right}{\longrightarrow}\
\left(\zeta',\frac{d\Phi}{dw}w'\right).
$$
Then $\mathcal L^0_{(\ell,w)}(\bar\Phi)=right(\ker(left))$.
Passing to a factor space if necessary we may assume that
$\ker(left)\cap\ker(right)=0$; this means that:
$$
\frac{d\Phi}{dw}w'\quad \&\quad Q(w',\cdot)=0\quad \Rightarrow\quad w'=0.
\eqno (17)
$$
Under this assumption, $\dim \mathcal L^0_{(\ell,w)}(\bar\Phi)=
\dim\ker(left)$. On the other hand, relations (17) imply that
the mapping $left:\mathbb R^{m*}\times\mathcal W\to\mathcal W^*$
is surjective. Indeed, if, on the contrary, the map $left$ is not
surjective then there exists a nonzero vector $v\in(\mathcal
W^*)^*=\mathcal W$ which annihilates the image of $left$; in other
words, $\zeta'\frac{d\Phi}{dw}v-Q(w',v)=0,\ \forall\zeta',w'$.
Hence $\frac{d\Phi}{dw}v=0\ \&\ Q(v,\cdot)=0$ that contradicts
(17). It follows that $\dim\mathcal L^0_{(\ell,w)}(\bar\Phi)=
\dim(\mathbb R^{m*}\times\mathcal W)-\dim\mathcal W^*=m.
\quad\square$

For infinite dimensional $\mathcal W$, the space
$\mathcal L^0_{(\ell,w)}(\bar\Phi)$ may have dimension smaller
than $m$ due to an ill-posedness of equations (15); to guarantee
dimension $m$ one needs certain coercivity of the form
$\zeta\frac{d^2\Phi}{dw^2}$. I am not going to discuss here what
kind of coercivity is sufficient, it can be easily reconstructed
from the proof of Proposition~I.4 (see also \cite{a}). Anyway,
independently on any coercivity one can take a finite dimensional
approximation of the original problem and obtain a Lagrangian
subspace $\mathcal L^0_{(\ell,w)}(\bar\Phi)$ guaranteed by
Proposition~I.4. What happens with these subspaces when the
approximation becomes better and better, do they have a
well-defined limit (which would be unavoidably Lagrangian)? A
remarkable fact is that such a limit does exist for any sharp
Lagrangian point. It contains $\mathcal L^0_{(\ell,w)}(\bar\Phi)$
and is called the $\mathcal L$-derivative of $\bar\Phi$ at
$(\ell,w)$. To formulate this result we need some basic
terminology from set theoretic topology.

A partially ordered set $(\mathfrak A,\prec)$ is a {\it directed
set} if $\forall\alpha_1,\alpha_2\in\mathfrak A$
$\exists\beta\in\mathfrak A$ such that $\alpha_1\prec\beta$ and
$\alpha_2\prec\beta$. A family $\{x_\alpha\}_{\alpha\in\mathfrak
A}$ of points of a topological space $\mathcal X$ indexed by the
elements of $\mathfrak A$ is a generalized sequence in $\mathcal
X$. A point $x\in\mathcal X$ is the limit of the generalized
sequence $\{x_\alpha\}_{\alpha\in\mathfrak A}$ if for any
neighborhood $\mathcal O_x$ of $x$ in $\mathcal X$
$\exists\alpha\in\mathfrak A$ such that $x_\beta\in\mathcal O_x,\
\forall\beta\succ\alpha$; in this case we write $x=\lim\limits_{\mathfrak
A} x_\alpha$.

Let $\mathfrak w$ be a finite dimensional submanifold of $\mathcal
W$ and $w\in\mathfrak w$. If $(\ell,w)$ is a Lagrangian point for
$\bar\Phi=(J,\Phi)$, then it is a Lagrangian point for
$\bar\Phi|_{\mathfrak w}$. A straightforward calculation shows
that the Lagrangian subspace
$\mathcal L^0_{(\ell,w)}(\bar\Phi|_{\mathfrak w})$ depends on the
tangent space $W=T_w\mathfrak w$ rather than on $\mathfrak w$,
i.e. $\mathcal L^0_{(\ell,w)}(\bar\Phi|_{\mathfrak w})=
\mathcal L^0_{(\ell,w)}(\bar\Phi|_{\mathfrak w'})$ as soon as
$T_w\mathfrak w=T_w\mathfrak w'=W$. We denote $\Lambda_W
=\mathcal L^0_{(\ell,w)}(\bar\Phi|_{\mathfrak w})$. Recall that
$\Lambda_W$ is an $m$-dimensional subspace of the
$2m$-dimensional space $T_\ell(T^*N)$, i.e. $\Lambda_W$ is a
point of the Grassmann manifold of all $m$-dimensional
subspaces in $T_\ell(T^*N)$.

Finally, we denote by $\mathfrak W$ the set of all finite
dimensional subspaces of $T_w\mathcal W$ partially ordered by the
inclusion ``$\subset$". Obviously, $(\mathfrak W,\subset)$ is a
directed set and $\{\Lambda_W\}_{W\in\mathfrak W}$ is a generalized
sequence indexed by the elements of this directed set. It is easy
to check that there exists $W_0\in\mathfrak W$ such that
$\Lambda_W\supset\mathcal L^0_{(\ell,w)}(\bar\Phi),\ \forall
W\supset W_0$. In particular, if $\mathcal
L^0_{(\ell,w)}(\bar\Phi)$ is $m$-dimensional, then
$\Lambda_{W_0}=\mathcal L^0_{(\ell,w)}(\bar\Phi),\ \forall
W\supset W_0$, the sequence $\Lambda_W$ is stabilizing and
$\mathcal L^0_{(\ell,w)}(\bar\Phi)=\lim\limits_{\mathfrak W}\Lambda_W$. In
general, the sequence $\Lambda_W$ is not stabilizing, nevertheless
the following important result is valid.

\begin{theorem} If $(\ell,w)$ is a sharp Lagrangian point, then
there exists \linebreak $\mathcal L_{(\ell,w)}(\bar\Phi)=
\lim\limits_{\mathfrak W}\Lambda_W$.
\end{theorem}

We omit the proof of the theorem, you can find this proof in paper
\cite{a} with some other results which allow to efficiently
compute $\lim\limits_{\mathfrak W}\Lambda_W$. Lagrangian subspace
$\mathcal L_{(\ell,w)}(\bar\Phi)=\lim\limits_{\mathfrak
W}\Lambda_W$ is called the {\it $\mathcal L$-derivative} of
$\bar\Phi=(J,\Phi)$ at the Lagrangian point $(\ell,w)$.

Obviously, $\mathcal L_{(\ell,w)}(\bar\Phi)\supset\mathcal
L^0_{(\ell,w)}(\bar\Phi)$. One should think on
$\mathcal L_{(\ell,w)}(\bar\Phi)$ as on a completion of
$\mathcal L^0_{(\ell,w)}(\bar\Phi)$ by means of a kind of
weak solutions to system (15) which could be missed due to the
ill-posedness of the system.

Now we should explain the connection between
$\mathcal L_{(\ell,w)}(\bar\Phi)$ and \linebreak
$\mathrm{Hess}_w(J\bigr|_{\Phi^{-1}(z)})$. We start from the
following simple observation:

\begin{lemma} Assume that $\dim\mathcal W<\infty$, $w$ is
a regular point of $\Phi$ and \linebreak $\ker D_w\Phi\cap\ker(D^2_wJ-\ell
D^2_w\Phi)=0$. Then
$$
\ker\mathrm{Hess}_w(J\bigr|_{\Phi^{-1}(z)})=0\quad
\Leftrightarrow \quad \mathcal L_{(\ell,w)}(\bar\Phi)\cap
T_\ell(T^*_zN)=0,
$$
i.e. quadratic form $\mathrm{Hess}_w(J\bigr|_{\Phi^{-1}(z)})$
is nondegenerate if and only if the subspace $\mathcal
L_{(\ell,w)}(\bar\Phi)$ is transversal to the fiber $T^*_zN$.
\end{lemma}
{\bf Proof.} We make computations in coordinates. First,
$T_\ell(T^*_zN)=\{(\zeta',0):\zeta'\in\mathbb R^{n*}\}$; then,
according to equations (15), $(\zeta',0)\in\mathcal
L_{(\ell,w)}(\bar\Phi)$ if and only if there exists $w\in\mathcal
W$ such that
$$
\frac{d\Phi}{dw} w'=0,\quad
\frac{d^2J}{dw^2}(w',\cdot)-\ell\frac{d^2\Phi}{dw^2}(w',\cdot)=
\zeta'\frac{d\Phi}{dw}. \eqno (18)
$$
Regularity of $w$ implies that $\zeta'\frac{d\Phi}{dw}\ne 0$ and
hence $w'\ne 0$ as soon as $\zeta'\ne 0$. Equalities (18) imply:
$\frac{d^2J}{dw^2}(w',v)-\ell\frac{d^2\Phi}{dw^2}(w',v)=0,
\quad\forall v\in\ker\frac{d\Phi}{dw}$, i.e.
$w'\in\ker\mathrm{Hess}_w(J\bigr|_{\Phi^{-1}(z)})$. Moreover, our
implications are invertible: we could start from a nonzero vector
$w'\in\ker\mathrm{Hess}_w(J\bigr|_{\Phi^{-1}(z)})$
and arrive to a
nonzero vector $(\zeta',0)\in\mathcal L_{(\ell,w)}(\bar\Phi). \quad
\square $

\medskip\noindent{\bf Remark.} Condition
$\ker D_w\Phi\cap\ker(D^2_wJ-\ell D^2_w\Phi)=0$ from Lemma~I.1 is not
heavy. Indeed, a pair $(J,\Phi)$ satisfies this condition at all
its Lagrangian points if and only if 0 is a regular value of the
mapping $(\zeta,w)\mapsto\zeta\frac{d\Phi}{dw}-\frac{dJ}{dw}$.
Standard Transversality Theorem implies that this is true for
generic pair $(J,\Phi)$.

\section{Maslov index}

Lemma~I.1 is a starting point for a far going theory which allows to
effectively compute the Morse index of the Hessians in terms of the
$\mathcal L$-derivatives.

How to do it?
Normally, extremal problems depend on some parameters. Actually,
$z\in N$ is such a parameter and there could be other ones,
which we do not explicitly add to the constraints. In the optimal
control problems a natural parameter is the time interval
$t_1-t_0$. Anyway, assume that we have a continuous family of the
problems and their sharp Lagrangian points: $\ell_\tau
D_{w_\tau}\Phi_\tau=d_{w_\tau}J_\tau,\ \tau_0\le\tau\le\tau_1$;
let $\Lambda(\tau)=\mathcal L_{(\ell_\tau,w_\tau)}(\bar\Phi_\tau)$.
Our goal is to compute the difference
$\mathrm{ind}\,\mathrm{Hess}_{w_{\tau_1}}
(J_{\tau_1}\bigr|_{\Phi_{\tau_1}^{-1}(z_{\tau_1})})-
\mathrm{ind}\,\mathrm{Hess}_{w_{\tau_0}}
(J_{\tau_0}\bigr|_{\Phi_{\tau_0}^{-1}(z_{\tau_0})})
$ in terms of the family of Lagrangian subspaces $\Lambda(\tau)$;
that is to get a tool to follow the evolution of the Morse index
under a continuous change of the parameters. This is indeed very
useful since for some special values of the parameters the index
could be known a'priori. It concerns, in particular, optimal
control problems with the parameter $\tau=t_1-t_0$. If $t_1-t_0$
is very small then sharpness of the Lagrangian point almost
automatically implies the positivity or negativity of the Hessian.

First we discuss the finite-dimensional case: Theorem~I.1 indicates
that finite-dimensional approximations may already contain all
essential information. Let $Q_\tau$ be a continuous family of
quadratic forms defined on a finite-dimensional vector space. If
$\ker Q_\tau=0,\ \tau_0\le\tau\le\tau_1$, then
$\mathrm{ind}Q_\tau$ is constant on the segment $[\tau_0,\tau_1]$.
This is why Lemma~I.1 opens the way to follow evolution of the index
in terms of the $\mathcal L$-derivative: it locates values of the
parameter where the index may change. Actually, $\mathcal
L$-derivative allows to evaluate this change as well; the
increment of $\mathrm{ind}Q_\tau$ is computed via so called Maslov
index of a family of Lagrangian subspaces. In order to define this
index we have to recall some elementary facts about symplectic
spaces.

Let $\Sigma,\sigma$ be a symplectic space, i.e. $\Sigma$ is a
$2n$-dimensional vector space and $\sigma$ be a nondegenerate
anti-symmetric bilinear form on $\Sigma$. The skew-orthogonal
complement to the subspace $\Gamma\subset\Sigma$ is the subspace
$\Gamma^\angle=\{x\in\Sigma : \sigma(x,\Gamma)=0\}$. The
nondegeneracy of $\sigma$ implies that
$\dim\Gamma^\angle=2n-\dim\Gamma$. A subspace $\Gamma$ is
isotropic if and only if $\Gamma^\angle\supset\Gamma$; it is
Lagrangian if and only if $\Gamma^\angle=\Gamma$.

Let $\Pi=span\{e_1,\ldots,e_n\}$ be a lagrangian subspace
of $\Sigma$. Then there exist vectors $f_1,\ldots,f_n\in\Sigma$
such that $\sigma(e_i,f_j)=\delta_{ij}$, where $\delta_{ij}$ is
the Kronecker symbol. We show this using induction with respect
to $n$. Skew-orthogonal complement to the space
$span\{e_1,\ldots,e_{n-1}\}$ contains an element $f$ which is
not skew-orthogonal to $e_n$; we set $f_n=\frac
1{\sigma(e_n,f)}f$. We have
$$
span\{e_n,f_n\}\cap span\{e_n,f_n\}^\angle=0
$$
and the restriction of $\sigma$ to
$span\{e_n,f_n\}^\angle$ is a nondegenerate bilinear
form. Hence $span\{e_n,f_n\}^\angle$ is a
$2(n-1)$-dimensional symplectic space with a Lagrangian subspace
$span\{e_1,\ldots,e_{n-1}\}$. According to the induction
assumption, there exist $f_1,\ldots,f_{n-1}$ such that
$\sigma(e_i,f_j)=\delta_{ij}$ and we are done.

Vectors $e_1,\ldots,e_n,f_1,\ldots,f_n$ form a basis of $\Sigma$;
in particular,\linebreak $\Delta=span\{f_1,\ldots,f_n\}$ is a
transversal to $\Pi$ Lagrangian subspace,
$\Sigma=\Pi\oplus\Delta$. If
$x_i=\sum\limits_{j=1}^n(\zeta_i^je_j+z_i^jf_j),\ i=1,2,$ and
$\zeta_i=(\zeta^1_i,\ldots,\zeta^n_i)$,
$z_i=(z_i^1,\ldots,z_i^n)^\top$, then
$\sigma(x_1,x_2)=\zeta_1z_2-\zeta_2z_1$. The coordinates $\zeta,z$
identify $\Sigma$ with $\mathbb R^{n*}\times\mathbb R^n$; any
transversal to $\Delta$ $n$-dimensional subspace
$\Lambda\subset\Sigma$ has the following presentation in these
coordinates:
$$
\Lambda=\{z^\top,S_\Lambda z):z\in\mathbb R^n\},
$$
where $S_\Lambda$ is an $n\times n$-matrix. The subspace $\Lambda$
is Lagrangian if and only if $S^*_\Lambda=S_\Lambda$. We have:
$$
\Lambda\cap\Pi=\{(z^\top,0):z\in\ker S_\Lambda\},
$$
the subspace $\Lambda$ is transversal to $\Pi$ if and only if
$S_\Lambda$ is nondegenerate.

That's time to introduce some notations. Let $L(\Sigma)$ be the
set of all Lagrangian subspaces, a closed subset of the
Grassmannian $G_n(\Sigma)$ of $n$-dimensional subspaces in
$\Sigma$. We set
$$
\Delta^\pitchfork=\{\Lambda\in L(\Sigma):\Lambda\cap\Delta=0\},
$$
an open subset of $L(\Sigma)$. The mapping $\Lambda\mapsto
S_\Lambda$ gives a regular parametrization of $\Delta^\pitchfork$
by the $n(n+1)/2$-dimensional space of symmetric $n\times
n$-matrices. Moreover, above calculations show that
$L(\Sigma)=\bigcup\limits_{\Delta\in L(\Sigma)}\Delta^\pitchfork$.
Hence $L(\Sigma)$ is a $n(n+1)/2$-dimensional submanifold of the
Grassmannian $G_n(\Sigma)$ covered by coordinate charts
$\Delta^\pitchfork$. The manifold $L(\Sigma)$ is called Lagrange
Grassmannian associated to the symplectic space $\Sigma$. It is
not hard to show that any coordinate chart $\Delta^\pitchfork$ is
everywhere dense in $L(\Sigma)$;
our calculations give also a local parametrization of its
complement.

Given $\Pi\in L(\Sigma)$, the subset
$$
\mathcal M_\Pi=L(\Sigma)\setminus\Pi^\pitchfork=\{\Lambda\in
L(\Sigma):\Lambda\cap\Pi\ne 0\}
$$
is called the {\it train} of $\Pi$.
Let $\Lambda_0\in\mathcal M_\Pi,\ \dim(\Lambda_0\cap\Pi)=k$. Assume
that $\Delta$ is transversal to both $\Lambda_0$ and $\Pi$ (i.e.
$\Delta\in\Lambda_0^\pitchfork\cap\Pi^\pitchfork$). The mapping
$\Lambda\mapsto S_\Lambda$ gives a regular parametrization of the
neighborhood of $\Lambda_0$ in $\mathcal M_\Pi$ by a neighborhood
of a corank $k$ matrix in the set of all degenerate symmetric
$n\times n$-matrices. A basic perturbation theory for symmetric
matrices now implies that a small enough neighborhood of
$\Lambda_0$ in $\mathcal M_\Pi$ is diffeomorphic to the Cartesian
product of a neighborhood of the origin of the cone of all degenerate
symmetric $k\times k$-matrices and a
$(n(n+1)-k(k+1))/2$-dimensional smooth manifold (see \cite[Lemma
2.2]{qu} for details). We see that $\mathcal M_\Pi$ is not a
smooth submanifold of $L(\Sigma)$ but a union of smooth strata,
$\mathcal M_\Pi=\bigcup\limits_{k>0}\mathcal M_\Pi^{(k)}$, where
$\mathcal M_\Pi^{(k)}=\{\Lambda\in
L(\Sigma):\dim(\Lambda\cap\Pi)=k\}$ is a smooth submanifold of
$L(\Sigma)$ of codimension $k(k+1)/2$.

Let $\Lambda(\tau),\ \tau\in[t_0,t_1]$ be a smooth family of
Lagrangian subspaces (a smooth curve in $L(\Sigma)$) and
$\Lambda(t_0),\Lambda(t_1)\in\Pi^\pitchfork$. We are going to
define the intersection number of $\Lambda(\cdot)$ and $\mathcal
M_\Pi$. It is called the Maslov index and is denoted
$\mu_\Pi(\Lambda(\cdot))$. Crucial property of this index is its
homotopy invariance: given a homotopy $\Lambda^s(\cdot)$,
$s\in[t_0,t_1]$ such that
$\Lambda^s(t_0),\Lambda^s(t_1)\in\Pi^\pitchfork \ \forall
s\in[0,1]$, we have
$\mu_\Pi(\Lambda^0(\cdot))=\mu_\Pi(\Lambda^1(\cdot))$.

It is actually enough to define $\mu_\Pi(\Lambda(\cdot))$ for the
curves which have empty intersection with
$\mathcal M_\Pi\setminus\mathcal M^{(1)}_\Pi$; the desired index
would have a well-defined extension to other curves by continuity.
Indeed, generic curves have empty intersection with
$\mathcal M_\Pi\setminus\mathcal M^{(1)}_\Pi$ and, moreover,
generic homotopy has empty intersection with
$\mathcal M_\Pi\setminus\mathcal M^{(1)}_\Pi$ since any of
submanifolds $\mathcal M_\Pi^{(k)},\ k=2,\ldots n$ has codimension
greater or equal to 3 in $L(\Sigma)$. Putting any curve in general
position by a small perturbation, we obtain the curve which
bypasses \linebreak $\mathcal M_\Pi\setminus\mathcal M^{(1)}_\Pi$,
and the invariance with respect to generic homotopies of the Maslov
index defined for generic curves
would imply that the value of the index does not depend on the
choice of a small perturbation.

What remains is to fix a ``coorientation" of the smooth
hypersurface $\mathcal M^{(1)}_\Pi$ in $L(\Sigma)$, i.~e. to
indicate the ``positive and negative sides" of the hypersurface.
As soon as we have a coorientation, we may compute
$\mu_\Pi(\Lambda(\cdot))$ for any curve $\Lambda(\cdot)$ which is
transversal to $\mathcal M^{(1)}_\Pi$ and has empty intersection
with $\mathcal M_\Pi\setminus\mathcal M^{(1)}_\Pi$. Maslov index
of $\Lambda(\cdot)$ is just the number of points where
$\Lambda(\cdot)$ intersects $\mathcal M^{(1)}_\Pi$ in the positive
direction minus the number of points where this curve intersects
$\mathcal M^{(1)}_\Pi$ in the negative direction. Maslov index of
any curve with endpoints out of $\mathcal M_\Pi$ is defined by
putting the curve in general position. Proof of the homotopy invariance
is the same as for usual intersection number of a curve with a
closed cooriented hypersurface (see, for instance, the nice elementary
book by J.~Milnor ``Topology from the differential viewpoint", 1965).

The coorientation is a byproduct of the following important structure on the
tangent spaces to $L(\Sigma)$. It happens that any tangent vector
to $L(\Sigma)$ at the point $\Lambda\in L(\Sigma)$ can be
naturally identified with a quadratic form on $\Lambda$. Her we
use the fact that $\Lambda$ is not just a point in the
Grassmannian but an $n$-dimensional linear space. To associate a
quadratic form on $\Lambda$ to the velocity
$\dot\Lambda(t)\in T_{\Lambda(t)}L(\Sigma)$ of a smooth curve
$\Lambda(\cdot)$ we proceed as follows: given $x\in\Lambda(t)$ we
take a smooth curve $\tau\mapsto x(\tau)$ in $\Sigma$ in such a
way that $x(\tau)\in\Lambda(\tau),\ \forall\tau$ and $x(\tau)=x$.
Then we define a quadratic form $\underline{\dot\Lambda}(t)(x),\
x\in\Lambda(t)$, by the formula
$\underline{\dot\Lambda}(t)(x)=\sigma(x,\dot x(t))$.

The point is that $\sigma(x,\dot x(t))$ does not depend on the
freedom in the choice of the curve $\tau\mapsto x(\tau)$, although
$\dot x(t)$ depends on this choice. Let us check the required
property in the coordinates.  We have $x=(z^\top,S_{\Lambda(t)}z)$
for some $z\in\mathbb R^n$ and
$x(\tau)=(z(\tau)^\top,S_{\Lambda(\tau)}z(\tau))$. Then
$$
\sigma(x,\dot x(t))=z^\top(\dot S_{\Lambda(t)}z+ S_{\Lambda(t)}\dot
z)-\dot z^\top  S_{\Lambda(t)}z=z^\top\dot S_{\Lambda(t)}z;
$$
vector $\dot z$ does not show up. We have obtained a coordinate
presentation of $\underline{\dot\Lambda}(t)$:
$$
\underline{\dot\Lambda}(t)(z^\top,S_{\Lambda(t)}z)=z^\top\dot
S_{\Lambda(t)}z,
$$
which implies that $\dot\Lambda\mapsto \underline{\dot\Lambda},\
\dot\Lambda\in T_\Lambda L(\Sigma)$ is an isomorphism of
$T_\Lambda L(\Sigma)$ on the linear space of quadratic forms on
$\Lambda$.

We are now ready to define the coorientation of $\mathcal
M^{(1)}_\Pi$. Assume that $\Lambda(t)\in\mathcal
M^{(1)}_\Pi$, i.~e. $\Lambda(t)\cap\Pi=\mathbb R x$ for some
nonzero vector $x\in\Sigma$. In coordinates, $x=(z^\top,0)$, where
$\mathbb R x=\ker S_{\Lambda(t)}$. It is easy to see that
$\dot\Lambda(t)$ is transversal to $\mathcal
M^{(1)}_\Pi$ (i.~e. $\dot S_{\Lambda(t)}$ is transversal to the
cone of degenerate symmetric matrices) if and only if
$\underline{\dot\Lambda}(t)(x)\ne 0$ (i.~e. $z^\top\dot S_{\Lambda(t)}z\ne
0$). Vector $x$ is defined up to a scalar multiplier and
$\underline{\dot\Lambda}(t)(\alpha x)=\alpha^2\underline{\dot\Lambda}(t)(x)$ so
that the sign of $\underline{\dot\Lambda}(t)(x)$ does not depend on the
selection of $x$.

\medskip\noindent {\bf Definition.} We say that $\Lambda(\cdot)$
intersects $\mathcal M^{(1)}_\Pi$ at the point $\Lambda(t)$ in the
positive (negative) direction if $\underline{\dot\Lambda}(t)(x)>0$ ($<0$).

\medskip This definition completes the construction of the Maslov
index. A weak point of the construction is the necessity to put
the curve in general position in order to compute the intersection
number. This does not look as an efficient way to thinks since
putting the curve in general position is nothing else but a
deliberate spoiling of a maybe nice and symmetric original object
that makes even more involved the nontrivial problem of the
localization
of its intersection with $\mathcal M_\Pi$. Fortunately, just
the fact that Maslov index is homotopy invariant leads to a very
simple and effective way of its computation without putting things
in general position and without looking for the intersection
points with $\mathcal M_\Pi$.
\begin{lemma} Assume that
$\Pi\cap\Delta=\Lambda(\tau)\cap\Delta=0,\
\forall\tau\in[t_0,t_1]$. Then
$\mu_\Pi(\Lambda(\cdot))=\mathrm{ind}S_{\Lambda(t_0)}-
\mathrm{ind}S_{\Lambda(t_1)}$, where $\mathrm{ind}S$ is the Morse
index of the quadratic form $z^\top Sz,\ z\in\mathbb R^n$.
\end{lemma}
{\bf Proof.} The matrices $S_{\Lambda(t_0)}$ and $S_{\Lambda(t_0)}$
are nondegenerate since
$\Lambda(t_0)\cap\Pi=\Lambda(t_1)\cap\Pi=0$ (we define the Maslov
index only for the curves whose endpoins are out of $\mathcal
M_\Pi$). The set of nondegenerate quadratic forms with a
prescribed value of the Morse index is a connected open subset of
the linear space of all quadratic forms in $n$ variables. Hence
homotopy invariance of the Maslov index implies that
$\mu_\Pi(\Lambda(\cdot))$ depends only on
$\mathrm{ind}S_{\Lambda(t_0)}$ and $\mathrm{ind}S_{\Lambda(t_1)}$.
It remains to compute $\mu_\Pi$ of sample curves in
$\Delta^\pitchfork$, say, for segments of the curve
$\Lambda(\cdot)$ such that
$$
S_{\Lambda(\tau)}=\left(\begin{smallmatrix} \tau-1& 0 &\ldots & 0\\
0 & \tau-2 &\ldots & 0\\ \vdots& \vdots & \ddots & \vdots\\
0 & 0 &\ldots & \tau-n \end{smallmatrix}\right). \eqno \square
$$

In general, given curve is not contained in the fixed coordinate
neighborhood $\Delta^\pitchfork$ but any curve can be divided into
segments $\Lambda(\cdot)|_{[\tau_i,\tau_{i+1}]},\ i=0,\ldots, l$, in
such a way that $\Lambda(\tau)\in\Delta_i^\pitchfork\ \
\forall\tau\in[\tau_i,\tau_{i+1}]$, where $\Delta_i\in\Pi^\pitchfork,\
i=0,\ldots,l$; then
$
\mu_\Pi(\Lambda(\cdot))=
\sum\limits_i\mu_{\Pi}\left(\Lambda(\cdot)|_{[\tau_i,\tau_{i+1}]}\right).
$

Lemma~I.2 implies the following useful formula which is valid for
the important class of {\it monotone increasing curves} in
the Lagrange Grassmannian, i.e. the curves $\Lambda(\cdot)$ such that
$\underline{\dot\Lambda}(t)$ are nonnegative quadratic forms:
$\underline{\dot\Lambda}(t)\ge 0,\ \forall t$.

\begin{corollary} Assume that $\underline{\dot\Lambda}(\tau)\ge
0,\ \forall\tau\in[t_0,t_1]$ and $\{\tau\in[t_0,t_1]:\Lambda(\tau)\cap\Pi\ne
0\}$ is a finite subset of $(t_0,t_1)$. Then
$$
\mu_\Pi\left(\Lambda(\cdot)\right)=\sum\limits_{\tau\in(t_0,t_1)}
\dim(\Lambda(\tau)\cap\Pi). \eqno \square
$$
\end{corollary}
Corollary I.1 can be also applied to the case of monotone
decreasing curves defined by the inequality
$\underline{\dot\Lambda}(t)\le 0,\ \forall t$; the change of
parameter $t\mapsto t_0+t_1-t$ makes the curve monotone
increasing and and change sign of the Maslov index.

\medskip Let me now recall that our interest to these symplectic
playthings was motivated by the conditional minimum problems. As
it was mentioned at the beginning of the section, we are going to
apply this stuff to the case $\Sigma=T_{\ell_\tau}(T^*M)$, $\ell_\tau\in
T^*_{z_\tau}M$, $\Pi=T_{\ell_\tau}(T^*_{z_\tau}M)$,
$\Lambda(\tau)=\mathcal L_{(\ell_\tau,w_\tau)}(\bar\Phi_\tau)$,
where $z_\tau=\Phi_\tau(w_\tau)$. In this case, not only $\Lambda$
but also $\Pi$ and even symplectic space $\Sigma$ depend on
$\tau$. We thus have to define Maslov index in such situation.
This is easy. We consider the bundle
$$
\{(\xi,\tau):\xi\in T_{\ell_\tau}(T^*M),\ t_0\le\tau\le t_1\}
\eqno (19)
$$
over the segment $[t_0,t_1]$ induced from $T(T^*M)$ by the mapping
$\tau\mapsto\ell_\tau$. Bundle (19) endowed with the symplectic
structure and its subbundle
$$\{(\xi,\tau):\xi\in T_{\ell_\tau}(T^*_{z_\tau}M)\}$$
are trivial as any bundle over a segment. More precisely, let
$t\in[t_0,t_1]$, $\Sigma_t=T_{\ell_t}(T^*M)$,
$\Pi_t=T_{\ell_t}(T_{z_t}^*M)$; then there exists a
continuous with respect to $\tau$ family of linear symplectic
mappings $\Xi_\tau:T_{\ell_\tau}(T^*M)\to\Sigma_t$ such that
$\Xi_\tau(T_{\ell_\tau}(T_{z_\tau}^*M))=\Pi_t,\
t_0\le\tau\le t_1,\ \Xi_t=\mathrm{Id}$. To any continuous family of
Lagrangian subspaces $\Lambda(\tau)\subset T_{\ell_\tau}(T^*M)$,
where $\Lambda(t_i)\cap\Pi_{t_i}=0,\ i=0,1$, we associate a curve
$\Xi_.\Lambda(\cdot):\tau\mapsto\Xi_\tau\Lambda(\tau)$ in the
Lagrange Grassmannian $L(\Sigma_t)$ and set
$\mu(\Lambda(\cdot))\stackrel{def}{=}\mu_{\Pi_t}(\Xi_.\Lambda(\cdot))$.
Homotopy invariance of the Maslov index implies that
$\mu_{\Pi_t}(\Xi_.\Lambda(\cdot))$ does not depend on the choice
of $t$ and $\Xi_\tau$.

\begin{theorem} Assume that $\dim\mathcal W<\infty$,
$$
\bar\Phi_\tau=(J_\tau,\Phi_\tau):\mathcal W\to\mathbb R\times
M,\quad \tau\in[t_0,t_1]
$$
is a continuous one-parametric family of smooth mappings and
$(\ell_\tau,w_\tau)$ is a continuous family of their Lagrangian points
such that $\ell_\tau\ne 0$, $w_\tau$ is a regular point of $\Phi_\tau$, and
$\ker D_{w_\tau}\Phi_\tau\cap\ker(D^2_{w_\tau}J_\tau-
\ell_\tau D^2_{w_\tau}\Phi_\tau)=0$, $t_0\le\tau\le t_1$. Let
$z_\tau=\Phi(w_\tau)$, $\Lambda(\tau)=\mathcal
L_{(\ell_\tau,w_\tau)}(\bar\Phi_\tau)$. If $ \mathrm{Hess}_{w_{t_i}}
(J_{t_i}\bigr|_{\Phi_{t_i}^{-1}(z_{t_i})}),\ i=1,2$, are nondegenerate, then
$$
\mathrm{ind}\,\mathrm{Hess}_{w_{t_0}}
(J_{t_0}\bigr|_{\Phi_{t_0}^{-1}(z_{t_0})})-
\mathrm{ind}\,\mathrm{Hess}_{w_{t_1}}
(J_{t_1}\bigr|_{\Phi_{t_1}^{-1}(z_{t_1})})=\mu(\Lambda(\cdot)).
$$
\end{theorem}
{\bf Remark.} If $\ell_\tau=0$, then $w_\tau$ is a critical point
of $J_\tau$ (without restriction to the level set of $\Phi_\tau$).
Theorem~I.2 can be extended to this situation (with the same
proof) if we additionally assume that
$\ker\mathrm{Hess}_{w_\tau}J_\tau=0$ for any $\tau$ such that
$\ell_\tau=0$.

\medskip\noindent {\bf Proof.} We introduce simplified notations:
$A_\tau=D_{w_\tau}\Phi_\tau$, $Q_\tau=D^2_{w_\tau}J_\tau-
\ell_\tau D^2_{w_\tau}\Phi_\tau$; the $\mathcal L$-derivative $\mathcal
L_{(\ell_\tau,w_\tau)}(\bar\Phi_\tau)=\Lambda(\tau)$ is uniquely
determined by the linear map $A_\tau$ and the symmetric bilinear form $Q_\tau$.
Fix local coordinates in the neighborhoods of $w_\tau$ and $z_\tau$ and set:
$$
\Lambda(A,Q)=\{(\zeta,A v) : \zeta A+Q(v,\cdot)=0\}\in L(\mathbb
R^{n*}\times\mathbb R^n);
$$
then $\Lambda_\tau=\Lambda(A_\tau,Q_\tau)$.

The assumption $\ker A_\tau\cap\ker Q_\tau=0$ implies the
smoothness of the mapping $(A,Q)\mapsto\Lambda(A,Q)$ for $(A,Q)$
close enough to $(A_\tau,Q_\tau)$. Indeed, as it is shown in the
proof of Proposition~I.4, this assumption implies that the mapping
$left_\tau:(\zeta,v)\mapsto\zeta A_\tau+Q_\tau(v,\cdot)$ is surjective.
Hence the kernel of the mapping $$
(\zeta,v)\mapsto\zeta A+Q(v,\cdot) \eqno (20)$$
smoothly depends on $(A,Q)$ for $(A,Q)$ close to
$(A_\tau,Q_\tau)$. On the other hand, $\Lambda(A,Q)$ is the image
of the mapping $(\zeta,v)\mapsto(\zeta,A v)$ restricted to the
kernel of map (20).

Now we have to disclose a secret which the attentive reader
already knows and is perhaps indignant with our lightness:
$Q_\tau$ is not a well-defined bilinear form on
$T_{w_\tau}\mathcal W$, it essentially depends on the choice of
local coordinates in $M$. What are well-defined is the mapping
$Q_\tau\bigr|_{\ker A_\tau}:\ker A_\tau\to T^*_{w_\tau}\mathcal W$
(check this by yourself or see \cite[Subsec.~2.3]{sy}), the map
$A_\tau:T_{w_\tau}\mathcal W\to T_{z_\tau}M$ and, of course, the
Lagrangian subspace $\Lambda(\tau)=\mathcal
L_{(\ell_\tau,w_\tau)}(\bar\Phi_\tau)$.
By the way, the fact that
$Q_\tau\bigr|_{\ker A_\tau}$ is well-defined guarantees that assumptions
of Theorem~I.2 do not depend on the coordinates choice.

Recall that any local coordinates $\{z\}$ on $M$ induce coordinates
$\{(\zeta,z):\zeta\in\mathbb R^{n*},z\in\mathbb R^n\}$ on $T^*M$
and $T^*_zM=\{(\zeta,0):\zeta\in\mathbb R^{n*}\}$ in the induced coordinates.

\begin{lemma} Given $\hat z\in M$, $\ell\in T^*_{\hat z}M\setminus\{0\}$, and a Lagrangian
subspace $\Delta\in T_\ell(T^*_{\hat z}M)^\pitchfork\subset L(T_\ell(T^*M))$, there exist
centered at $\hat z$ local coordinates on $M$ such that
$\Delta=\{(0,z):z\in\mathbb R^n\}$ in the induced coordinates on $T_\ell(T^*M)$.
\end{lemma}
{\bf Proof.} Working in arbitrary local coordinates we have $\ell=(\zeta_0,0)$,
$\Delta=\{(Sz,z):z\in\mathbb R^n\}$, where $S$ is a symmetric matrix. In other words,
$\Delta$ is the tangent space at $(\zeta_0,0)$ to the graph of the differential of the
function $a(z)=\zeta_0z+\frac 12z^\top Sz$. any smooth function with a nonzero differential
can be locally made linear by a smooth change of variables. To prove the lemma it is enough
to make a coordinates change  which kills second derivative of the function $a$, for instance:
$z\mapsto z+\frac 1{2|\zeta_0|^2}(z^\top Sz)\zeta_0^\top. \quad \square$

We continue the proof of Theorem~I.2. Lemma~I.3 gives us the way to
take advantage of the fact that $Q_\tau$ depends on the choice of
local coordinates in $M$. Indeed, bilinear form $Q_\tau$ is
degenerate if and only if
$\Lambda_\tau\cap\{(0,z):z\in\mathbb R^n\}\ne 0$; this immediately
follows from the relation
$$
\Lambda_\tau=\{(\zeta,A_\tau v):\zeta A_\tau+Q_\tau(v,\cdot)=0\}.
$$
Given $t\in[t_0,t_1]$ take a transversal to
$T_{\ell_t}(T^*_{z_t}M)$ and $\Lambda(t)$ Lagrangian subspace
$\Delta_t\subset T_{\ell_t}(T^*M)$ and centered at $z_t$ local
coordinates in $M$ such that $\Delta_t=\{(0,z):z\in\mathbb R^n\}$
in these coordinates. Then $\Lambda(\tau)$ is transversal to
$\{(0,z):z\in\mathbb R^n\}$ for all $\tau$ from a neighborhood
$O_t$ of $t$ in $[t_0,t_1]$. Selecting an appropriate finite
subcovering from the covering $O_t,\ t\in[t_0,t_1]$ of $[t_0,t_1]$
we can construct a subdivision
$t_0=\tau_0<\tau_1<\ldots<\tau_k<\tau_{k+1}=t_1$ of $[t_0,t_1]$
with the following property: $\forall i\in\{0,1,\ldots,k\}$ the
segment $\{z_\tau:\tau\in[\tau_i,\tau_{i+1}]\}$ of the curve
$z_\tau$ is contained in a coordinate neighborhood $\mathcal O^i$
of $M$ such that
$\Lambda_\tau\cap\{(0,z):z\in\mathbb R^n\}=0\
\forall\tau\in[\tau_i,\tau_{i+1}]$ in the correspondent local
coordinates.

We identify the form $Q_\tau$ with its symmetric matrix, i.e.
$Q_\tau(v_1,v_2)=v_1^\top Q_\tau v_2$.
Then $Q_\tau$ is a nondegenerate symmetric matrix and
$$
\Lambda(\tau)=\{(\zeta,-A_\tau Q^{-1}_\tau A^\top_\tau\zeta^\top\},
\quad\tau_i\le\tau\le\tau_{i+1}. \eqno (21)
$$
Now focus on the subspace $\Lambda(\tau_i)$; it has a nontrivial
intersection with
$\{(\zeta,0):\zeta\in\mathbb
R^{n*}\}=T_{\ell_{\tau_i}}(T^*_{z_{\tau_i}}M)$ if and only if the
matrix $A_{\tau_i}Q^{-1}_{\tau_i}A_{\tau_i}^\top$ is degenerate.
This is the matrix of the restriction of the nondegenerate
quadratic form $v\mapsto v^\top Q^{-1}_{\tau_i}v$ to the image of
the linear map $A^\top_{\tau_i}$. Hence
$A_{\tau_i}Q^{-1}_{\tau_i}A_{\tau_i}^\top$ can be made
nondegenerate by the arbitrary small perturbation of the map
$A_{\tau_i}:T_{w_{\tau_i}}\mathcal W\to T_{z_{\tau_i}}M$.
Such perturbations can be realized simultaneously for
$i=1,\ldots,k$\,\footnote{We do not need to perturb $A_{t_0}$
and $A_{t_{k+1}}$: assumption of the theorem and Lemma~I.1
guarantee the required nondegeneracy  property.}
by passing to a continuous family $\tau\mapsto A'_\tau,\
t_0\le\tau\le t_1$, arbitrary close and homotopic to the family
$\tau\mapsto A_\tau$. In fact, $A'_\tau$ can be chosen equal to
$A_\tau$ out of an arbitrarily small neighborhood of
$\{\tau_1,\ldots,\tau_k\}$. Putting now $A'_\tau$ instead of
$A_\tau$ in the expression for $\Lambda(\tau)$ we obtain a family
of Lagrangian subspaces $\Lambda'(\tau)$. This family is
continuous (see the paragraph containing formula (20)) and
homotopic to $\Lambda(\cdot)$. In particular, it has the same
Maslov index as $\Lambda(\cdot)$. In other words, we can assume
without lack of generality that $\Lambda(\tau_i)\cap
T_{\ell_{\tau_i}}(T^*_{z_{\tau_i}}M)=0,\ i=0,1,\ldots,k+1$.
Then
$\mu(\Lambda(\cdot))=
\sum\limits_{i=0}^k\mu\left(\Lambda(\cdot)\bigr|_{[\tau_i,\tau_{i+1}]}\right).$
 Moreover, it follows from (21) and Lemma~I.2 that
$$
\mu\left(\Lambda(\cdot)\bigr|_{[\tau_i,\tau_{i+1}]}\right)=
\mathrm{ind}(A_{\tau_{i+1}}Q^{-1}_{\tau_{i+1}}A^\top_{\tau_{i+1}})-
\mathrm{ind}(A_{\tau_i}Q^{-1}_{\tau_i}A^\top_{\tau_i}).
$$
Besides that, $\mathrm{ind}Q_{\tau_i}=\mathrm{ind}Q_{\tau_{i+1}}$
since $Q_{\tau}$ is nondegenerate for all
$\tau\in[\tau_i,\tau_{i+1}]$ and continuously depends on $\tau$.

Recall that
$\mathrm{Hess}_{w_\tau}\left(J_\tau\bigr|_{\Phi^{-1}(z_\tau)}\right)=
Q_\tau\bigr|_{\ker A_\tau}.$ In order to complete proof of the
theorem it remains to show that
$$
\mathrm{ind}Q_\tau=\mathrm{ind}\left(Q_\tau\bigr|_{\ker
A_\tau}\right)+\mathrm{ind}(A_\tau Q^{-1}_\tau A^\top_\tau) \eqno
(22)
$$
for $\tau=\tau_i,\tau_{i+1}$.

Let us rearrange the second term in the right-hand side of (22).
The change of variables $v=Q^{-1}_\tau A^\top_\tau z,\ z\in\mathbb
R^n$, implies: $\mathrm{ind}\left(A_\tau Q^{-1}_\tau A^\top_\tau\right)=
\mathrm{ind}\left(Q_\tau\bigr|_{\{Q^{-1}_\tau A^\top_\tau
z\,:\,z\in\mathbb R^n\}}\right).$ We have: $Q_\tau(v,\ker A_\tau)=0$
if and only if $Q_\tau(v,\cdot)=z^\top A_\tau$ for some $z\in\mathbb
R^n$, i.e. $v^\top Q_\tau=z^\top A_\tau$, $v=Q^{-1}_\tau A^\top_\tau
z$. Hence the right-hand side of (22) takes the form
$$
\mathrm{ind}Q_\tau=\mathrm{ind}\left(Q_\tau\bigr|_{\ker
A_\tau}\right)+\mathrm{ind}\left(Q_\tau\bigr|_{\{v\,:\,Q_\tau(v,\ker
A_\tau)=0\}}\right)
$$
and $Q_\tau\bigr|_{\{v:Q_\tau(v,\ker A_\tau)=0\}}$ is a
nondegenerate form for $\tau=\tau_i,\tau_{i+1}$. Now equality (22)
is reduced to the following elementary fact of linear algebra: If
$Q$ is a nondegenerate quadratic form on $\mathbb R^m$ and
$E\subset\mathbb R^m$ is a linear subspace, then
$\mathrm{ind}Q=\mathrm{ind}\left(Q|_E\right)+
\mathrm{ind}\left(Q|_{E^\bot_Q}\right)+\dim(E\cap E^\bot_Q),
$
where $E^\bot_Q=\{v\in\mathbb R^m:Q(v,E)=0\}$ and $E\cap
E^\bot_Q=\ker\left(Q|_E\right)=\ker\left(Q|_{E^\bot_Q}\right).
\qquad \square$

\medskip\noindent{\bf Remark.} Maslov index $\mu_\Pi$ is somehow more than
just the intersection number with $\mathcal M_\Pi$. It can be
extended, in a rather natural way, to all continuous curves in the
Lagrange Grassmannian including those whose endpoint belong to
$\mathcal M_\Pi$. This extension allows to get rid of the annoying
nondegeneracy assumption for
$\mathrm{Hess}_{w_{t_i}}
(J_{t_i}\bigr|_{\Phi_{t_i}^{-1}(z_{t_i})})$ in the statement of
Theorem~I.2. In general, Maslov index computes 1/2 of the difference
of the signatures of the Hessians which is equal to the difference
of the Morse indices in the degenerate case (see \cite{sy} for
this approach).

\section{Regular extremals}

A combination of the finite-dimensional Theorem~I.2 with the
limiting procedure of Theorem~I.1 and with homotopy invariance of
the Maslov index allows to efficiently compute Morse indices of
the Hessians for numerous infinite-dimensional problems. Here we
restrict ourselves to the simplest case of a regular extremal of
the optimal control problem.

We use notations and definitions of Sections 3, 4. Let $h(\lambda,u)$
be the Hamiltonian of a smooth optimal control system and
$\lambda_t,\ t_0\le t\le t_1$, be an extremal contained in the
regular domain $\mathcal D$ of $h$. Then $\lambda_t$ is a solution
of the Hamiltonian system $\dot\lambda=\vec H(\lambda)$, where
$H(\lambda)=h(\lambda,\bar u(\lambda)),\ \frac{\partial h}{\partial
u}h(\lambda,\bar u(\lambda))=0$.

Let $q(t)=\pi(\lambda_t),\ t_0,\le t\le t_1$ be the extremal path.
Recall that the pair $(\lambda_{t_0},\lambda_t)$ is a Lagrangian
multiplier for the conditional minimum problem defined on an open
subset of the space
$$
M\times L_\infty([t_0,t_1],U)=
\{(q_t),u(\cdot)):q\in M, u(\cdot)\in L_\infty([t_0,t_1],U)\},
$$
where $u(\cdot)$ is control and $q_t$ is the value at  $t$ of the
solution to the differential equation $\dot q=f(q,u(\tau)),\
\tau\in[t_0,t_1]$. In particular, $F_t(q_t,u(\cdot))=q_t$. The
cost is $J_{t_0}^{t_1}(q_t,u(\cdot))$ and constraints  are
$F_{t_0}(q_t,u(\cdot))=q(0),\ q_t=q(t)$.

Let us set $J_t(u)=J_{t_0}^t(q(t),u(\cdot)),\
\Phi_t(u)=F_{t_0}(q(t),u(\cdot))$. A covector $\lambda\in
T^*M$ is a Lagrange multiplier for the problem $(J_t,\Phi_t)$ if
and only if there exists an extremal $\hat\lambda_\tau,\
t_0\le\tau\le t$, such that $\lambda_{t_0}=\lambda,\
\hat\lambda_t\in T^*_{q(t)}M$. In particular, $\lambda_{t_0}$ is a
Lagrange multiplier for the problem $(J_t,\Phi_t)$ associated to the
control $u(\cdot)=\bar u(\lambda_.)$. Moreover, all sufficiently
close to $\lambda_{t_0}$ Lagrange multipliers for this problem are
values at $t_0$ of the solutions $\lambda(\tau),\ t_0\le\tau\le t$ to the
Hamiltonian system $\dot\lambda=\vec H(\lambda)$ with the boundary
condition $\lambda(t)\in T^*_{q(t)}M$.

We'll use exponential notations for one-parametric groups of
diffeomorphisms generated ordinary differential equations. In
particular, $e^{\tau\vec H}:T^*M\to T^*M,\ \tau\in\mathbb R$, is a
flow generated by the equation $\dot\lambda=\vec H(\lambda)$, so
that $\lambda(\tau')=e^{(\tau'-\tau)\vec H}(\lambda(\tau),\
\tau,\tau'\in\mathbb R,$ and Lagrange multipliers for the problem
$(J_t,\Phi_t)$ fill the $n$-dimensional submanifold
$e^{(t_0-t)\vec H}\left(T^*_{q(t)}M\right)$.

We set $\bar\Phi_t=(J_t,\Phi_t)$; it is easy to see that the
$\mathcal L$-derivative $\mathcal L_{(\lambda_{t_0},u)}(\bar\Phi_t)$
is the tangent space to $e^{(t_0-t)\vec
H}\left(T^*_{q(t)}M\right)$, i.e. $\mathcal
L_{(\lambda_{t_0},u)}(\bar\Phi_t)=e_*^{(t_0-t)\vec H}T_{\lambda_t}
\left(T^*_{q(t)}M\right)$. Indeed, let us recall the construction
of the $\mathcal L$-derivative. First we linearize the equation
for Lagrange multipliers at $\lambda_{t_0}$. Solutions of the
linearized equation form an isotropic subspace
$\mathcal L^0_{(\lambda_{t_0},u)}(\bar\Phi_t)$ of the symplectic
space $T_{\lambda_{t_0}}(T^*M)$. If
$\mathcal L^0_{(\lambda_{t_0},u)}(\bar\Phi_t)$ is a Lagrangian
subspace (i.e. $\dim\mathcal L^0_{(\lambda_{t_0},u)}(\bar\Phi_t)
=\dim M$), then $\mathcal L_{(\lambda^0_{t_0},u)}(\bar\Phi_t)=
\mathcal L_{(\lambda_{t_0},u)}(\bar\Phi_t)$, otherwise we need a
limiting procedure to complete the Lagrangian subspace. In the
case under consideration, $\mathcal L^0_{(\lambda_{t_0},u)}(\bar\Phi_t)
=e_*^{(t_0-t)\vec H}T_{\lambda_t}\left(T^*_{q(t)}M\right)$
has a proper dimension and thus coincides with \linebreak
$\mathcal L_{(\lambda_{t_0},u)}(\bar\Phi_t)$. We can check
independently that $e_*^{(t_0-t)\vec H}T_{\lambda_t}
\left(T^*_{q(t)}M\right)$ is Lagrangian: indeed,
$T_{\lambda_t}\left(T^*_{q(t)}M\right)$ is Lagrangian and
$e_*^{(t_0-t)\vec H}:T_{\lambda_t}(T^*M)\to T_{\lambda_{t_0}}(T^*M)$
is an isomorphism of symplectic spaces since Hamiltonian flows
preserve the symplectic form.

So $t\mapsto\mathcal L_{(\lambda_{t_0},u)}(\bar\Phi_t)$ is a
smooth curve in the Lagrange Grassmannian
$L\left(T_{\lambda_{t_0}}(T^*M)\right)$ and we can try to compute
Morse index of
$$
\mathrm{Hess}_u\left(J_{t_1}\bigr|_{\Phi^{-1}_{t_1}(q(t_0))}\right)=
\mathrm{Hess}_u\left(J_{t_0}^{t_1}\bigr|_{F^{-1}_{t_0}(q(t_0))\cap F^{-1}_{t_1}
(q(t_1))}\right)
$$
via the Maslov index of this curve. Of course, such a computation
has no sense if the index is infinite.

\begin{prop}{\rm (Legendre condition)} If quadratic form
$\frac{\partial^2h}{\partial u^2}(\lambda_t,u(t))$ is negative definite for
any $t\in[t_0,t_1]$, then $\mathrm{ind}\,
\mathrm{Hess}_u\left(J_{t_1}\bigr|_{\Phi^{-1}_{t_1}(q(t_0))}\right)<\infty$
and
$\mathrm{Hess}_u\left(J_{t}\bigr|_{\Phi^{-1}_{t}(q(t_0))}\right)$
is positive definite for any $t$ sufficiently close to (and
strictly greater than) $t_0$. If
$\frac{\partial^2h}{\partial u^2}(\lambda_t,u(t))\nleq 0$
 for some $t\in[t_0,t_1]$, then
$\mathrm{ind}\,
\mathrm{Hess}_u\left(J_{t_1}\bigr|_{\Phi^{-1}_{t_1}(q(t_0))}\right)=\infty$.
\end{prop}
We do not give here the proof of this well-known result; you can
find it in many sources (see, for instance, the textbook
\cite{as}). It is based on the fact that
$
\frac{\partial^2h}{\partial u^2}(\lambda_t,u(t))=
\lambda(\frac{\partial^2f}{\partial u^2}(q(t),u(t)))-
\frac{\partial^2\varphi}{\partial u^2}(q(t),u(t))
$
is the infinitesimal (for the ``infinitesimally small interval" at $t$)
version of $\lambda_{t_0}D^2_u\Phi_{t_1}-D^2_uJ_{t_1}$ while
$\mathrm{Hess}_u\left(J_{t_1}\bigl|_{\Phi_{t_1}^{-1}(q(t_0))}\right)=
(D^2_uJ_{t_1}-\lambda_{t_0} D^2_w\Phi_{t_1})\bigl|_{\ker D_u\Phi_{t_1}}$.

Next theorem shows that in the `regular' infinite dimensional
situation of this section we may compute the Morse index similarly
to the finite dimensional case. The proof of the theorem requires
some information about second variation of optimal control
problems which is out of the scope of these notes. The required
information can be found in Chapters 20, 21 of \cite{as}.
Basically, it implies that finite dimensional arguments used in
the proof of Theorem~I.2 are legal also in our infinite dimensional
case.

We set: $\Lambda(t)=e_*^{(t_0-t)\vec H}T_{\lambda_t}
\left(T^*_{q(t)}M\right)$.
\begin{theorem} Assume that $\frac{\partial^2h}{\partial
u^2}(\lambda_t,u(t))$ is a negative definite quadratic form and
$u$ is a regular point of $\Phi_t,\ \forall t\in(t_0,t_1].$ Then:
\begin{itemize} \item The form
$\mathrm{Hess}_u\left(J_{t_1}\bigr|_{\Phi^{-1}_{t_1}(q(t_0))}\right)$ is
degenerate if and only if \newline
$\Lambda(t_1)\cap\Lambda(t_0)\ne 0$;
\item If $\Lambda(t_1)\cap\Lambda(t_0)=0$, then there exists $\bar
t>t_0$ such that
$$
\mathrm{ind}\,\mathrm{Hess}_u\left(J_{t_1}\bigr|_{\Phi^{-1}_{t_1}(q(t_0))}\right)=
-\mu\left(\Lambda(\cdot)\bigr|_{[\tau,t_1]}\right),\quad
\forall\tau\in(t_0,\bar t). \eqno \square
$$
\end{itemize}
\end{theorem}

Note that Legendre condition implies monotonicity of the curve
$\Lambda(\cdot)$; this property simplifies the evaluation of the Maslov
index. Fix some local coordinates in $M$ so that
$T^*M\cong\{(p,q)\in\mathbb R^{n*}\times\mathbb R^n\}$.

\begin{lemma} Quadratic form $\underline{\dot\Lambda}(t)$ is
equivalent (with respect to a linear change of variables) to the
form $-\frac{\partial^2H}{\partial p^2}(\lambda_t)=\frac{\partial\bar u}{\partial
p}^\top\frac{\partial^2h}{\partial u^2}(\lambda_t,\bar u(\lambda_t))
\frac{\partial\bar u}{\partial p}$.
\end{lemma}
{\bf Proof.} Equality $\frac{\partial^2H}{\partial p^2}=-\frac{\partial\bar u}{\partial
p}^*\frac{\partial^2h}{\partial u^2}\frac{\partial\bar u}{\partial
p}$ is an easy corollary of the identities $H(p,q)=h(p,q,\bar
u(p,q)),\ \frac{\partial h}{\partial u}\bigr|_{u=\bar u(p,q)}=0$.
Indeed, $\frac{\partial^2H}{\partial p^2}=
2\frac{\partial^2h}{\partial u\partial p}\frac{\partial\bar
u}{\partial p}+\frac{\partial\bar u}{\partial p}^\top
\frac{\partial^2h}{\partial u^2}\frac{\partial\bar u}{\partial p}$
and $\frac\partial{\partial p}\left(\frac{\partial h}{\partial
u}\right)=\frac{\partial^2h}{\partial p\partial
u}+\frac{\partial^2h}{\partial u^2}\frac{\partial\bar u}{\partial
p}=0$.
Further, we have:
$$
\frac d{dt}\Lambda(t)=\frac d{dt}e_*^{(t_0-t)\vec H}T_{\lambda_t}
\left(T^*_{q(t)}M\right)=e_*^{(t_0-t)\vec H}\frac
d{d\varepsilon}\Bigr|_{\varepsilon=0}e_*^{-\varepsilon\vec H}
T_{\lambda_{t+\varepsilon}}\left(T^*_{q(t+\varepsilon)}M\right).
$$
Set $\Delta(\varepsilon)=e_*^{-\varepsilon\vec H}
T_{\lambda_{t+\varepsilon}}\left(T^*_{q(t+\varepsilon)}M\right)\in
L\left(T_{\lambda(t)}(T^*M)\right)$. It is enough to prove that
$\underline{\dot\Delta(0)}$ is equivalent to
$-\frac{\partial^2H}{\partial p^2}(\lambda_t)$.
Indeed, $\dot\Lambda(t)=
e_*^{(t_0-t)\vec H}T_{\lambda_t}\dot\Delta(0)$, where
$$
e_*^{(t_0-t)\vec H}:T_{\lambda_t}(T^*M)\to
T_{\lambda_{t_0}}(T^*M)
$$
is a symplectic isomorphism. The association of the quadratic form
$\underline{\dot\Lambda}(t)$ on the subspace $\Lambda(t)$ to the
tangent vector $\dot\Lambda(t)\in
L\left(T_{\lambda_{t_0}}(T^*M)\right)$ is intrinsic, i.e. depends only
on the symplectic structure on $(T_{\lambda_{t_0}}(T^*M)$. Hence
$\underline{\dot\Delta}(0)(\xi)=\underline{\dot\Lambda}(t)\left(
e_*^{(t_0-t)\vec H}\xi\right)$, $\forall\xi\in\Delta(0)=
T_{\lambda_t}\left(T^*_{q(t)}M\right)$.

What remains, is to compute $\underline{\dot\Delta}(0)$; we do it in
coordinates. We have:
$$
\Delta(\varepsilon)=\left\{(\xi(\varepsilon),\eta(\varepsilon)):
\begin{array}{rcl}\dot\xi(\tau)&=&\xi\frac{\partial^2H}{\partial
p\partial q}(\lambda_{t-\tau})+\eta^\top\frac{\partial^2H}{\partial q^2}
(\lambda_{t-\tau}),\\ \dot\eta(\tau)&=&-\frac{\partial^2H}{\partial
p^2}(\lambda_{t-\tau})\xi^\top-\frac{\partial^2H}{\partial
q\partial p}(\lambda_{t-\tau})\eta,
\end{array} {\xi(0)\in\mathbb R^{n*}\atop \eta(0)=0} \right\},
$$
$$
\underline{\dot\Delta}(0)(\xi(0))=
\sigma\left((\xi(0),0),(\dot\xi(0),\dot\eta(0))\right)
=\xi(0)\dot\eta(0)=-\xi(0)\frac{\partial^2H}{\partial
p^2}(\lambda_t)\xi(0)^\top. \eqno \square
$$

Now combining Lemma~I.4 with Theorem~I.3 and Corollary~I.1 we obtain the
following version of the classical ``Morse formula"

\begin{corollary} Under conditions of Theorem~I.3, if
$\{\tau\in(t_0,t_1]:\Lambda(\tau)\cap\Lambda(t_0)\ne 0\}$ is a
finite subset of $(t_0,t_1)$, then
$$
\mathrm{ind}\,\mathrm{Hess}J_{t_1}\bigr|_{\Phi^{-1}_{t_1}(q(t_0))}=
\sum\limits_{\tau\in(t_0,t_1)}\dim(\Lambda(\tau)\cap\Lambda(t_0)).
$$
\end{corollary}

\part{Geometry of Jacobi curves}

\section{Jacobi curves}

Computation of the $\mathcal L$-derivative for regular extremals in
the last section has led us to the construction of curves in the Lagrange
Grassmannians which works for all Hamiltonian systems on the
cotangent bundles, independently on any optimal control problem.
Set $\Delta_\lambda=T_\lambda(T^*_qM)$, where $\lambda\in
T^*_qM,\ q\in M$. The curve $\tau\mapsto e^{-\tau\vec
H}_*\Delta_{e^{\tau\vec H}(\lambda)}$ in the Lagrange Grassmannian
$L\left(T_\lambda(T^*M)\right)$ is the result of the action of the
flow $e^{t\vec H}$ on the vector distribution
$\{\Delta_\lambda\}_{\lambda\in T^*M}$.
Now we are going to study differential geometry of these curves;
their geometry will provide us with a canonical connection on
$T^*M$ associated with the Hamiltonian system and with
curvature-type invariants. All that gives a far going
generalization (and a dynamical interpretation) of classical
objects from Riemannian geometry.

In fact, construction of the basic invariants does not need
symplectic structure and the Hamiltonian nature of the flow, we
may deal with more or less arbitrary pairs ({\sl vector field},
{\sl rank $n$ distribution}) on a $2n$-dimensional manifold $N$. The
resulting curves belong to the usual Grassmannian of all
$n$-dimensional subspaces in the $2n$-dimensional one. We plan to
work for some time in this more general situation and then come
back to the symplectic framework.

In these notes we mainly deal with the case of involutive
distributions (i.e. with $n$-foliations) just because our main
motivation and applications satisfy this condition. The reader can
easily recover more general definitions and construction by
himself.

So we consider a $2n$-dimensional smooth manifold $N$ endowed with a
smooth foliation of rank $n$. Let $z\in N$, by $E_z$ we denote the
passing through $z$ leaf of the foliation; then $E_z$ is an
$n$-dimensional submanifold of $N$. Point $z$ has a
coordinate neighborhood $O_z$ such that the restriction of the
foliation to $O_z$ is a (trivial) fiber bundle and the fibers
$E_{z'}^{loc},\ z'\in O_z,$ of this fiber bundle are connected
components of $E_{z'}\cap O_z$. Moreover, there exists a
diffeomorphism $O_z\cong\mathbb R^n\times\mathbb R^n$, where
$\mathbb R^n\times\{y\},\ y\in\mathbb R^n,$ are identified with the
fibers so that both the typical fiber and the base are
diffeomorphic to $\mathbb R^n$. We denote by $O_z/E^{loc}$ the
base of this fiber bundle and by $\pi:O_z\to O_z/E^{loc}$ the
canonical projection.

 Let $\zeta$ be a smooth vector field on $N$. Then
$z'\mapsto\pi_*\zeta(z')$, $z'\in E_z^{loc}$ is a smooth mapping
of $E_z^{loc}$ to $T_{\pi(z)}(O_z/E^{loc})$.
We denote the last mapping by
$\Pi_z(\zeta):E_z^{loc}\to T_{\pi(z)}(O_z/E^{loc})$.

\medskip
\noindent{\bf Definition.} We call $\zeta$ a {\it lifting} field
if $\Pi_z(\zeta)$ is a constant mapping $\forall z\in N$; The
field $\zeta$ is called {\it regular} if $\Pi_z(\zeta)$ is a
submersion, $z\in N$.

\medskip
The flow generated by the lifting field maps leaves of the
foliation in the leaves, in other words it is leaves-wise. On the
contrary, the flow generated by the regular field "smears" the
fibers over $O_z/E^{loc}$; basic examples are second order
differential equations on a manifold $M$ treated as the vector fields
on the tangent bundle $TM=N$.

Let us write things in coordinates:
We fix local coordinates acting in the domain $O\subset
N$, which turn the foliation into the Cartesian product of
vector spaces: $O\cong\{(x,y):x,y\in\mathbb R^n\}$,
$\pi:(x,y)\mapsto y$. Then vector field $\zeta$ takes the form
$\zeta=\sum\limits_{i=1}^n\left(a^i\frac\partial{\partial
x_i}+b^i\frac\partial{\partial y_i}\right)$, where $a^i,b^i$ are smooth
functions on $\mathbb R^n\times\mathbb R^n$.
The coordinate representation of the map $\Pi_z$ is:
$\Pi_{(x,y)}:x\mapsto \left(b^1(x,y),\ldots,b^n(x,y)\right)^\top$.
Field $\zeta$ is regular if and only if $\Pi_{(x,y)}$ are submersions;
in other words, if and only if
$\left(\frac{\partial b^i}{\partial x_j}\right)_{i,j=1}^n$
is a nondegenerate matrix. Field $\zeta$ is lifting if and only
if $\frac{\partial b^i}{\partial x_j}\equiv 0,\ i,j=1,\ldots,n$.

Now turn back to the coordinate free setting. The fibers $E_z$, $z\in N$
are integral manifolds of the
involutive distribution ${\cal E}=\{T_zE_z:z\in N\}$. Given
a vector field $\zeta$ on $N$, the (local) flow $e^{t\zeta}$
generated by $\zeta$, and $z\in N$ we define the family of
subspaces $$ J_z(t)=\left( e^{-t\zeta}\right)_*{\cal E}|_z\subset
T_zN. $$ In other words, $J_z(t)=\left(
e^{-t\zeta}\right)_*T_{e^{t\zeta}(z)}E_{e^{t\zeta}(z)}$,
$J_z(0)=T_zE_z$.

$J_x(t)$ is an $n$-dimensional subspace of $T_zN$, i.e. an element
of the Grassmannian $G_n(T_zN)$. We thus have (the germ of) a
curve $t\mapsto J_z(t)$ in $G_n(T_zN)$ which is called a {\it
Jacobi curve}.

\medskip\noindent{\bf Definition.} We say that field $\zeta$
is {\it k-ample} for an interger $k$ if $\forall z\in N$ and for
any curve $t\mapsto \hat J_z(t)$ in $G_n(T_zN)$ with the same
$k$-jet as $J_z(t)$ we have $\hat J_z(0)\cap\hat J_z(t)=0$ for all
$t$ close enough but not equal to 0. The field is called {\it
ample} if it is $k$-ample for some $k$.

\medskip It is easy to show that a field is 1-ample if and only if
it is regular.

\section{The cross-ratio}
Let $\Sigma$ be a $2n$-dimensional vector space, $v_0,v_1\in
G_n(\Sigma),\ v_0\cap v_1=0$. Than $\Sigma=v_0+v_1$. We denote by
$\pi_{v_0v_1}:\Sigma\to v_1$ the projector of $\Sigma$ onto $v_1$
parallel to $v_0$. In other words, $\pi_{v_0v_1}$ is a linear
operator on $\Sigma$ such that $\pi_{v_0v_1}\bigr|_{v_0}=0$,
$\pi_{v_0v_1}\bigr|_{v_1}=\mbox{id}$. Surely, there is a one-to-one
correspondence between pairs of transversal $n$-dimensional
subspaces of $\Sigma$ and rank $n$ projectors in $\mbox{gl}(\Sigma)$.

\begin{lemma} Let $v_0\in G_n(\Sigma)$; we set
$v_0^\pitchfork=\{v\in G_n(\Sigma):v\cap v_0=0\}$,
an open dense subset of $G_n(\Sigma)$.
Then $\{\pi_{vv_0} : v\in v_0^\pitchfork\}$ is an affine
subspace of $\mbox{gl}(\Sigma)$.
\end{lemma}
Indeed, any operator of the form
$\alpha\pi_{vv_0}+(1-\alpha)\pi_{wv_0}$, where $\alpha\in\mathbb
R$, takes values in $v_0$ and its restriction to $v_0$ is the
identity operator. Hence $\alpha\pi_{vv_0}+(1-\alpha)\pi_{wv_0}$
is the projector of $\Sigma$ onto $v_0$ along some subspace.

The mapping $v\mapsto\pi_{vv_0}$ thus serves as a local coordinate
chart on $G_n(\Sigma)$. These charts indexed by $v_0$ form a
natural atlas on $G_n(\Sigma)$.

Projectors $\pi_{vw}$ satisfy the following basic
relations:\footnote{Numbering of formulas is separate in each of two
parts of the paper}
$$
\pi_{v_0v_1}+\pi_{v_1v_0}=id, \quad
\pi_{v_0v_2}\pi_{v_1v_2}=\pi_{v_1v_2},\quad
\pi_{v_0v_1}\pi_{v_0v_2}=\pi_{v_0v_1}, \eqno (1)
$$
where $v_i\in G_n(\Sigma),\ v_i\cap v_j=0$ for $i\ne j$.
If $n=1$, then $G_n(\Sigma)$ is just the projective line
$\mathbb{RP}^1$; basic geometry of $G_n(\Sigma)$ is somehow
similar to geometry of the projective line for arbitrary $n$ as
well. The group $\mbox{GL}(\Sigma)$ acts transitively on $G_n(\Sigma)$.
Let us consider its standard action on $(k+1)$-tuples of points in
$G_n(\Sigma)$:
$$
A(v_0,\ldots,v_k)\stackrel{def}{=}
(Av_0,\ldots,Av_k),\quad A\in\mbox{GL}(\Sigma),\
v_i\in G_n(\Sigma).
$$
It is an easy exercise to check that the only invariants of a
triple $(v_0,v_1,v_2)$ of points of $G_n(\Sigma)$ for such an
action are dimensions of the intersections: $\dim(v_i\cap
v_j),\ 0\le i\le 2$, and $\dim(v_0\cap v_1\cap v_2)$. Quadruples
of points possess a more interesting invariant: a multidimensional
version of the classical cross-ratio.

\medskip\noindent{\bf Definition.} Let $v_i\in G_n(\Sigma),\ i=0,1,2,3$,
and $v_0\cap v_1=v_2\cap v_3=0.$ The cross-ratio of $v_i$ is the
operator $[v_0,v_1,v_2,v_3]\in\mbox{gl}(v_1)$ defined by the formula:
$$
[v_0,v_1,v_2,v_3]=\pi_{v_0v_1}\pi_{v_2v_3}\bigr|_{v_1}.
$$
{\sl Remark.} We do not lose information when restrict the product
$\pi_{v_0v_1}\pi_{v_2v_3}$ to $v_1$; indeed, this product takes
values in $v_1$ and its kernel contains $v_0$.

\medskip
For $n=1$, $v_1$ is a line and $[v_0,v_1,v_2,v_3]$ is a real
number. For general $n$, the Jordan form of the operator provides
numerical invariants of the quadruple $v_i,\ i=0,1,2,3$.

We will mainly use an infinitesimal version of the cross-ratio
that is an invariant $[\xi_0,\xi_1]\in\mbox{gl}(v_1)$ of a pair of
tangent vectors $\xi_i\in T_{v_i}G_n(\Sigma),\ i=0,1,$ where
$v_0\cap v_1=0$. Let $\gamma_i(t)$ be curves in $G_n(\Sigma)$ such
that $\gamma_i(0)=v_i,\ \frac d{dt}\gamma_i(t)\bigr|_{t=0}=\xi_i$,
$i=0,1$. Then the cross-ratio:
$[\gamma_0(t),\gamma_1(0),\gamma_0(\tau),\gamma_1(\theta)]$ is a well
defined operator on $v_1=\gamma_1(0)$ for all $t,\tau,\theta$ close enough to 0.
Moreover, it follows from
(1) that
$[\gamma_0(t),\gamma_1(0),\gamma_0(0),\gamma_1(0)]=$\linebreak
$[\gamma_0(0),\gamma_1(0),\gamma_0(t),\gamma_1(0)]=
[\gamma_0(0),\gamma_1(0),\gamma_0(0),\gamma_1(t)]=id$. We set
$$
[\xi_0,\xi_1]=\frac{\partial^2}{\partial t\partial\tau}
[\gamma_0(t),\gamma_1(0),\gamma_0(0),\gamma_1(\tau)]\bigr|_{v_1}\Bigr|_{t=\tau=0}
\eqno (2)
$$
It is easy to check that the right-hand side of (2) depends only
on $\xi_0,\xi_1$ and that $(\xi_0,\xi_1)\mapsto[\xi_0,\xi_1]$ is a
bilinear mapping from $T_{v_0}G_n(\Sigma)\times
T_{v_1}G_n(\Sigma)$ onto $gl(v_1)$.

\begin{lemma} Let $v_0,v_1\in G_n(\Sigma),\ v_0\cap v_1=0,\
\xi_i\in T_{v_i}G_n(\Sigma),\ and\
\xi_i=\frac d{dt}\gamma_i(t)\bigr|_{t=0},\ i=0,1$. Then
$
[\xi_0,\xi_1]=\frac{\partial^2}{\partial
t\partial\tau}
\pi_{\gamma_1(t)\gamma_0(\tau)}\bigr|_{v_1}\Bigr|_{t=\tau=0}
$
and $v_1,v_0$ are invariant subspaces of the operator
$\frac{\partial^2}{\partial
t\partial\tau}
\pi_{\gamma_1(t)\gamma_0(\tau)}\bigr|_{v_1}\Bigr|_{t=\tau=0}$.
\end{lemma}
{\bf Proof.} According to the definition,
$
[\xi_0,\xi_1]=\frac{\partial^2}{\partial
t\partial\tau}
(\pi_{\gamma_0(t)\gamma_1(0)}\pi_{\gamma_0(0)\gamma_1(\tau)})\bigr|_{v_1}\Bigr|_{t=\tau=0}.
$
The differentiation of the identities
$\pi_{\gamma_0(t)\gamma_1(0)}\pi_{\gamma_0(t)\gamma_1(\tau)}=
\pi_{\gamma_0(t)\gamma_1(0)},$\linebreak
$\pi_{\gamma_0(t)\gamma_1(\tau)}\pi_{\gamma_0(0)\gamma_1(\tau)}=
\pi_{\gamma_0(0)\gamma_1(\tau)}$ gives the equalities:
$$
\frac{\partial^2}{\partial
t\partial\tau}
(\pi_{\gamma_0(t)\gamma_1(0)}\pi_{\gamma_0(0)\gamma_1(\tau)})\Bigr|_{t=\tau=0}=
-\pi_{v_0v_1}\frac{\partial^2}{\partial
t\partial\tau}
\pi_{\gamma_0(t)\gamma_1(\tau)}\Bigr|_{t=\tau=0}
$$
$$
=-\frac{\partial^2}{\partial
t\partial\tau}
\pi_{\gamma_0(t)\gamma_1(\tau)}\Bigr|_{t=\tau=0}\pi_{v_0v_1}.
$$
It remains to mention that $\frac{\partial^2}{\partial
t\partial\tau}\pi_{\gamma_1(t)\gamma_0(\tau)}=
-\frac{\partial^2}{\partial t\partial\tau}
\pi_{\gamma_0(\tau)\gamma_1(t)}$.\quad $\square$

\section{Coordinate setting} Given $v_i\in G_n(\Sigma)$,
$i=0,1,2,3$, we coordinatize $\Sigma=\mathbb R^n\times\mathbb
R^n=\{(x,y) : x\in\mathbb R^n, y\in\mathbb R^n\}$ in such a way
that $v_i\cap\{(0,y):y\in\mathbb R^n\}=0$. Then there exist
$n\times n$-matrices $S_i$ such that $$ v_i=\{(x,S_ix):x\in\mathbb
R^n\},\quad i=0,1,2,3. \eqno (3) $$ The relation $v_i\cap v_j=0$
is equivalent to $\det(S_i-S_j)\ne 0$. If $S_0=0$, then the
projector $\pi_{v_0v_1}$ is represented by the $2n\times
2n$-matrix $\left(\begin{array}{cc} 0 & S_1^{-1}\\ 0 & I
\end{array}
\right).$ In general, we have
$$
\pi_{v_0v_1}=\left(\begin{array}{cc}
S_{01}^{-1}S_0 & -S_{01}^{-1}\\
S_1S_{01}^{-1}S_0 & -S_1S_{01}^{-1}
\end{array}
\right),
$$
where $S_{01}=S_0-S_1$. Relation (3) provides coordinates $\{x\}$
on the spaces $v_i$. In these coordinates, the operator
$[v_0,v_1,v_2,v_3]$ on $v_1$ is represented by the matrix:
$$
[v_0,v_1,v_2,v_3]=S_{10}^{-1}S_{03}S_{32}^{-1}S_{21},
$$
where $S_{ij}=S_i-S_j$.

We now compute the coordinate representation of the
infinitesimal cross-ratio. Let $\gamma_0(t)=\{(x,S_tx) :
x\in\mathbb R^n\}$, $\gamma_1(t)=\{(x,S_{1+t}x) :
x\in\mathbb R^n\}$ so that $\xi_i=\frac
d{dt}\gamma_i(t)\bigr|_{t=0}$ is represented by the matrix
$\dot S_i=\frac d{dt}S_t\bigr|_{t=i},\ i=0,1.$ Then
$[\xi_0,\xi_1]$ is represented by the matrix
$$
\frac{\partial^2}{\partial t\partial\tau}S^{-1}_{1t}
S_{t\tau}S^{-1}_{\tau 0}S_{01}\Bigr|_{\frac{t=0}{\tau=1}}=
\frac\partial{\partial t}S^{-1}_{1t}\dot S_1\Bigr|_{t=0}=
S^{-1}_{01}\dot S_0S^{-1}_{01}\dot S_1.
$$
So
$$[\xi_0,\xi_1]=S^{-1}_{01}\dot S_0S^{-1}_{01}\dot S_1. \eqno (4)$$

There is a canonical isomorphism
$T_{v_0}G_n(\Sigma)\cong\mbox{Hom}(v_0,\Sigma/v_0)$; it is defined
as follows. Let $\xi\in T_{v_0}G_n(\Sigma),\ \xi=\frac
d{dt}\gamma(t)|_{t=0}$, and $z_0\in v_0$. Take a smooth curve
$z(t)\in\gamma(t)$ such that $z(0)=z_0$. Then the residue class
$(\dot z(0)+v_0)\in\Sigma/v_0$ depends on $\xi$ and $z_0$ rather
than on a particular choice of $\gamma(t)$ and $z(t)$. Indeed, let
$\gamma'(t)$ be another curve in $G_n(\Sigma)$ whose velocity at
$t=0$ equals $\xi$. Take some smooth with respect to $t$ bases of
$\gamma(t)$ and $\gamma'(t)$:
$\gamma(t)=span\{e_1(t),\ldots,e_n(t)\},\
\gamma'(t)=span\{e'_1(t),\ldots,e'_n(t)\}$, where
$e_i(0)=e'_i(0),\ i=1,\ldots, n$; then
$\left(\dot e_i(0)-\dot e'_i(0)\right)\in v_0,\ i=1,\ldots,n$. Let
$z(t)=\sum\limits_{i=1}^n\alpha_i(t)e_i(t),\
z'(t)=\sum\limits_{i=1}^n\alpha'_i(t)e'_i(t)$, where
$\alpha_i(0)=\alpha'_i(0)$. We have:
$$
\dot z(0)-\dot
z'(0)=\sum\limits_{i=1}^n\left((\dot\alpha_i(0)-\dot\alpha'_i(0))e_i(0)+
\alpha'_i(0)(\dot e_i(0)-\dot e'_i(0))\right)\in v_0,
$$
i.e. $\dot z(0)+v_0=\dot z'(0)+v_0$.

We associate to $\xi$ the mapping $\bar\xi:v_0\to\Sigma/v_0$
defined by the formula $\bar\xi z_0=\dot z(0)+v_0$. The fact that
$\xi\to\bar\xi$ is an isomorphism of the linear spaces
$T_{v_0}G_n(\Sigma)$ and $\mbox{Hom}(v_0,\Sigma/v_0)$ can be
easily checked in coordinates. The matrices $\dot S_i$ above are
actually coordinate presentations of $\bar\xi_i,\ i=0,1$.

The standard action of the group $\mbox{GL}(\Sigma)$ on
$G_n(\Sigma)$ induces the action of $\mbox{GL}(\Sigma)$ on the
tangent bundle $TG_n(\Sigma)$. It is easy to see that the only
invariant of a tangent vector $\xi$ for this action is
$\mbox{rank}\bar\xi$ (tangent vectors are just ``double points" or
``pairs of infinitesimaly close points" and number
$(n-\mbox{rank}\bar\xi)$ is the infinitesimal version of the
dimension of the intersection for a pair of points in the
Grassmannian).
Formula (4) implies:
$$
\mbox{rank}[\xi_0,\xi_1]\le\min\{\mbox{rank}\bar\xi_0,\mbox{rank}\bar\xi_1\}.
$$

\section{Curves in the Grassmannian}
Let $t\mapsto v(t)$ be a germ at $\bar t$ of a smooth curve in the
Grassmannian $G_n(\Sigma)$.

\medskip\noindent {\bf Definition.} We say that the germ
$v(\cdot)$ is {\it ample} if $v(t)\cap v(\bar t)=0\ \forall
t\ne\bar t$ and the operator-valued function $t\mapsto\pi_{v(t)
v(\bar t)}$ has a pole at $\bar t$. We say that the germ $v(\cdot)$
is {\it regular} if the function $t\mapsto\pi_{v(t)v(\bar t)}$ has
a simple pole at $\bar t$. A smooth curve in $G_n(\Sigma)$ is
called ample (regular) if all its germs are ample (regular).

\medskip Assume that $\Sigma=\{(x,y):x,y\in\mathbb R^n\}$ is
coordinatized in such a way that $v(\bar t)=\{(x,0):x\in\mathbb
R^n\}$. Then $v(t)=\{(x,S_tx):x\in\mathbb R^n\}$, where $S(\bar
t)=0$ and $\pi_{v(t)v(\bar t)}=\left(\begin{array}{cc}
I & -S_t^{-1}\\
0 & 0
\end{array}
\right).$ The germ $v(\cdot)$ is ample if and only if the scalar
function $t\mapsto\det S_t$ has a finite order root at $\bar t$. The germ
$v(\cdot)$ is regular if and only if the matrix $\dot S_{\bar t}$
is not degenerate. More generally, the curve
$\tau\mapsto\{(x,S_\tau x):x\in\mathbb R^n\}$ is ample if and only if $\forall t$ the
function $\tau\mapsto\det(S_\tau-S_t)$ has a finite order root at
$t$. This curve is regular if and only if $\det\dot S_t\ne 0,\
\forall t.$
The intrinsic version of this coordinate characterization of
regularity reads: the curve $v(\cdot)$ is regular if and only if
the map $\bar{\dot v}(t)\in\mbox{Hom}(v(t),\Sigma/v(t))$ has rank
$n,\ \forall t$.

Coming back to the vector fields and their Jacobi curves (see
Sec.~8) one can easily check that a vector field is ample
(regular) if and only if its Jacobi curves are ample (regular).

Let $v(\cdot)$ be an ample curve in $G_n(\Sigma)$. We consider the
Laurent expansions at $t$ of the operator-valued function
$\tau\mapsto\pi_{v(\tau)v(t)}$,
$$
\pi_{v(\tau)v(t)}=\sum\limits_{i=-k_t}^m
(\tau-t)^i\pi^i_t+O(\tau-t)^{m+1}.
$$
Projectors of $\Sigma$ on the subspace $v(t)$ form an affine subspace
of $\mbox{gl}(\Sigma)$ (cf. Lemma~II.1). This fact implies that $\pi^0_t$ is
a projector of $\Sigma$ on $v(t)$; in other words,
$\pi^0_t=\pi_{v^\circ(t)v(t)}$ for some
$v^\circ(t)\in v(t)^\pitchfork$. We thus obtain another curve
$t\mapsto v^\circ(t)$ in $G_n(\Sigma)$, where
$\Sigma=v(t)\oplus v^\circ(t),\ \forall t$. The curve $t\mapsto
v^\circ(t)$ is called the {\it derivative curve} of the ample
curve $v(\cdot)$.

The affine space $\{\pi_{wv(t)}:w\in v(t)^\pitchfork\}$ is a
translation of the linear space $\frak N(v(t))=\{\frak n:\Sigma\to
v(t) \mid \frak n|_{v(t)}=0\}\subset \mbox{gl}(\Sigma)\}$ containing only
nilpotent operators. It is easy to see that $\pi^i_t\in\frak
N(v(t))$ for $i\ne 0$.

The derivative curve is not necessary ample. Moreover, it may be
nonsmooth and even discontinuous.

\begin{lemma} If $v(\cdot)$ is regular then $v^{\circ}(\cdot)$ is
smooth.
\end{lemma}
{\bf Proof.} We'll find the coordinate representation of
$v^{\circ}(\cdot)$. Let $v(t)=\{(x,S_tx):x\in\mathbb R^n\}$.
Regularity of $v(\cdot)$ is equivalent to the nondegeneracy of
$\dot S_t$. We have:
$$
\pi_{v(\tau)v(t)}=\left(\begin{array}{cc}
S_{\tau t}^{-1}S_\tau & -S_{\tau t}^{-1}\\
S_tS_{\tau t}^{-1}S_\tau & -S_tS_{\tau t}^{-1}
\end{array}
\right),
$$
where $S_{\tau t}=S_\tau-S_t$. Then
$S^{-1}_{\tau t}=(\tau-t)^{-1}\dot S_t^{-1}-\frac 12\dot S_t^{-1}\ddot
S_t\dot S_t^{-1}+O(\tau-t)$ as $\tau\to t$ and
$$
\pi_{v(\tau)v(t)}=(\tau-t)^{-1}\left(\begin{array}{cc}\dot
S_t^{-1}S_t & -\dot S_t^{-1} \\ S_t\dot S_t^{-1}S_t & -S_t\dot
S_t^{-1}
\end{array}\right)+
$$ $$ \left(\begin{array}{cc}I-\frac 12\dot S_t^{-1}\ddot S_t\dot
S_t^{-1}S_t & \frac 12 \dot S_t^{-1}\ddot S_t\dot S_t^{-1}\\
S_t-\frac 12S_t\dot S_t^{-1}\ddot S_t\dot S_t^{-1}S_t & \frac
12S_t\dot S_t^{-1}\ddot S_t\dot
S_t^{-1}\end{array}\right)+O(\tau-t). $$ We set $A_t=-\frac 12\dot
S_{t}^{-1}\ddot S_{t}\dot S_{t}^{-1}$; then
$\pi_{v^\circ(t)v(t)}=\left(\begin{array}{cc}I+A_tS_t & -A_t\\
S_t+S_tA_tS_t &-S_tA_t\end{array}\right)$ is smooth with respect
to $t$. Hence $t\mapsto v^\circ(t)$ is smooth. We obtain: $$
v^\circ(t)=\left\{(A_ty,y+S_tA_ty):y\in\mathbb R^n\right\}. \eqno
(5) $$

\section{The curvature}

\noindent{\bf Definition.} Let $v$ be an ample curve and $v^\circ$
be the derivative curve of $v$. Assume that $v^\circ$ is
differentiable at $t$ and set $R_v(t)=[\dot v^\circ(t),\dot
v(t)]$. The operator $R_v(t)\in gl(v(t))$ is called the
{\it curvature} of the curve $v$ at $t$.

\medskip If $v$ is a regular curve, then $v^\circ$ is smooth, the
curvature is well-defined and has a simple coordinate
presentation. To find this presentation, we'll use formula (4)
applied to $\xi_0=\dot v^\circ(t),\ \xi_1=\dot v(t)$. As before,
we assume that $v(t)=\{(x,S_tx):x\in\mathbb R^n\}$; in particular,
$v(t)$ is transversal to the subspace $\{(0,y):y\in\mathbb R^n\}$.
In order to apply (4) we need an extra assumption on the coordinatization of
$\Sigma$: the subspace $v^\circ(t)$ has to be transversal to
$\{(0,y):y\in\mathbb R^n\}$ for given $t$. The last property is
equivalent to the nondegeneracy of the matrix $A_t$ (see (6)). It
is important to note that the final expression for $R_v(t)$ as
a differential operator of $S$ must be valid without this extra
assumption since the definition of $R_v(t)$ is intrinsic! Now we
compute: $v^\circ(t)=\{(x,(A_t^{-1}+S_t)x):x\in\mathbb R^n\},\quad
R_v(t)=[\dot v^\circ(t),\dot v(t)]=A_t\frac
d{dt}(A^{-1}_t+S_t)A_t\dot S_t
=(A_t\dot S_t)^2-\dot A_t\dot S_t=\frac 14(\dot S_t^{-1}\ddot
S_t)^2-\dot A_t\dot S_t.
$
We also have $\dot A\dot S=-\frac 12\frac d{dt}(\dot S^{-1}\ddot
S\dot S^{-1})\dot S=(\dot S^{-1})^2-\frac 12\dot S^{-1}\stackrel{\ldots}{S}$.
Finally,
$$
R_v(t)=\frac 12\dot S_t^{-1}\stackrel{\ldots}{S}_t-\frac 34(\dot S_t^{-1}\ddot S_t)^2=
\frac d{dt}\left((2\dot S_t)^{-1}\ddot S_t\right)-\left((2\dot
S_t)^{-1}\ddot S_t\right)^2, \eqno (6)
$$
the matrix version of the Schwartzian derivative.

Curvature operator is a fundamental invariant of the curve in the
Grassmannian. One more intrinsic construction of this operator,
without using the derivative curve, is provided by the following

\begin{prop} Let $v$ be a regular curve in $G_n(\Sigma)$. Then
$$
[\dot v(\tau),\dot v(t)]=(\tau-t)^{-2}\mbox{id} +\frac 13R_v(t)+O(\tau-t)
$$
as $\tau\to t$.
\end{prop}
{\bf Proof.} It is enough to check the identity in some
coordinates. Given $t$ we may assume that
$$
v(t)=\{(x,0):x\in\mathbb R^n\},\quad
v^\circ(t)=\{(0,y):y\in\mathbb R^n\}.
$$
Let $v(\tau)=\{(x,S_\tau x:x\in\mathbb R^n\}$, then $S_t=\ddot
S_t=0$ (see (5)). Moreover, we may assume that the bases of the subspaces
$v(t)$ and $v^\circ(t)$ are coordinated in such a way that $\dot
S_t=I$. Then $R_v(t)=\frac 12\stackrel{\ldots}{S}_t$ (see (6)). On the other
hand, formula (4) for the infinitesimal cross-ratio implies:
$$
[\dot v(\tau),\dot v(t)]=S^{-1}_\tau\dot S_\tau S^{-1}_\tau=
-\frac d{d\tau}(S^{-1}_\tau)=
$$
$$
-\frac
d{d\tau}\left((\tau-t)I+\frac{(\tau-t)^3}6\stackrel{\ldots}{S}_t\right)^{-1}
+O(\tau-t)=
$$
$$
-\frac
d{d\tau}\left((\tau-t)^{-1}I-\frac{(\tau-t)}6\stackrel{\ldots}{S}_t\right)
+O(\tau-t)=(\tau-t)^{-2}I+\frac 16\stackrel{\ldots}{S}_t+O(\tau-t).
$$ $\square$

Curvature operator is an invariant of the curves in
$G_n(\Sigma)$ with fixed parametrizations. Asymptotic presentation
obtained in Proposition~II.1
implies a nice chain rule for the curvature of the reparametrized
curves.

Let $\varphi:\mathbb R\to\mathbb R$ be a regular change of variables,
i.e. $\dot\varphi\ne 0,\ \forall t$. The standard imbedding $\mathbb
R\subset\mathbb{RP}^1=G_1(\mathbb R^2)$ makes $\varphi$ a regular
curve in $G_1(\mathbb R^2)$. As we know (see (6)), the curvature
of this curve is the Schwartzian of $\varphi$:
$$
R_\varphi(t)=\frac{\stackrel{\ldots}{\varphi}(t)}{2\dot\phi(t)}-\frac
34\left(\frac{\ddot\varphi(t)}{\dot\varphi(t)}\right)^2.
$$
We set $v_\varphi(t)=v(\varphi(t))$ for any curve $v$ in $G_n(\Sigma)$.

\begin{prop} Let $v$ be a regular curve in $G_n(\Sigma)$ and
$\varphi:\mathbb R\to\mathbb R$ be a regular change of variables.
Then
$$
R_{v_\varphi}(t)=\dot\varphi^2(t)R_v(\varphi(t))+R_\varphi(t). \eqno (7)
$$
\end{prop}
{\bf Proof.} We have
$$
[\dot v_\varphi(\tau),\dot v_\varphi(t)]=
(\tau-t)^{-2}\mbox{id} +\frac 13R_{v_\varphi}(t)+O(\tau-t).
$$
On the other hand,
$$
[\dot v_\varphi(\tau),\dot v_\varphi(t)]=
[\dot\varphi(\tau)\dot v(\varphi(\tau)),\dot\varphi(t)\dot v(\varphi(t))]=
\dot\varphi(\tau)\dot\varphi(t)[\dot v(\varphi(\tau)),\dot v(\varphi(t))]=
$$
$$
\dot\varphi(\tau)\dot\varphi(t)
\left((\varphi(\tau)-\varphi(t))^{-2}id +\frac
13R_v(\varphi(t))+O(\tau-t)\right)=
$$
$$
\frac{\dot\varphi(\tau)\dot\varphi(t)}
{((\varphi(\tau)-\varphi(t))^2}\mbox{id} +\frac
{\dot\varphi^2(t)}3R_v(\varphi(t))+O(\tau-t).
$$
We treat $\varphi$ as a curve in $\mathbb{RP}^1=G_1(\mathbb R^2)$.
Then $[\dot\varphi(\tau),\dot\varphi(t)]=\frac{\dot\varphi(\tau)\dot\varphi(t)}
{(\varphi(\tau)-\varphi(t))^2}$, see~(4). The one-dimensional
version of Proposition~II.1 reads:
$$
[\dot\varphi(\tau),\dot\varphi(t)]=(t-\tau)^{-2}+\frac 13R_\varphi(t)+O(\tau-t).
$$
Finally,
$$
[\dot v_\varphi(\tau),\dot v_\varphi(t)]=
(t-\tau)^{-2}+\frac 13\left(R_\varphi(t)+\dot\varphi^2(t)R_v(\varphi(t))\right)+O(\tau-t).
\quad \square
$$

The following identity is an immediate corollary of Proposition~II.2:
$$
\left(R_{v_\varphi}-\frac 1n(\mbox{tr}R_{v_\varphi})\mbox{id}\right)(t)=
\dot\varphi^2(t)\left(R_v-\frac 1n(\mbox{tr}R_v)
\mbox{id}\right)(\varphi(t)). \eqno (8)
$$

\medskip\noindent{\bf Definition.} An ample curve $v$ is called flat
if $R_v(t)\equiv 0$.

\medskip It follows from Proposition~II.1 that any small enough piece of
a regular curve can be
made flat by a reparametrization if and only if the curvature of
the curve is a scalar operator, i.e.
$R_v(t)=\frac 1n(\mbox{tr}R_v(t))\mbox{id}$. In the case of a nonscalar
curvature, one can use equality (8) to define a distinguished
parametrization of the curve and then derive invariants which do
not depend on the parametrization.

\medskip{\sl Remark.} In this paper we are mainly focused on the
regular curves. See paper \cite{az} for the version of the chain
rule which is valid for any ample curve and for basic invariants
of unparametrized ample curves.

\section{Structural equations}

Assume that $v$ and $w$ are two smooth curves in
$G_n(\Sigma)$ such that \linebreak $v(t)\cap w(t)=0,\ \forall t$.

\begin{lemma} For any $t$ and any $e\in v(t)$ there exists a unique
$f_e\in w(t)$ with the following property: $\exists$ a smooth
curve $e_\tau\in v(\tau),\ e_t=e$, such that
$\frac d{d\tau}e_\tau\bigr|_{\tau=t}=f_e$.
Moreover, the mapping $\Phi_t^{vw}:e\mapsto f_t$ is linear and for any
$e_0\in v(0)$ there exists a unique smooth curve $e(t)\in v(t)$
such that $e(0)=e_0$ and
$$\dot e(t)=\Phi_t^{vw}e(t),\quad \forall t. \eqno (9)$$
\end{lemma}
{\bf Proof.} First we take any curve $\hat e_\tau\in v(\tau)$ such
that $e_t=e$. Then $\hat e_\tau=a_\tau+b_\tau$ where $a_\tau\in
v(t),\ b_\tau\in w(t)$. We take $x_\tau\in v(\tau)$ such that $x_t=\dot a_t$
and set $e_\tau=\hat e_\tau+(t-\tau)x_\tau$. Then $\dot e_t=\dot
b_t$ and we put $f_e=\dot b_t$.

Let us prove that $\dot b_t$ depends only on $e$ and not on
the choice of $e_\tau$. Computing the difference of two admissible
$e_\tau$ we reduce the lemma to the following statement:
if $z(\tau)\in v(\tau),\ \forall\tau$ and $z(t)=0$, then $\dot
z(t)\in v(t)$.

To prove the last statement we take smooth
vector-functions $e_\tau^i\in v(\tau),\ i=1,\ldots,n$ such that
$v(\tau)=span\{e^1_\tau,\ldots,e^n_\tau\}$. Then
$z(\tau)=\sum\limits_{i=1}^n\alpha_i(\tau)e^i_\tau,\
\alpha_i(t)=0$. Hence $\dot
z(t)=\sum\limits_{i=1}^n\dot\alpha_i(t)e^i_t\in v_t.$

Linearity of the map $\Phi_t^{vw}$ follows from the uniqueness of
$f_e$. Indeed, if $f_{e^i}=\frac d{d\tau}e^i_\tau\bigr|_{\tau=t}$,
then $\frac
d{d\tau}(\alpha_1e^1_\tau+\alpha_2e^2_\tau)\bigr|_{\tau=t}=
\alpha_1f_{e^1}+\alpha_2f_{e^2}$; hence
$\alpha_1f_{e^1}+\alpha_2f_{e^2}=f_{\alpha_1e^1+\alpha_2e^2},\
\forall e^i\in v(t),\ \alpha_i\in\mathbb R,\ i=1,2$.

Now consider the smooth submanifold $V=\{(t,e) : t\in\mathbb R,\
e\in v(t)\}$ of $\mathbb R\times\Sigma$. We have
$(1,\Phi_t^{vw}e)\in T_{(t,e)}V$ since $(1,\Phi_t^{vw}e)$ is the
velocity of a curve $\tau\mapsto(\tau,e_\tau)$ in $V$. So
$(t,e)\mapsto(1,\Phi_t^{vw}e),\ (t,e)\in V$ is a smooth vector
field on $V$.
The curve $e(t)\in v(t)$ satisfies (9) if and only if $(t,e(t))$
is a trajectory of this vector field. Now the standard existence
and uniqueness theorem for ordinary differential equations
provides the existence of a unique solution to the Cauchy problem
for small enough $t$ while the linearity of the equation guarantees
that the solution is defined for all $t. \quad \square $

\medskip
It follows from the proof of the lemma that
$\Phi^{vw}_te=\pi_{v(t)w(t)}\dot e_\tau\bigr|_{\tau=t}$ for any
$e_\tau\in v(\tau)$ such that $v_t=e$. Let
$v(t)=\{(x,S_{vt}x):x\in\mathbb R^n\},\ w(t)=\{(x,S_{wt}x):x\in\mathbb
R^n\}$; the matrix presentation of $\Phi_t^{vw}$ in coordinates
$x$ is $(S_{wt}-S_{vt})^{-1}\dot S_{vt}$. Linear mappings $\Phi_t^{vw}$
and $\Phi_t^{wv}$ provide a factorization of the infinitesimal
cross-ratio $[\dot w(t),\dot v(t)]$. Indeed, equality (4) implies:
$$
[\dot w(t),\dot v(t)]=-\Phi_t^{wv}\Phi_t^{vw}. \eqno (10)
$$
Equality (9) implies one more useful
presentation of the infinitesimal cross-ratio: if $e(t)$ satisfies
(9), then
$$
[\dot w(t),\dot
v(t)]e(t)=-\Phi_t^{wv}\Phi_t^{vw}e(t)=-\Phi_t^{wv}\dot e(t)=
-\pi_{w(t)v(t)}\ddot e(t). \eqno (11)
$$
Now let $w$ be the derivative curve of $v$, $w(t)=v^\circ(t)$. It
happens that $\ddot e(t)\in v(t)$ in this case and (11) is reduced
to the {\it structural equation}:
$$
\ddot e(t)=-[\dot v^\circ(t),\dot v(t)]e(t)=-R_v(t)e(t),
$$
where $R_v(t)$ is the curvature operator. More precisely, we have
the following
\begin{prop} Assume that $v$ is a regular curve in $G_n(\Sigma)$,
$v^\circ$ is its derivative curve, and $e(\cdot)$ is a smooth
curve in $\Sigma$ such that $e(t)\in v(t), \ \forall t$.
Then $\dot e(t)\in v^\circ(t)$ if and only if $\ddot e(t)\in v(t).$
\end{prop}
{\bf Proof.} Given $t$, we take coordinates in such a way
that $v(t)=\{(x,0):x\in\mathbb R^n\},\
v^\circ(t)=\{(0,y):y\in\mathbb R^n\}$.
Then $v(\tau)=\{(x,S_\tau x):x\in\mathbb R^n\}$ for $\tau$ close enough to $t$,
where $S_t=\ddot S_t=0$ (see (5)).

Let $e(\tau)=\{(x(\tau),S_\tau x(\tau))\}$. The inclusion
$\dot e(t)\in v^\circ(t)$ is equivalent to the equality $\dot x(t)=0$.
Further,
$$
\ddot e(t)=\{\ddot x(t),\ddot S_tx(t)+2\dot S_t\dot x(t)+S_t\ddot
x(t)\}=\{\ddot x(t),2\dot S\dot x\}\in v(t).
$$
Regularity of $v$ implies the nondegeneracy of $\dot S(t)$. Hence
$\ddot e(t)\in v(t)$ if and only if $\dot x(t)=0. \quad \square$

Now equality (11) implies
\begin{corollary} If $\dot e(t)=\Phi_t^{vv^\circ}e(t)$, then
$\ddot e(t)+R_v(t)e(t)=0$.
\end{corollary}

Let us consider invertible linear mappings $V_t:v(0)\to v(t)$
defined by the relations $V_te(0)=e(t),\ \dot
e(\tau)=\Phi_{\tau}^{vv^\circ}e(\tau),\ 0\le\tau\le t$. It follows from the structural
equation that the curve $v$ is uniquely reconstructed from $\dot
v(0)$ and the curve $t\mapsto V_t^{-1}R_V(t)$ in
$\mbox{gl}(v(0))$. Moreover, let $v_0\in G_n(\Sigma)$ and
$\xi\in T_{v_0}G_n(\Sigma)$, where the map
$\bar\xi\in\mbox{Hom}(v_0,\Sigma/v_0)$ has rank $n$;
then for any smooth curve $t\mapsto A(t)$ in
$\mbox{gl}(v_0)$ there exists a unique regular curve $v$ such that
$\dot v(0)=\xi$ and $V^{-1}_tR_v(t)V_t=A(t)$. Indeed, let
$e_i(0),\ i=1,\ldots,n$, be a basis of $v_0$ and
$A(t)e_i(0)=\sum\limits_{j=1}^na_{ij}(t)e_j(0)$. Then
$v(t)=span\{e_1(t),\ldots,e_n(t)\}$, where
$$
\ddot e_i(\tau)+\sum
\limits_{j=1}^na_{ij}(\tau)e_j(\tau)=0,\ 0\le\tau\le t,  \eqno (12)
$$ are uniquely defined
by fixing the $\dot v(0)$.

The obtained classification of regular curves in terms of the
curvature is particularly simple in the case of a scalar curvature
operators $R_v(t)=\rho(t)\mbox{id}$. Indeed, we have
$A(t)=V_t^{-1}R_v(t)V_t=\rho(t)\mbox{id}$ and system (12) is
reduced to $n$ copies of the Hill equation $\ddot
e(\tau)+\rho(\tau)e(\tau)=0$.

Recall that all $\xi\in TG_n(\Sigma)$ such that
$\mbox{rank}\bar\xi=n$ are equivalent under the action of
$\mbox{GL}(\Sigma)$ on $TG_n(\Sigma)$ induced by the standard
action on the Grassmannian $G_n(\Sigma)$. We thus obtain

\begin{corollary} For any smooth scalar function $\rho(t)$ there
exists a unique, up to the action of $\mbox{GL}(\Sigma)$, regular
curve $v$ in $G_n(\Sigma)$ such that $R_v(t)=\rho(t)\mbox{id}$.
\end{corollary}

Another important special class is that of symmetric curves.

\medskip\noindent{\bf Definition.} A regular curve $v$ is called
{\it symmetric} if $V_tR_v(t)=R_v(t)V_t,\ \forall t$.

\medskip\noindent In other words, $v$ is symmetric if and only
the curve $A(t)=V^{-1}_tR_v(t)V_t$ in $\mbox{gl}(v(0))$ is
constant and coincides with $R_v(0)$. The structural equation
implies

\begin{corollary} For any $n\times n$-matrix $A_0$, there
exists a unique, up to the action of $\mbox{GL}(\Sigma)$,
symmetric curve $v$ such that $R_v(t)$ is similar to $A_0$.
\end{corollary}

The derivative curve $v^\circ$ of a regular curve $v$ is not
necessary regular. The formula $R_v(t)=\Phi_t^{v^\circ
v}\Phi_t^{vv^\circ}$ implies that $v^\circ$ is regular if and only
if the curvature operator $R_v(t)$ is nondegenerate for any $t$. Then we may
compute the second derivative curve
$v^{\circ\circ}=(v^\circ)^\circ$.

\begin{prop} A regular curve $v$ with nondegenerate curvature
operators is symmetric if and only if $v^{\circ\circ}=v$.
\end{prop}
{\bf Proof.} Let us consider system (12). We are going to apply
Proposition~II.3 to the curve $v^\circ$ (instead of $v$) and the
vectors $\dot e_i(t)\in v^\circ(t)$. According to Proposition~II.3,
$v^{\circ\circ}=v$ if and only if $\frac{d^2}{dt^2}\dot
e_i(t)\,\in v^\circ(t)$. Differentiating (12) we obtain that
$v^{\circ\circ}=v$ if and only if the functions $\alpha_{ij}(t)$
are constant. The last property is none other than a
characterization of symmetric curves.$\quad \square$

\section{Canonical connection}
Now we apply the developed theory of curves in the Grassmannian to
the Jacobi curves $J_z(t)$ (see Sec.~8).

\begin{prop} All Jacobi curves $J_z(\cdot),\ z\in N$, associated
to the given vector field $\zeta$ are regular (ample) if and only if
the field $\zeta$ is regular (ample).
\end{prop}
{\bf Proof.} The definition of the regular (ample) field is
actually the specification of the definition of the regular
(ample) germ of the curve in the Grassmannian: general definition
is applied to the germs at $t=0$ of the curves $t\mapsto J_z(t)$.
What remains is to demonstrate that other germs of these curves
are regular (ample) as soon as the germs at 0 are. The latter fact
follows from the identity
$$
J_z(t+\tau)=e_*^{-t\zeta}J_{e^{t\zeta}(z)}(\tau) \eqno (13)
$$
(which, in turn, is an immediate corollary of the identity
$e_*^{-(t+\tau)\zeta}=e_*^{-t\zeta}\circ e_*^{-\tau\zeta}$).
Indeed, (13) implies that the germ of $J_z(\cdot)$ at $t$ is the
image of the germ of $J_{e^{t\zeta}(\tau)}(\cdot)$ at 0 under the
fixed linear transformation $e_*^{-t\zeta}:T_{e^{t\zeta}(z)}N\to
T_zN$. The properties of the germs to be regular or ample survive
linear transformations since they are intrinsic properties. \quad
$\square$

Let $\zeta$ be an ample field. Then the derivative curves
$J^\circ_z(t)$ are well-defined. Moreover, identity (13) and the
fact that the construction of the derivative curve is intrinsic
imply:
$$
J^\circ_z(t)=e_*^{-t\zeta}J^\circ_{e^{t\zeta}(z)}(0). \eqno (14)
$$
The value at 0 of the derivative curve provides the splitting
$T_zM=J_z(0)\oplus J^\circ_z(0)$, where the first summand is the
tangent space to the fiber, $J_z(0)=T_zE_z$.

Now assume that $J^\circ_z(t)$ smoothly depends on $z$; this assumption
is automatically fulfilled in the case of a regular $\zeta$,
where we have the explicit coordinate presentation for
$J^\circ_z(t)$. Then the subspaces $J^\circ_z(0)\subset T_zN,\ z\in
N,$ form a smooth vector distribution, which is the direct
complement to the vertical distribution ${\cal E}=\{T_zE_z:z\in
N\}$. Direct complements to the vertical distribution are called
Ehresmann connections (or just nonlinear connections, even if
linear connections are their special cases). The Ehresmann
connection ${\cal E}_\zeta=\{J^\circ_z(0) : z\in N\}$ is called
the {\it canonical connection} associated with $\zeta$ and the
correspondent splitting $TN={\cal E}\oplus{\cal E}_\zeta$ is
called the {\it canonical splitting}. Our nearest goal is to give
a simple intrinsic characterization of ${\cal E}_\zeta$ which does
not require the integration of the equation $\dot z=\zeta(z)$ and
is suitable for calculations not only in local coordinates but
also in moving frames.

Let ${\cal F}=\{F_z\subset T_zN : z\in N\}$ be an Ehresmann
connection. Given a vector field $\xi$ on $E$ we denote
$\xi_{ver}(z)=\pi_{F_zJ_z(0)}\xi,\ \xi_{hor}(z)=\pi_{J_z(0)F_z}\xi$, the
``vertical" and the ``horizontal" parts of $\xi(z)$. Then
$\xi=\xi_{ver}+\xi_{hor}$, where $\xi_{ver}$ is a section of the
distribution ${\cal E}$ and $\xi_{hor}$ is a section of the
distribution ${\cal F}$. In general, sections of ${\cal E}$ are
called vertical fields and sections of ${\cal F}$ are called
horizontal fields.

\begin{prop} Assume that $\zeta$ is a regular field. Then ${\cal
F}={\cal E}_\zeta$ if and only if the equality
$$
[\zeta,[\zeta,\nu]]_{hor}=2[\zeta,[\zeta,\nu]_{ver}]_{hor}
\eqno (15)
$$
holds for any vertical vector field $\nu$. Here $[\, ,\,]$ is Lie
bracket of vector fields.
\end{prop}
{\bf Proof.} The deduction of identity (15) is based on
the following classical expression:
$$
\frac d{dt}e_*^{-t\zeta}\xi=e_*^{-t\zeta}[\zeta,\xi], \eqno (16)
$$
for any vector field $\xi$.

Given $z\in N$, we take coordinates in $T_zN$ in such a way that
$T_zN=\{(x,y): x,y\in\mathbb R^n\}$, where
$J_z(0)=\{(x,0): x\in\mathbb R^n\},\ J^\circ_z(0)=
\{(0,y): y\in\mathbb R^n\}$. Let $J_z(t)=\{(x,S_tx): x\in\mathbb
R^n\}$, then $S_0=\ddot S_0=0$ and $\det\dot S_0\ne 0$ due to the
regularity of the Jacobi curve $J_z$.

Let $\nu$ be a vertical vector field, $\nu(z)=(x_0,0)$ and
$\left(e_*^{-t\zeta}\nu\right)(z)=(x_t,y_t)$. Then
$(x_t,0)=\left(e_*^{-t\zeta}\nu\right)_{ver}(z),\
(0,y_t)=\left(e_*^{-t\zeta}\nu\right)_{hor}(z)$. Moreover,
$y_t=S_tx_t$ since $\left(e_*^{-t\zeta}\nu\right)(z)\in J_z(t)$.
Differentiating the identity $y_t=S_tx_t$ we obtain:
$\dot y_t=\dot S_tx_t+S_t\dot x_t.$
In particular, $\dot y_0=\dot
S_0x_0$. It follows from (16) that $(\dot
x_0,0)=[\zeta,\nu]_{ver},\ (0,\dot y_0)=[\zeta,\nu]_{hor}$. Hence
$(0,\dot S_0x_0)=[\zeta,\nu]_{hor}(z)$, where, I recall, $\nu$ is
any vertical field. Now we differentiate once more and evaluate
the derivative at 0:
$$
\ddot y_0=\ddot S_0x_0+2\dot S_0\dot
x_0+S_0\ddot x_0=2\dot S_0\dot x_0. \eqno (17)
$$
The Lie bracket presentations of the left and right hand sides of
(17) are: $(0,\ddot y_0)=[\zeta,[\zeta,\nu]]_{hor},\
(0,\dot S_0\dot x_0)=[\zeta,[\zeta,\nu]_{ver}]_{hor}$. Hence (17)
implies identity (15).

Assume now that $\{(0,y): y\in\mathbb R^n\}\ne J^\circ_z(0)$; then
$\ddot S_0x_0\ne 0$ for some $x_0$. Hence $\ddot y_0\ne 2\dot
S_0\dot x_0$ and equality (15) is violated. $\quad \square$

Inequality (15) can be equivalently written in the following form
that is often more convenient for the computations:
$$
\pi_*[\zeta,[\zeta,\nu]](z)=2\pi_*[\zeta,[\zeta,\nu]_{ver}](z),
\quad \forall z\in N. \eqno (18)
$$

Let $R_{J_z}(t)\in\mbox{gl}(J_z(t))$ be the curvature of the
Jacobi curve $J_z(t)$. Identity (13) and the fact that
construction of the Jacobi curve is intrinsic imply that
$$
R_{J_z}(t)=e_*^{-t\zeta}R_{J_{e^{t\zeta}(z)}}(0)e_*^{t\zeta}\bigr|_{J_z(t)}.
$$
Recall that $J_z(0)=T_zE_z$; the operator
$R_{J_z}(0)\in\mbox{gl}(T_zE_z)$ is called {\it the
curvature operator of the field} $\zeta$ at $z$. We introduce the
notation: $R_\zeta(z)\stackrel{def}{=}R_{J_z}(0)$; then
$R_\zeta=\left\{R_\zeta(z)\right\}_{z\in E}$ is an automorphism of
the ``vertical" vector bundle $\left\{T_zE_z\right\}_{z\in
M}$.

\begin{prop} Assume that $\zeta$ is an ample vector field and
$J_z^\circ(0)$ is smooth with respect to $z$. Let $TN={\cal
E}\oplus{\cal E}_\zeta$ be the canonical splitting. Then
$$
R_\zeta\nu=-[\zeta,[\zeta,\nu]_{hor}]_{ver} \eqno (19)
$$
for any vertical field $\nu$.
\end{prop}
{\bf Proof.} Recall that $R_{J_z}(0)=[\dot J_z^\circ(0),\dot
J_z(0)]$, where $[\cdot,\cdot]$ is the infinitesimal cross--ratio
(not the Lie bracket!). The presentation (10) of the infinitesimal
cross--ratio implies:
$$
R_\zeta(z)=R_{J_z}(0)=-\Phi_0^{J_z^\circ J_z}\Phi_0^{J_zJ_z^\circ},
$$
where $\Phi_0^{vw}e=\pi_{v(0)w(0)}\dot e_0$ for any smooth curve
$e_\tau\in v(\tau)$ such that $e_0=e$. Equalities (14) and (16)
imply:
$
\Phi_0^{J_zJ^\circ_z}\nu(z)=[\zeta,\nu]_{ver}(z), \ \forall
z\in M.$
Similarly, $\Phi_0^{J_z^\circ J_z}\mu(z)=[\zeta,\mu]_{hor}(z)$
for any horizontal field $\mu$ and any $z\in M$. Finally,
$$
R_\zeta(z)\nu(z)=-\Phi_0^{J_z^\circ J_z}\Phi_0^{J_zJ^\circ_z}=
-[\zeta,[\zeta,\nu]_{hor}]_{ver}(z). \eqno \square
$$

\section{Coordinate presentation}

We fix local coordinates acting in the domain ${\mathcal O}\subset
N$, which turn the foliation into the Cartesian product of
vector spaces: ${\mathcal O}\cong\{(x,y):x,y\in\mathbb R^n\}$,
$\pi:(x,y)\mapsto y$. Then vector field $\zeta$ takes the form
$\zeta=\sum\limits_{i=1}^n\left(a^i\frac\partial{\partial
x_i}+b^i\frac\partial{\partial y_i}\right)$, where $a^i,b^i$ are smooth
functions on $\mathbb R^n\times\mathbb R^n$. Below we use abridged
notations: $\frac\partial{\partial x_i}=\partial_{x_i},\
\frac{\partial\varphi}{\partial x_i}=\varphi_{x_i}$ etc. We also
use the standard summation agreement for repeating indices.

Recall the coordinate characterization of the regularity
property for the vector field $\zeta$. Intrinsic definition of
regular vector fields is done in Section~8; it is based on the
mapping $\Pi_z$ whose coordinate presentation is:
$\Pi_{(x,y)}:x\mapsto \left(b^1(x,y),\ldots,b^n(x,y)\right)^\top$.
Field $\zeta$ is regular if and only if $\Pi_y$ are submersions;
in other words, if and only if $\left(b^i_{x_j}\right)_{i,j=1}^n$
is a non degenerate matrix.

Vector fields $\partial_{x_i},\ i=1,\ldots,n$, provide a basis of
the space of vertical fields. As soon as coordinates are fixed,
any Ehresmann connection finds a unique basis of the form:
$$
\left(\partial_{y_i}\right)_{hor}=\partial_{y_i}+c_i^j\partial_{x_j},
$$
where $c_i^j,\ i,j=1,\ldots,n$, are smooth functions on
$\mathbb R^n\times\mathbb R^n$. To characterize a connection in
coordinates thus means to find functions $c_i^j$. In the case of
the canonical connection of a regular vector field, the functions
$c_i^j$ can be easily recovered from identity (18) applied to
$\nu=\partial_{x_i},\ i=1,\ldots,n$. We'll do it explicitly for
two important classes of vector fields: second order ordinary
differential equations and Hamiltonian systems.

A second order ordinary differential equation
$$
\dot y=x,\quad \dot x=f(x,y) \eqno (20)
$$
there corresponds to the vector field $\zeta=
f^i\partial_{x_i}+x_i\partial_{y_i}$, where
$f=(f_1,\ldots,f_n)^\top$. Let $\nu=\partial_{x_i}$; then
$$
[\zeta,\nu]=-\partial_{y_i}-f^j_{x_i}\partial_{x_j},\
[\zeta,\nu]_{ver}=(c_i^j-f^j_{x_i})\partial_{x_j},
$$
$$
\pi_*[\zeta,[\zeta,\nu]]=f_{x_i}^j\partial_{y_j},\
\pi_*[\zeta,[\zeta,\nu]_{ver}]=(f_{x_i}^j-c^j_i)\partial_{y_j}.
$$
Hence, in virtue of equality (18) we obtain that $c_i^j=\frac
12f^j_{x_i}$ for the canonical connection associated with the
second order differential equation (20).

Now consider a Hamiltonian vector field
$\zeta=-h_{y_i}\partial_{x_i}+h_{x_i}\partial_{y_i}$, where $h$ is
a smooth function on $\mathbb R^n\times\mathbb R^n$ (a
Hamiltonian). The field $\zeta$ is regular if and only if the
matrix $h_{xx}=\left(h_{x_ix_j}\right)_{i,j=1}^n$ is non
degenerate. We are going to characterize the canonical connection
associated with $\zeta$. Let $C=\left(c_i^j\right)_{i,j=1}^n$; the
straightforward computation similar to the computation made for
the second order ordinary differential equation gives the
following presentation for the matrix $C$:
$$
2\left(h_{xx}Ch_{xx}\right)_{ij}=h_{x_k}h_{x_ix_jy_k}-h_{y_k}h_{x_ix_jx_k}
-h_{x_iy_k}h_{x_kx_j}-h_{x_ix_k}h_{y_kx_j}
$$
or, in the matrix form:
$$
2h_{xx}Ch_{xx}=\{h,h_{xx}\}-h_{xy}h_{xx}-h_{xx}h_{yx},
$$
where $\{h,h_{xx}\}$ is the Poisson bracket:
$\{h,h_{xx}\}_{ij}=\{h,h_{x_ix_j}\}=h_{x_k}h_{x_ix_jy_k}-h_{y_k}h_{x_ix_jx_k}$.

Note that matrix $C$ is symmetric in the Hamiltonian case
(indeed, \linebreak $h_{xx}h_{yx}=(h_{xy}h_{xx})^\top$). This is not occasional
and is actually guaranteed by the fact that Hamiltonian flows
preserve symplectic form $dx_i\wedge dy_i$. See Section~17 for the
symplectic version of the developed theory.

As soon as we found the canonical connection, formula (19) gives
us the presentation of the curvature operator although the
explicit coordinate expression can be bulky. Let us specify the vector field more.
In the case of the Hamiltonian of a natural mechanical system,
$h(x,y)=\frac 12|x|^2+U(y)$, the canonical connection is trivial:
$c_i^j=0$; the matrix of the curvature operator is just $U_{yy}$.

Hamiltonian vector field associated to the Hamiltonian\linebreak
$h(x,y)=g^{ij}(y)x_ix_j$ with a non degenerate symmetric matrix
$\left(g^{ij}\right)_{i,j=1}^n$ generates a (pseudo-)Rie\-mann\-ian
geodesic flow. Canonical connection in this case is classical Levi
Civita connection and the curvature operator is Ricci operator of
(pseudo-)Rie\-mann\-ian geometry (see \cite[Sec.~5]{ag} for details).
Finally, Hamiltonian $h(x,y)=g^{ij}(y)x_ix_j+U(y)$ has the same
connection as Hamiltonion $h(x,y)=g^{ij}(y)x_ix_j$ while its
curvature operator is sum of Ricci operator and second
covariant derivative of $U$.

\section{Affine foliations}

Let $\left[{\cal E}\right]$ be the sheaf of germs of sections of
the distribution ${\cal E}=\{T_zE_z : z\in N\}$ equipped with the
Lie bracket operation. Then $\left[{\cal E}\right]_z$ is just the Lie
algebra of germs at $z\in M$ of vertical vector fields. Affine
structure on the foliation $E$ is a sub-sheaf $\left[{\cal
E}\right]^a\subset\left[{\cal E}\right]$ such that $\left[{\cal
E}\right]^a_z$ is an Abelian sub-algebra of $\left[{\cal
E}\right]_z$ and $\{\varsigma(z) : \varsigma\in\left[{\cal
E}\right]^a_z\}=T_zE_z,\ \forall z\in N$. A foliation with a fixed
affine structure is called the {\it affine foliation}.

The notion of the affine foliation generalizes one of the vector bundle.
In the case of the vector bundle, the sheaf $\left[{\cal
E}\right]^a$ is formed by the germs of vertical vector fields
whose restrictions to the fibers are constant (i.e. translation
invariant) vector fields on the fibers. In the next section we will
describe an important class of affine foliations which is not reduced
to the vector bundles.

\begin{lemma} Let ${\cal E}$ be an affine foliation,
$\varsigma\in\left[{\cal E}\right]^a_z$ and $\varsigma(z)=0$. Then
$\varsigma|_{E_z}=0$.
\end{lemma}
{\bf Proof.} Let
$\varsigma_1,\ldots,\varsigma_n\in\left[{\cal E}\right]^a_z$ be such
that $\varsigma_1(z),\ldots,\varsigma_n(z)$ form a basis of
$T_zE_z$. Then $\varsigma=b_1\varsigma_1+\cdots+b_n\varsigma_n$,
where $b_i$ are germs of smooth functions vanishing at $z$.
Commutativity of $\left[{\cal E}\right]^a_z$ implies:
$0=[\varsigma_i,\varsigma]=(\varsigma_ib_1)\varsigma_1+\cdots+
(\varsigma_ib_n)\varsigma_n$. Hence functions $b_i|_{E_z}$ are constants,
i.e. $b_i|_{E_z}=0,\ i=1,\ldots, n. \quad \square$

Lemma~II.5 implies that $\varsigma\in\left[{\cal E}\right]^a_z$ is
uniquely reconstructed from $\varsigma(z)$. This property permits
to define the {\it vertical derivative} of any vertical vector field
$\nu$ on $M$. Namely, $\forall v\in T_zE_z$ we set
$$
D_v\nu=[\varsigma,\nu](z),\ \mbox{where}\
\varsigma\in\left[{\cal E}\right]^a_z,\ \varsigma(z)=v.
$$
Suppose $\zeta$ is a regular vector field on the manifold $N$
endowed with the affine $n$-foliation. The canonical Ehresmann
connection ${\cal E}_\zeta$ together with the vertical derivative
allow to define a canonical linear connection $\nabla$ on the
vector bundle ${\cal E}$. Sections of the vector bundle ${\cal E}$
are just vertical vector fields. We set
$$
\nabla_\xi\nu=[\xi,\nu]_{ver}+D_\nu(\xi_{ver}),
$$
where $\xi$ is any vector field on $N$ and $\nu$ is a vertical
vector field. It is easy to see that $\nabla$ satisfies axioms of
a linear connection. The only non evident one is:
$\nabla_{b\xi}\nu=b\nabla_\xi\nu$ for any smooth function $b$. Let
$z\in N \ ,\varsigma\in\left[{\cal E}\right]^a_z$, and
$\varsigma(z)=\nu(z)$. We have
$$
\nabla_{b\xi}\nu=[b\xi,\nu]_{ver}+[\varsigma,b\xi_{ver}]=
$$
$$
b\left([\xi,\nu]_{ver}+[\varsigma,\xi_{ver}]\right)-(\nu
b)\xi_{ver}+(\varsigma b)\xi_{ver}.
$$
Hence
$$
(\nabla_{b\xi}\nu)(z)=
b(z)\left([\xi,\nu]_{ver}(z)+[\varsigma,\xi_{ver}](z)\right)=
(b\nabla_\xi\nu)(z).
$$

Linear connection $\nabla$ gives us the way to express Pontryagin
characteristic classes of the vector bundle ${\cal E}$ via the
regular vector field $\zeta$. Indeed, any linear connection
provides an expression for Pontryagin classes. We are going to
briefly recall the correspondent classical construction (see \cite{mi}
for details). Let
$R^\nabla(\xi,\eta)=[\nabla_\xi,\nabla_\eta]-\nabla_{[\xi,\eta]}$
be the curvature of linear connection $\nabla$. Then
$R^\nabla(\xi,\eta)\nu$ is $C^\infty(M)$-linear with respect to
each of three arguments $\xi,\eta,\nu$. In particular,
$R^\nabla(\cdot,\cdot)\nu(z)\in\bigwedge^2(T^*_zN)\otimes T_zE_z,\
z\in N.$ In other words,
$R^\nabla(\cdot,\cdot)\in\mbox{Hom}\left({\cal E},\bigwedge^2(T^*N)\otimes
{\cal E}\right)$.

Consider the commutative exterior algebra
$$
\bigwedge\nolimits^{ev}N=C^\infty(N)\oplus
\bigwedge\nolimits^2(T^*N)\oplus\cdots
\oplus\bigwedge\nolimits^{2n}(T^*N)
$$
of the even order differential forms
on $N$. Then $R^\nabla$ can be treated as an endomorphism of the
module $\bigwedge^{ev}N\otimes{\cal E}$ over algebra
$\bigwedge^{ev}N$, i. e. \linebreak
$R^\nabla\in\mbox{End}_{\bigwedge^{ev}N}\left(\bigwedge^{ev}M
\otimes{\cal E}\right)$. Now consider characteristic polynomial
\linebreak $\det(tI+\frac
1{2\pi}R^\nabla)=t^n+\sum\limits_{i=1}^n\phi_it^{n-i}$, where the
coefficient $\phi_i$ is an order $2i$ differential form on $N$.
All forms $\phi_i$ are closed; the forms $\phi_{2k-1}$ are exact
and the forms $\phi_{2k}$ represent the Pontryagin characteristic
classes, $k=1,\ldots,[\frac n2]$.

\section{Symplectic setting}

Assume that $N$ is a symplectic manifold endowed with a symplectic
form $\sigma$. Recall that a symplectic form is just a closed
non degenerate differential 2-form. Suppose $E$ is a Lagrange
foliation on the symplectic manifold $(N,\sigma)$; this means that
$\sigma|_{E_z}=0,\ \forall z\in N$. Basic examples are cotangent
bundles endowed with the standard symplectic structure: $N=T^*M,\
E_z=T^*_{\pi(z)}M$, where $\pi:T^*M\to M$ is the canonical
projection. In this case $\sigma=d\tau$, where $\tau=\{\tau_z :
z\in T^*M \}$ is the Liouville 1-form on $T^*M$ defined by the
formula: $\tau_z=z\circ\pi_*$. Completely integrable Hamiltonian
systems provide another important class of Lagrange foliations.
We'll briefly recall the correspondent terminology. Details can be
found in any introduction to symplectic geometry (for instance, in
\cite{ar}).

Smooth functions on the symplectic manifold are called
Hamiltonians. To any Hamiltonian there corresponds a Hamiltonian
vector field $\vec h$ on $M$ defined by the equation: $d
h=\sigma(\cdot,\vec h)$. The Poisson bracket $\{h_1,h_2\}$ of the
Hamiltonians $h_1$ and $h_2$ is the Hamiltonian defined by the
formula: $\{h_1,h_2\}=\sigma(\vec h_1,\vec h_2)=\vec h_1h_2$.
Poisson bracket is obviously anti-symmetric and satisfies the
Jacobi identity: $\{h_1,\{h_2,h_3\}\}+\{h_3,\{h_1,h_2\}\}+
\{h_2,\{h_3,h_1\}\}=0$. This identity is another way to say that
the form $\sigma$ is closed. Jacobi identity implies one more
useful formula: $\overrightarrow{\{h_1,h_2\}}=[\vec h_1,\vec h_2]$.

We say that Hamiltonians $h_1,\ldots,h_n$ are in involution if
$\{h_i,h_j\}=0$; then $h_j$ is constant along trajectories of the
Hamiltonian equation $\dot z=\vec h_i(z),\ i,j=1,\ldots,n$. We say
that $h_1,\ldots,h_n$ are independent if $d_zh_1\wedge\cdots\wedge
d_zh_n\ne 0,\ z\in N$. $n$ independent Hamiltonians in involution
form a {\it completely integrable system}. More precisely, any of
Hamiltonian equations $\dot z=\vec h_i(z)$ is completely
integrable with first integrals $h_1,\ldots,h_n$.

\begin{lemma} Let Hamiltonians $h_1,\ldots,h_n$ form a completely
integrable system. Then the $n$-foliation
$E_z=\{z'\in M : h_i(z')=h_i(z),\ i=1,\ldots,n\},\quad z\in N$,
is Lagrangian.
\end{lemma}
{\bf Proof.} We have $\vec h_ih_j=0,\ i,j=1,\ldots,n$, hence $\vec
h_i(z)$ are tangent to $E_z$. Vectors $\vec h_1(z),\ldots,\vec
h_n(z)$ are linearly independent, hence
$$
span\{\vec h_1(z),\ldots,\vec h_n(z)\}=T_zE_z.
$$
Moreover, $\sigma(\vec h_i,\vec h_j)=\{h_i,h_j\}=0$, hence
$\sigma|_{E_z}=0. \quad \square$

Any Lagrange foliation possesses a canonical affine structure.
Let $\left[{\cal E}\right]$ be the sheaf of germs of the
distribution ${\cal E}=\{T_zE_z : z\in N\}$ as in Section~16; then
$\left[{\cal E}\right]^a$ is the intersection of $\left[{\cal
E}\right]$ with the sheaf of germs of Hamiltonian vector fields.

We have to check that Lie algebra $\left[{\cal E}\right]^a_z$ is
Abelian and generates $T_zE_z,\ \forall z\in N$. First check the
Abelian property. Let $\vec h_1,\vec h_2\in\left[{\cal
E}\right]^a_z$; we have $[\vec h_1,\vec h_2]=\overrightarrow{\{h_1,h_2\}},\
\{h_1,h_2\}=\sigma(\vec h_1,\vec h_2)=0$, since $\vec h_i$ are
tangent to $E_z$ and $\sigma|_{E_z}=0$. The second property
follows from the Darboux--Weinstein theorem (see \cite{ar})
which states that all Lagrange foliations are locally equivalent.
More precisely, this theorem states that any $z\in M$ possesses
a neighborhood $O_z$ and
local coordinates which turn the restriction of the Lagrange
foliation $E$ to $O_z$ into the trivial bundle $\mathbb
R^n\times\mathbb R^n=\{(x,y) : x,y\in\mathbb R^n\}$ and,
simultaneously, turn $\sigma|_{O_z}$ into the form
$\sum\limits_{i=1}^ndx_i\wedge dy_i$. In this special coordinates,
the fibers become coordinate subspaces $\mathbb R^n\times\{y\},\
y\in\mathbb R^n$, and the required property is obvious: vector
fields $\frac\partial{\partial x_i}$ are Hamiltonian fields
associated to the Hamiltonians $-y_i,\ i=1,\ldots,n$.

\medskip
Suppose $\zeta$ is a Hamiltonian field on the symplectic manifold
endowed with the Lagrange foliation, $\zeta=\vec h$. Let
$\varsigma\in\left[{\cal E}\right]^a_z,\ \varsigma=\vec s$; then
$\varsigma h=\{s,h\}$. The field $\vec h$ is regular if and only
if the quadratic form $s\mapsto\{s,\{s,h\}\}(z)$ has rank $n$.
Indeed, in the `Darboux--Weinstein coordinates' this quadratic
form has the matrix $\{\frac{\partial^2 h}{\partial x_i\partial
x_j}\}_{i,j=1}^n$.

Recall that the tangent space $T_zN$ to the symplectic manifold $N$ is a
symplectic space endowed with the symplectic structure $\sigma_z$.
An $n$-dimensional subspace $\upsilon\subset T_zN$ is a
Lagrangian subspace if $\sigma_z|_\upsilon=0$. The set
$$
L(T_zN)=\{\upsilon\in G_n(T_zM) : \sigma_z|_\upsilon=0\}
$$
of all Lagrange subspaces of $T_zM$ is a Lagrange
Grassmannian.

Hamiltonian flow $e^{t\vec h}$ preserves the symplectic form,
$\left(e^{t\vec h}\right)^*\sigma=\sigma$. Hence
$\left(e^{t\vec h}\right)_*:T_zN\to T_{e^{t\vec h}(z)}N$
transforms Lagrangian subspaces in the Lagrangian ones. It follows
that the Jacobi curve
$J_z(t)=\left(e^{-t\vec h}\right)_*T_{e^{t\vec h}(z)}E_{e^{t\vec
h}(z)}$ consists of Lagrangian subspaces, $J_z(t)\in L(T_zN)$.

We need few simple facts on Lagrangian Grassmannians
(see Sec.~6 for the basic information and \cite[Sec.~4]{sy} for a consistent
description of their geometry). Let
$(\Sigma,\bar\sigma)$ be a $2n$-dimensional symplectic space and
$\upsilon_0,\upsilon_1\in L(\Sigma)$ be a pair of transversal
Lagrangian subspaces, $\upsilon_0\cap\upsilon_1=0$. Bilinear form
$\bar\sigma$ induces a non degenerate pairing of the spaces
$\upsilon_0$ and $\upsilon_1$ by the rule
$(e,f)\mapsto\bar\sigma(e,f),\ e\in\upsilon_0, f\in\upsilon_1$. To
any basis $e_1,\ldots,e_n$ of $\upsilon_0$ we may associate a
unique dual basis $f_1,\ldots,f_n$ of $\upsilon_1$ such that
$\bar\sigma(e_i,f_j)=\delta_{ij}$. The form $\bar\sigma$ is
totally normalized in the basis $e_1,\ldots,e_n,f_1,\ldots,f_n$ of
$\Sigma$, since $\sigma(e_i,e_j)=\sigma(f_i,f_j)=0$. It follows
that symplectic group
$$
\mbox{Sp}(\Sigma)=\{A\in\mbox{GL}(\Sigma):
\bar\sigma(Ae,Af)=\bar\sigma(e,f),\ e,f\in\Sigma\}
$$
acts transitively on the pairs of transversal Lagrangian subspaces.

Next result is a `symplectic specification' of Lemma~II.1 from
Section~9.

\begin{lemma} Let $\upsilon_0\in L(\Sigma)$; then
$\{\pi_{\upsilon\upsilon_0} : \upsilon\in\upsilon_0^\pitchfork\cap
L(\Sigma)\}$ is an affine subspace of the affine space
$\{\pi_{v\upsilon_0} : v\in\upsilon_0^\pitchfork\}$ characterized
by the relation:
$$
v\in\upsilon_0^\pitchfork\cap L(\Sigma)\ \Leftrightarrow\
\bar\sigma(\pi_{v\upsilon_0}\cdot,\cdot)+\bar\sigma(\cdot,\pi_{v\upsilon_0}\cdot)
=\bar\sigma(\cdot,\cdot).
$$
\end{lemma}
{\bf Proof.} Assume that $\upsilon_1\in\upsilon^{\pitchfork}_0\cap
L(\Sigma)$. Let $e,f\in\Sigma,\ e=e_0+e_1,\ f=f_0+f_1$ where
$e_i,f_i\in\upsilon_i,\ i=0,1$; then
$$
\bar\sigma(e,f)=\bar\sigma(e_0+e_1,f_0+f_1)=\bar\sigma(e_0,f_1)+\bar\sigma(e_1,f_0)=
$$
$$
\bar\sigma(e_0,f)+\bar\sigma(e,f_0)=\bar\sigma(\pi_{\upsilon_1\upsilon_0}e,f)+
\bar\sigma(e,\pi_{\upsilon_1\upsilon_0}f).
$$
Conversely, let $v\in\upsilon^\pitchfork_0$ is not a Lagrangian
subspace. Then there exist $e,f\in v$ such that $\bar\sigma(e,f)\ne
0$, while
$\bar\sigma(\pi_{v\upsilon_0}e,f)=\bar\sigma(e,\pi_{v\upsilon_0}f)=0.\quad
\square $

\begin{corollary} Let $v(\cdot)$ be an ample curve in
$G_n(\Sigma)$ and $v^\circ(\cdot)$ be the derivative curve of
$v(\cdot)$. If $v(t)\in L(\Sigma),\ \forall t$, then $v^\circ(t)\in
L(\Sigma)$.
\end{corollary}
{\bf Proof.} The derivative curve $v^\circ$ was defined in
Section~11. Recall that $\pi_{v^\circ(t)v(t)}=\pi^0_t$, where
$\pi^0_t$ is the free term of the Laurent expansion
$$
\pi_{v(\tau)v(t)}\approx\sum\limits_{i=-k_t}^\infty
(\tau-t)^i\pi^i_t.
$$
The free term $\pi^0_t$ belongs to the affine hull of
$\pi_{v(\tau)v(t)}$, when $\tau$ runs a neighborhood of $t$. Since
$\pi_{v(\tau)v(t)}$ belongs to the affine space
$\{\pi_{vv_0} : v\in v_0^\pitchfork\cap L(\Sigma)\}$, then
$\pi^0_t$ belongs to this affine space as well. $\quad\square$

We call a {\it Lagrange distribution} any rank $n$  vector
distribution $\{\Lambda_z\subset T_zN : z\in N\}$ on the
symplectic manifold $N$ such that $\Lambda_z\in L(T_zN),\ z\in N$.

\begin{corollary} Canonical Ehresmann connection ${\cal
E}_\zeta=\{J^{\circ}_z(0):z\in N\}$
associated to an ample Hamiltonian field $\zeta=\vec h$ is a
Lagrange distribution. $\quad\square$
\end{corollary}

It is clearly seeing in coordinates how Lagrange Grassmanian is
sitting in the usual one. Let
$\Sigma=\mathbb R^{n*}\times\mathbb R^n=\{(\eta,y):\eta\in\mathbb
R^{n*},y\in\mathbb R^n\}$. Then any $v\in\left(\{0\}\times\mathbb
R^n\right)^\pitchfork$ has a form $v=\{(y^\top,Sy) :y\in\mathbb
R^n\}$, where $S$ is an $n\times n$-matrix. It is easy to see that
$v$ is a Lagrangian subspace if and only if $S$ is a symmetric
matrix, $S=S^\top$.

\section{Monotonicity}

We continue to study curves in the Lagrange Grassmannian
$L(T_zN)$, in particular, the Jacobi curves
$t\mapsto \left(e^{-t\vec H}\right)_*T_{e^{t\vec H}(z)}E_{e^{t\vec H}(z)}$.
In Section~6 we identified the velocity $\dot\Lambda(t)$ of any
smooth curve $\Lambda(\cdot)$ in $L(T_zN)$ with a quadratic form
$\underline{\dot\Lambda}(t)$ on the subspace $\Lambda(t)\subset
T_zN$. Recall that the curve $\Lambda(\cdot)$ was called monotone
increasing if $\underline{\dot\Lambda}(t)\ge 0,\ \forall t$; it is
called monotone decreasing if $\underline{\dot\Lambda}(t)\le 0$.
It is called monotone in both cases.

\begin{prop} Set $\Lambda(t)=
\left(e^{-t\vec H}\right)_*T_{e^{t\vec H}(z)}E_{e^{t\vec H}(z)}$;
then quadratic form $\underline{\dot\Lambda}(t)$ is
equivalent (up to a linear change of variables) to the form
$$
\varsigma\mapsto-(\varsigma\circ\varsigma H)(e^{t\vec H}(z)),\quad
\varsigma\in[\mathcal E]^a_{e^{t\vec H}(z)}, \eqno (21)
$$
on $E_{e^{t\vec H}(z)}$.
\end{prop}
{\bf Proof.} Let $z_t=e^{t\vec H}(z)$, then
$$
\frac d{dt}\Lambda(t)=\frac d{dt}e_*^{(t_0-t)\vec H}T_{z_t}
E_{z_t}=e_*^{(t_0-t)\vec H}\frac
d{d\varepsilon}\Bigr|_{\varepsilon=0}e_*^{-\varepsilon\vec H}
T_{z_{t+\varepsilon}}E_{z_{t+\varepsilon}}.
$$
Set $\Delta(\varepsilon)=e_*^{-\varepsilon\vec H}
T_{z_{t+\varepsilon}}E_{z_{t+\varepsilon}}\in
L\left(T_{z_t}N\right)$.
It is enough to prove that
$\underline{\dot\Delta}(0)$ is equivalent to form (21).
Indeed, $\dot\Lambda(t)=
e_*^{(t_0-t)\vec H}T_{z_t}\dot\Delta(0)$, where
$$
e_*^{(t_0-t)\vec H}:T_{z_t}N\to
T_{z_{t_0}}N
$$
is a symplectic isomorphism. The association of the quadratic form
$\underline{\dot\Lambda}(t)$ on the subspace $\Lambda(t)$ to the
tangent vector $\dot\Lambda(t)\in
L\left(T_{z_{t_0}}N\right)$ is intrinsic, i.e. depends only
on the symplectic structure on $T_{z_{t_0}}N$. Hence
$\underline{\dot\Delta}(0)(\xi)=\underline{\dot\Lambda}(t)\left(
e_*^{(t_0-t)\vec H}\xi\right)$, $\forall\xi\in\Delta(0)=
T_{z_t}E_{z_t}$.

What remains, is to compute $\underline{\dot\Delta}(0)$; we do it in
the Darboux--Weinstein coordinates $z=(x,y)$. We have:
$\Delta(\varepsilon)=$
$$
\left\{(\xi(\varepsilon),\eta(\varepsilon)):
\begin{array}{rcl}\dot\xi(\tau)&=&\xi(\tau)\frac{\partial^2H}{\partial
x\partial y}(z_{t-\tau})+\eta(\tau)^\top\frac{\partial^2H}{\partial y^2}
(z_{t-\tau}),\\ \dot\eta(\tau)&=&-\frac{\partial^2H}{\partial
x^2}(z_{t-\tau})\xi(\tau)^\top-\frac{\partial^2H}{\partial
y\partial x}(z_{t-\tau})\eta(\tau),
\end{array} {\xi(0)=\xi\in\mathbb R^{n*}\atop \eta(0)=0\in\mathbb R^n} \right\},
$$
$$
\underline{\dot\Delta}(0)(\xi)=
\sigma\left((\xi,0),(\dot\xi(0),\dot\eta(0))\right)
=\xi\dot\eta(0)=-\xi\frac{\partial^2H}{\partial
x^2}(z_t)\xi^\top.
$$
Recall now that  form (21)
has matrix $\frac{\partial^2H}{\partial x^2}(z_t)$
in the Darboux--Weinstein coordinates. \quad $\square $

This proposition clearly demonstrates the importance of monotone
curves. Indeed, monotonicity of Jacobi curves is equivalent to the
convexity (or concavity) of the Hamiltonian on each leaf of the
Lagrange foliation. In the case of a cotangent bundle this means
the convexity or concavity of the Hamiltonian with respect to the
impulses. All Hamiltonians (energy functions) of mechanical
systems are like that! This is not an occasional fact but a
corollary of the list action principle. Indeed, trajectories of
the mechanical Hamiltonian system are extremals of the least
action principle and the energy function itself is the Hamiltonian
of the correspondent regular optimal control problem as it was
considered in Section~7. Moreover, it was stated in Section~7 that
convexity of the Hamiltonian with respect to the impulses is
necessary for the extremals to have finite Morse index. It turns
out that the relation between finiteness of the Morse index and
monotonicity of the Jacobi curve has a fundamental nature. A
similar property is valid for any, not necessary regular, extremal
of a finite Morse index. Of course, to formulate this property we
have first to explain what are Jacobi curve for non regular extremals.
To do that, we come back to the very beginning; indeed, Jacobi
curves appeared first as the result of calculation of the
$\mathcal L$-derivative at the regular extremal (see Sections 7,
8). On the other hand, $\mathcal L$-derivative is well-defined for
any extremal of the finite Morse index as it follows from
Theorem~I.1. We thus come to the following construction in which
we use notations and definitions of Sections 3, 4.

Let $h(\lambda,u)$
be the Hamiltonian of a smooth optimal control system,
$\lambda_t,\ t_0\le t\le t_1$, an extremal, and
$q(t)=\pi(\lambda_t),\ t_0,\le t\le t_1$ the extremal path.
Recall that the pair $(\lambda_{t_0},\lambda_t)$ is a Lagrangian
multiplier for the conditional minimum problem defined on an open
subset of the space
$$
M\times L_\infty([t_0,t_1],U)=
\{(q_t,u(\cdot)):q\in M, u(\cdot)\in L_\infty([t_0,t_1],U)\},
$$
where $u(\cdot)$ is control and $q_t$ is the value at  $t$ of the
solution to the differential equation $\dot q=f(q,u(\tau)),\
\tau\in[t_0,t_1]$. In particular, $F_t(q_t,u(\cdot))=q_t$. The
cost is $J_{t_0}^{t_1}(q_t,u(\cdot))$ and constraints  are
$F_{t_0}(q_t,u(\cdot))=q(0),\ q_t=q(t)$.

Let us set $J_t(u)=J_{t_0}^t(q(t),u(\cdot)),\
\Phi_t(u)=F_{t_0}(q(t),u(\cdot))$. A covector $\lambda\in
T^*M$ is a Lagrange multiplier for the problem $(J_t,\Phi_t)$ if
and only if there exists an extremal $\hat\lambda_\tau,\
t_0\le\tau\le t$, such that $\lambda_{t_0}=\lambda,\
\hat\lambda_t\in T^*_{q(t)}M$. In particular, $\lambda_{t_0}$ is a
Lagrange multiplier for the problem $(J_t,\Phi_t)$ associated to the
control $u(\cdot)=\bar u(\lambda_.)$.

Assume that $\mathrm{ind}\,\mathrm{Hess}_u
\left(J_{t_1}\bigr|_{\Phi^{-1}_{t_1}(q(t_0))}\right)\le\infty,\
t_0\le t\le t_1$ and set $\bar\Phi_t=(J_t,\Phi_t)$.
The curve
$$
t\mapsto \mathcal L_{(\lambda_{t_0},u)}(\bar\Phi_t),\quad t_0\le
t\le t_1
$$
in the Lagrange Grassmannian
$L\left(T_{\lambda_{t_0}}(T^*M)\right)$ is called the Jacobi curve
associated to the extremal $\lambda_t,\ t_0\le t\le t_1$.

In general, the Jacobi curve
$t\mapsto \mathcal L_{(\lambda_{t_0},u)}(\bar\Phi_t)$ is not
smooth, it may even be discontinues, but it is monotone decreasing
in a sense we are going to briefly describe now. You can find more
details in \cite{syn} (just keep in mind that similar quantities
may have opposite signs in different papers; sign agreements vary
from paper to paper that is usual for symplectic geometry).
Monotone curves in the Lagrange Grassmannian have analytic properties
similar to scalar monotone functions: no more than a countable set
of discontinuity points, right and left limits at every point, and
differentiability almost everywhere with semi-definite derivatives
(nonnegative for monotone increasing curves and nonpositive for
decreasing ones). True reason for such a monotonicity is a natural
monotonicity of the family $\bar\Phi_t$. Indeed, let $\tau<t$,
then $\bar\Phi_\tau$ is, in fact, the restriction of $\bar\Phi_t$
to certain subspace: $\bar\Phi_\tau=\bar\Phi_t\circ\mathfrak
p_\tau$, where $\mathfrak p_\tau(u)(s)=
\left\{\begin{array}{rcl}u(s)&,&s<\tau\\ \tilde u(s)&,&
s>\tau\end{array}\right.$.
One can define the Maslov index of a (maybe discontinues) monotone
curve in the Lagrange Grassmannian and the relation between the
Morse and Maslov index indices from Theorem~I.3 remains true.

In fact, Maslov index is a key tool in the whole construction. The
starting point is the notion of a {\it simple curve}. A smooth
curve $\Lambda(\tau),\ \tau_0\le\tau\le\tau_1,$ in the Lagrange
Grassmannian $L(\Sigma)$ is called simple if there exists
$\Delta\in L(\Sigma)$ such that $\Delta\cap\Lambda(\tau)=0,\
\forall\tau\in[\tau_0,\tau_1]$; in other words, the entire curve
is contained in one coordinate chart. It is not hard to show that
any two points of $L(\Sigma)$ can be connected by a simple
monotone increasing (as well as monotone decreasing) curve. An
important fact is that the Maslov index $\mu(\Lambda_\Pi(\cdot))$ of a
simple monotone increasing curve $\Lambda(\tau),\
\tau_0\le\tau\le\tau_1$ is uniquely determined by the triple
$(\Pi,\Lambda(\tau_0),\Lambda(\tau_1))$; i.e. it has the same
value for all simple monotone increasing curves connecting
$\Lambda(\tau_0)$ with $\Lambda(\tau_1)$. A simple way to see this
is to find an intrinsic algebraic expression for the Maslov index
preliminary computed for some simple monotone curve in some
coordinates. We can use Lemma~I.2 for this computation since the
curve is simple. The monotonic increase of the curve implies that
$S_{\Lambda(t_1)}>S_{\Lambda(t_0)}$.

\medskip\noindent {\bf Exercise.} Let $S_0,S_1$ be nondegenerate
symmetric matrices and $S_1\ge S_0$. Then
$
\mathrm{ind}S_0-\mathrm{ind}S_1=\mathrm{ind}\left(S^{-1}_0-S^{-1}_1\right).
$

\medskip Let
$x\in\left(\Lambda(\tau_0)+\Lambda(\tau_1\right)\cap\Pi$ so that
$x=x_0+x_1$, where $x_i\in\Lambda(\tau_i),\ i=0,1$. We set
$\mathfrak q(x)=\sigma(x_1,x_0)$. If
$\Lambda(\tau_0)\cap\Lambda(\tau_1)=0$, then
$\Lambda(\tau_0)+\Lambda(\tau_1)=\Sigma$, $x$ is any element of
$\Pi$ and $x_0,x_1$ are uniquely determined by $x$. This is not
true if $\Lambda(\tau_0)\cap\Lambda(\tau_1)\ne 0$ but $\mathfrak
q(x)$ is well-defined anyway: $\sigma(x_1,x_2)$ depends only on
$x_0+x_1$ since $\sigma$ vanishes on $\Lambda(\tau_i),\ i=0,1.$

Now we compute $\mathfrak q$ in coordinates. Recall that
$$
\Lambda(\tau_i)=\{(y^\top, S_{\Lambda(\tau_i)}y):y\mathbb
R^n\},\ i=0,1,\ \Pi=\{y^\top,0):y\in\mathbb R^n\}.
$$ We have
$$\mathfrak
q(x)=y_1^\top S_{\Lambda(\tau_0)}y_0-y_0^\top S_{\Lambda(\tau_1)}y_1,$$
where $x=(y_0^\top+y_1^\top,0)$,
$S_{\Lambda(\tau_0)}y_0+S_{\Lambda(\tau_1)}y_1=0$. Hence
$y_1=-S^{-1}_{\Lambda(tau_1)}S_{\Lambda(\tau_0)}y_0$ and
$$
\mathfrak q(x)=-y_0^\top S_{\Lambda(\tau_0)}y_0-
\left(S_{\Lambda(\tau_0)}y_0\right)^\top S^{-1}_{\Lambda(\tau_1)}
S_{\Lambda(\tau_0)}y_0=y^\top\left(S^{-1}_{\Lambda(\tau_0)}-
S^{-1}_{\Lambda(\tau_1)}\right)y,
$$ where
$y=S_{\Lambda(\tau_0)}y_0$. We see that the form $\mathfrak q$ is
equivalent, up to a linear change of coordinates, to the quadratic
form defined by the matrix
$S^{-1}_{\Lambda(\tau_0)}-S^{-1}_{\Lambda(\tau_1)}$. Now we set
$$
\mathrm{ind}_\Pi(\Lambda(\tau_0),\Lambda(\tau_1))
\stackrel{def}{=}\mathrm{ind}\,\mathfrak
q.
$$
The above exercise and Lemma~I.2 imply the following:
\begin{lemma} If $\Lambda(\tau),\ \tau_0\le\tau\le\tau_1$, is a
simple monotone increasing curve, then
$$
\mu(\Lambda(\cdot))=\mathrm{ind}_\Pi(\Lambda(\tau_0),\Lambda(\tau_1)).
$$
\end{lemma}

Note that definition of the form $\mathfrak q$ does not require
transversality of $\Lambda(\tau_i)$ to $\Pi$. It is convenient to
extend definition of
$\mathrm{ind}_\Pi(\Lambda(\tau_0),\Lambda(\tau_1))$ to this case.
General definition is as follows:
$$
\mathrm{ind}_\Pi(\Lambda_0,\Lambda_1)=\mathrm{ind}\,\mathfrak
q+\frac 12(\dim(\Pi\cap\Lambda_0)+\dim(\Pi\cap\Lambda_1))-
\dim(\Pi\cap\Lambda_0\cap\Lambda_1).
$$
The Maslov index also has appropriate extension (see
\cite[Sec.4]{sy}) and Lemma II.8 remains true.

Index $\mathrm{ind}_\Pi(\Lambda_0,\Lambda_1)$ satisfies the
triangle inequality:
$$
\mathrm{ind}_\Pi(\Lambda_0,\Lambda_2)\le
\mathrm{ind}_\Pi(\Lambda_0,\Lambda_1)+\mathrm{ind}_\Pi(\Lambda_1,\Lambda_2).
$$
Indeed, the right-hand side of the inequality is equal to the
Maslov index of a monotone increasing curve connecting $\Lambda_0$
with $\Lambda_2$, i.e. of the concatenation of two simple monotone increasing
curves. Obviously, the Maslov index of a simple monotone
increasing curve is not greater than the Maslov index of any other
monotone increasing curve connecting the same endpoints.

The constructed index gives a nice presentation of the Maslov
index of any (not necessary simple) monotone increasing curve
$\Lambda(t),\ t_0\le t\le t_1$:
$$
\mu_\Pi (\Lambda(\cdot))=\sum\limits_{i=0}^l
\mathrm{ind}_\Pi(\Lambda(\tau_i),\Lambda(\tau_{i+1})), \eqno (22)
$$
where $t_0=\tau_0<\tau_1<\cdots<\tau_l<\tau_{l+1}=t_1$ and
$\Lambda\bigr|_{[\tau_i,\tau_{i+1}]}$ are simple pieces of the
curve $\Lambda(\cdot)$. If the pieces are not simple, then the
right-hand side of (22) gives a low bound for the Maslov index
(due to the triangle inequality).

Let now $\Lambda(t),\ t_0\le t\le t_1,$ be a smooth curve which
is {\sl not} monotone increasing. Take any subdivision
$t_0=\tau_0<\tau_1<\cdots<\tau_l<\tau_{l+1}=t_1$ and compute the
sum $\sum\limits_{i=0}^l
\mathrm{ind}_\Pi(\Lambda(\tau_i),\Lambda(\tau_{i+1}))$. This sum
inevitably goes to infinity when the subdivision becomes finer and
finer. The reason is as follows:
$\mathrm{ind}_\Pi(\Lambda(\tau_i),\Lambda(\tau_{i+1}))>0$ for any
simple piece $\Lambda\bigr|_{[\tau_i,\tau_{i+1}]}$ such that
$\underline{\dot\Lambda}(\tau)\ngeq 0,\
\forall\tau\in[\tau_i,\tau_{i+1}]$ and
$\mu_\Pi(\Lambda\bigr|_{[\tau_i,\tau_{i+1}]}=0$. I advise reader
to play with the one-dimensional case of the curve in $L(\mathbb
R^2)=S^1$ to see better what's going on.

This should now be clear how to manage in the general nonsmooth
case. Take a curve $\Lambda(\cdot)$ (an arbitrary mapping from
$[t_0,t_1]$ into $L(\Sigma)$). For any finite subset
$\mathcal T+\{\tau_1,\ldots,\tau_k\}\subset[t_0,t_1]$, where
$t_0=\tau_0<\tau_1<\cdots<\tau_l<\tau_{l+1}=t_1$, we compute the
sum $I_\Pi^{\mathcal T}=\sum\limits_{i=0}^l
\mathrm{ind}_\Pi(\Lambda(\tau_i),\Lambda(\tau_{i+1}))$ and then
find supremum of these sums for all finite subsets:
$I_\Pi(\Lambda(\cdot))=\sup\limits_{\mathcal T}I_\Pi^{\mathcal
T}$. The curve $\Lambda(\cdot)$ is called monotone increasing if
$I_\Pi^{\mathcal T}<\infty$; it is not hard to show that the last
property does not depend on $\Pi$ and that monotone increased
curves enjoy listed above analytic properties. A curve
$\Lambda(\cdot)$ is called monotone decreasing if inversion of the
parameter $t\mapsto t_0+t_1-t$ makes it monotone increasing.

We set $\mu(\Lambda(\cdot))=I_\Pi(\Lambda(\cdot))$ for any
monotone increasing curve and \linebreak $\mu(\Lambda(\cdot))=
-I_\Pi(\hat\Lambda(\cdot))$ for a monotone decreasing one, where
$\hat\Lambda(t)=\Lambda(t_0+t_1-t)$. The defined in this way
Maslov index of a discontinues monotone curve equals the Maslov
index of the continues curve obtained by gluing all
discontinuities with simple monotone curves of the same direction
of monotonicity.

If $\Lambda(t)=\mathcal L_{(\lambda_{t_0},u)}(\bar\Phi_t)$ is the
Jacobi curve associated to the extremal with a finite Morse index,
then $\Lambda(\cdot)$ is monotone decreasing and its Maslov index
computes
$\mathrm{ind}\,\mathrm{Hess}_u\left(J_{t_1}\bigr|_{\Phi^{-1}_{t_1}(q(t_0))}\right)$
in the way similar to Theorem~I.3. Of course, these nice things
have some value only if we can effectively find Jacobi curves for
singular extremals: their definition was too abstract.
Fortunately, this is not so hard; see \cite{a} for the
explicit expression of Jacobi curves for a wide class of singular
extremals and, in particular, for singular curves of rank 2 vector
distributions (these last Jacobi curves have found important
applications in the geometry of distributions, see \cite{dz,z}).

\medskip One more important property of monotonic curves is as
follows.
\begin{lemma} Assume that $\Lambda(\cdot)$ is monotone and
right-continues at $t_0$, i.e.
$\Lambda(t_0)=\lim\limits_{t\searrow t_0}\Lambda(t)$. Then
$\Lambda(t_0)\cap\Lambda(t)=\bigcap\limits_{t_0\le\tau\le
t}\Lambda(t)$ for any $t$ sufficiently close to (and greater than)
$t_0$.
\end{lemma}
{\bf Proof.} We may assume that $\Lambda(\cdot)$ is monotone
increasing. Take centered at $\Lambda(t_0)$ local coordinates in
the Lagrange Grassmannian; the coordinate presentation of
$\Lambda(t)$ is a symmetric matrix $S_{\Lambda(t)}$, where
$S_{\Lambda(t_0)}=0$ and $t\mapsto y^\top S_{\Lambda(t)}y$ is a
monotone increasing scalar function $\forall y\in\mathbb R^n$. In
particular, $\ker S_{\Lambda(t)}=\Lambda(t)\cap\Lambda(t_0)$ is a
monotone decreasing family of subspaces. \quad $\square$

We set $\Gamma_t=\bigcap\limits_{t_0\le\tau\le t}\Lambda(\tau)$,
a monotone decreasing family of isotropic subspaces. Let
$\Gamma=\max\limits_{t>t_0}\Gamma_t$, then $\Gamma_t==\Gamma$ for
all $t>t_0$ sufficiently close to $t_0$. We have:
  $\Lambda(t)=\Lambda(t)^\angle$ and $\Lambda(t)\supset\Gamma$
for all $t>t_0$ close enough to $t_0$;
hence $\Gamma_t^\angle\supset\Lambda(t)$. In particular,
$\Lambda(t)$ can be treated as a Lagragian subspace of the
symplectic space $\Gamma^\angle/\Gamma$.
Moreover, Lemma~II.9 implies that
$\Lambda(t)\cap\Lambda(t_0)=\Gamma$. In other words, $\Lambda(t)$
is transversal to $\Lambda(t_0)$ in $\Gamma^\angle/\Gamma$. In the
case of a real-analytic monotone curve $\Lambda(\cdot)$ this
automatically implies that $\Lambda(\cdot)$ is an ample curve
in $\Gamma^\angle/\Gamma$. Hence any nonconstant monotone
analytic curve is reduced to an ample monotone curve. It becoms ample
after the factorization by a fixed (motionless) subspace.

\section{Comparizon theorem}

We come back to smooth regular curves after the deviation devoted
to a more general perspective.

\begin{lemma} Let $\Lambda(t),\ t\in[t_0,t_1]$ be a regular monotone
increasing curve in the Lagrange Grassmannian $L(\Sigma)$. Then
$\{t\in[t_0,t_1]:\Lambda(t)\cap\Pi\ne 0\}$ is a finite subset of
$[t_0,t_1]\ \forall \Pi\in L(\Sigma)$. If $t_0$ and $t_1$ are out
of this subset, then
$$
\mu_\Pi(\Lambda(\cdot))=\sum\limits_{t\in(t_0,t_1)}\dim(\Lambda(t)\cap\Pi).
$$
\end{lemma}
{\bf Proof.} We have to proof that $\Lambda(t)$ may have a nontrivial
intersection with $\Pi$ only for isolated values of $t$; the rest
is Lemma~I.1. Assume that $\Lambda(t)\cap\Pi\ne 0$ and take a centered
at $\Pi$ coordinate neighborhood in $L(\Sigma)$ which contains
$\Lambda(t)$. In these coordinates, $\Lambda(\tau)$ is presented by
a symmetric matrix $S_\Lambda(\tau)$ for any $\tau$ sufficiently
close to $t$ and $\Lambda(\tau)\cap\Pi=\ker S_{\Lambda(\tau)}.$
Monotonicity and regularity properties are equivalent to the
inequality $\dot S_{\Lambda(\tau)}>0$. In particular,
$y^\top\dot S_{\Lambda(t)}y>0\ \forall y\in\ker S_{\Lambda(t)}
\setminus\{0\}$. The last inequality implies that $ S_{\Lambda(\tau)}$
is a nondegenerate for all $\tau$ sufficiently close and not equal
to $t$.

\medskip\noindent{\bf Definition.} Parameter values $\tau_0,\tau_1$
are called conjugate for the continues curve $\Lambda(\cdot)$ in the
Lagrange Grassmannian if $\Lambda(\tau_0)\cap\Lambda(\tau_1)\ne 0$;
the dimension of $\Lambda(\tau_0)\cap\Lambda(\tau_1)$ is the
 {\it multiplicity} of the conjugate parameters.

\medskip If $\Lambda(\cdot)$ is a regular monotone increasing curve,
then, according to Lemma II.9, conjugate points are isolated and
the Maslov index
$\mu_{\Lambda(t_0)}\left(\Lambda\bigr|_{[t,t_1]}\right)$ equals the
sum of multiplicities of the conjugate to $t_0$ parameter values located
in $(t,t_1)$.
If $\Lambda(\cdot)$ is the Jacobi curve of an extremal of an
optimal control problem, then this Maslov index equals the Morse
index of the extremal; this is why conjugate points are so important.

Given a regular monotone curve $\Lambda(\cdot)$, the quadratic
form $\underline{\dot\Lambda}(t)$ defines an Euclidean structure
$\langle\cdot,\cdot\rangle_{\dot\Lambda(t)}$ on $\Lambda(t)$ so
that
$\underline{\dot\Lambda}(t)(x)=\langle x,x\rangle_{\dot\Lambda(t)}$.
Let $R_\Lambda(t)\in\mathrm{gl}(\Lambda(t))$ be the curvature
operator of the curve $\Lambda(\cdot)$; we define the {\it curvature
quadratic form}
$r_\lambda(t)$ on $\Lambda(t)$ by the formula:
$$
r_\Lambda(t)(x)=\langle R_\Lambda(t)x,x\rangle_{\dot\Lambda(t)},\quad
x\in\Lambda(t).
$$

\begin{prop} The curvature operator $R_\Lambda(t)$ is a
self-adjoint operator for the Euclidean structure
$\langle\cdot,\cdot\rangle_{\dot\Lambda(t)}$. The form
$r_\Lambda(t)$ is equivalent (up to linear changes of variables)
to the form $\underline{\dot\Lambda}^\circ(t)$, where
$\Lambda^\circ(\cdot)$ is the derivative curve.
\end{prop}
{\bf Proof.} The statement is intrinsic and we may check it in any
coordinates. Fix $t$ and take Darboux coordinates
$\{(\eta,y):\eta\in\mathbb R^{n*},y\in\mathbb R^n\}$ in $\Sigma$
in such a way that $\Lambda(t)=\{(y^\top,0):y\in\mathbb R^n\}$,
$\Lambda^\circ(t)=\{(0,y):y\in\mathbb R^n\}$,
$\underline{\dot\Lambda}(t)(y)=y^\top y$. Let
$\Lambda(\tau)=\{(y^\top,S_\tau y):y\in\mathbb R^n\}$, then
$S_t=0$. Moreover, $\dot S(t)$ is the matrix of the form
$\underline{\dot\Lambda}(t)$ in given coordinates, hence $\dot
S_t=I$. Recall that $\Lambda^\circ(\tau)=\{(y^\top A_\tau,y+S_\tau
A_\tau y):y\in\mathbb R^n\}$, where $A_\tau=-\frac 12\dot
S^{-1}_\tau\ddot S_\tau\dot S^{-1}_\tau$ (see (5)). Hence $\ddot
S_t=0$. We have: $R_\Lambda(t)=\frac 12\stackrel{\ldots}{S}_t$,
$r_\Lambda(t)(y)=\frac 12y^\top\stackrel{\ldots}{S}_ty$,
$$
\underline{\dot\Lambda}^\circ(t)(y)=\sigma\left((0,y),(y^\top\dot
A_t,0)\right)=-y^\top\dot A_ty=\frac
12y^\top\stackrel{\ldots}{S}_ty.
$$
So $r_\Lambda(t)$ and $\underline{\dot\Lambda}^\circ(t)$ have
equal matrices for our choice of coordinates in $\Lambda(t)$
and $\Lambda^\circ(t)$. The curvature operator is self-adjoint
since it is presented by a symmetric matrix in coordinates where
form $\underline{\dot\Lambda}(t)$ is the standard inner product.
$\quad \square$

Proposition II.9 implies that the curvature operators of regular
monotone curves in the Lagrange Grassmannian are diagonalizable
and have only real eigenvalues.

\begin{theorem} Let $\Lambda(\cdot)$ be a regular monotone curve
in the Lagrange Grassmannian $L(\Sigma)$, where $\dim\Sigma=2n$.
\begin{itemize}\item
If all eigenvalues of $R_\Lambda(t)$ do not exceed a constant
$c\ge 0$ for any $t$ from the domain of $\Lambda(\cdot)$, then
$|\tau_1-\tau_0|\ge\frac\pi{\sqrt c}$ for any pair of conjugate
parameter values $\tau_0,\tau_1$. In particular, If all
eigenvalues of $R_\Lambda(t)$ are nonpositive $\forall t$, then
$\Lambda(\cdot)$ does not possess conjugate parameter values.
\item If $\mathrm{tr}R_\Lambda(t)\ge nc$ for some constant $c>0$
and $\forall t$, then, for arbitrary $\tau_0\le t$, the segment
$[t,t+\frac\pi{\sqrt c}]$ contains a conjugate to $\tau_0$
parameter value as soon as this segment is contained in the domain
of $\Lambda(\cdot)$.
\end{itemize}
Both estimates are sharp.
\end{theorem}
{\bf Proof.} We may assume without lack of generality that
$\Lambda(\cdot)$ is ample monotone increasing. We start with the
case of nonpositive eigenvalues of $R_\Lambda(t)$. The absence of
conjugate points follows from Proposition II.9 and the following
\begin{lemma} Assume that $\Lambda(\cdot)$ is an ample monotone
increasing (decreasing) curve and $\Lambda^\circ(\cdot)$ is a
continues monotone decreasing (increasing) curve. Then
$\Lambda(\cdot)$ does not possess conjugate parameter values and
there exists a $\lim\limits_{t\to+\infty}\Lambda(t)=\Lambda_\infty$.
\end{lemma}
{\bf Proof.} Take some value of the parameter $\tau_0$; then
$\Lambda(\tau_0)$ and $\Lambda^\circ(\tau_0)$ is a pair of
transversal Lagrangian subspaces. We may choose coordinates in the
Lagrange Grassmannian in such a way that $S_{\Lambda(\tau))}=0$
and $S_{\Lambda^\circ(\tau_0)}=I$, i.e. $\Lambda(\tau_0)$ is
represented by zero $n\times n$-matrix and $\Lambda^\circ(\tau_0)$
by the unit matrix. Monotonicity assumption implies that $t\mapsto
S_{\Lambda(t)}$ is a monotone increasing curve in the space
of symmetric matrices and $t\mapsto S_{\Lambda^\circ(t)}$ is a
monotone decreasing curve. Moreover, transversality of
$\Lambda(t)$ and $\Lambda^\circ(t)$ implies that
$S_{\Lambda^\circ(t)}-S_{\Lambda(t)}$ is a nondegenerate matrix.
Hence $0<S_{\Lambda(t)}<S_{\Lambda^\circ(t)}\le I$ for any
$t>\tau_0$. In particular, $\Lambda(t)$ never leaves the
coordinate neighborhood under consideration for $T>\tau_0$, the
subspace $\Lambda(t)$ is always transversal to $\Lambda(\tau_0)$
and has a limit $\Lambda_\infty$, where
$S_{\Lambda_\infty}=\sup\limits_{t\ge\tau_0}S_{\Lambda(t)}. \qquad \square$

Now assume that the eigenvalues of $R_{\Lambda}(t)$ do not exceed
a constant $c>0$. We are going to reparametrize the the curve
$\Lambda(\cdot)$ and to use the chain rule (7). Take some $\bar t$
in the domain of $\Lambda(\cdot)$ and set
$$
\varphi(t)=\frac 1{\sqrt c}\left(\mathrm{arctan}(\sqrt
ct)+\frac\pi{2}\right)+\bar t, \quad
\Lambda_\varphi(t)=\Lambda(\varphi(t)).
$$
We have: $\varphi(\mathbb R)=\left(\bar t,\bar t+\frac\pi{\sqrt
c}\right)$, $\dot\varphi(t)=\frac 1{ct^2+1}$,
$R_{\varphi}(t)=-\frac c{(ct^2+1)^2}$. Hence, according to the
chain rule (7), the operator
$$
R_{\Lambda_\varphi}(t)=\frac
1{(ct^2+1)^2}\left(R_\Lambda(\varphi(t))-cI\right)
$$
has only nonpositive eigenvalues. Already proved part of the
theorem implies that $\Lambda_\varphi$ does not possess conjugate
values of the parameter. In other words, any length
$\frac\pi{\sqrt c}$ interval in the domain of $\Lambda(\cdot)$ is
free of conjugate pairs of the parameter values.

Assume now that $\mathrm{tr}R_\Lambda(t)\ge nc$. We will prove
that the existence of $\Delta\in L(\Sigma)$ such that
$\Delta\cap\Lambda(t)=0$ for all $t\in[\bar t,\tau]$ implies that
$\tau-\bar t<\frac\pi{\sqrt c}$. We'll prove it by contradiction.
If there exists such a $\Delta$, then $\Lambda\bigr|_{[\bar
t,\tau]}$ is completely contained in a fixed coordinate
neighborhood of $L(\Sigma)$, therefore the curvature operator
$R_\Lambda(t)$ is
defined by the formula (6). Put $B(t)=(2\dot S_t)^{-1}\ddot S_t$,
$b(t)=\mathrm{tr}B(t)$, $t\in[\bar t,\tau]$. Then
$$
\dot B(t)=B^2(t)+R_\Lambda(t), \quad \dot
b(t)=\mathrm{tr}B^2(t)+\mathrm{tr}R_\Lambda(t).
$$
Since for an arbitrary symmetric $n\times n$-matrix $A$ we have
$\mathrm{tr}A^2\ge\frac 1n(\mathrm{tr}A)^2$, the inequality $\dot
b\ge\frac{b^2}n+nc$ holds. Hence $b(t)\ge\beta(t),\ \bar t\le
t\le\tau$, where $\beta(\cdot)$ is a solution of the equation
$\dot\beta=\frac{\beta^2}n+nc$, i.e. $\beta(t)=n\sqrt c\tan(\sqrt
c(t-t_0))$. The function $b(\cdot)$ together with $\beta(\cdot)$
are bounded on the segment $[\bar t,\tau]$. Hence
$\tau-t\le\frac\pi{\sqrt c}$.

To verify that the estimates are sharp, it is enough to consider
regular monotone curves of constant curvature. $\quad \square$

\section{Reduction}

We consider a Hamiltonian system on a symplectic manifold $N$
endowed with a fixed Lagrange foliation $E$. Assume that
$g:N\to\mathbb R$ is a first integral of our Hamiltonian system,
i.e. $\{h,g\}=0$.

\begin{lemma} Let $z\in N,\ g(z)=c$. The leaf $E_z$ is transversal
to $g^{-1}(c)$ at $z$ if and only if $\vec g(z)\notin T_zE_z$.
\end{lemma}
{\bf Proof.} Hypersurface $g^{-1}(c)$ is not transversal to
 $g^{-1}(c)$ at $z$ if and only if
 $$
 d_zg(T_zE_z)=0\ \Leftrightarrow\ \sigma(\vec g(z),T_zE_z)=0\
 \Leftrightarrow\ \vec g(z)\in(T_zE_z)^\angle=T_zE_z. \eqno\square
$$

 If all points of some level $g^{-1}(c)$ satisfy conditions of
 Lemma II.12, then $g^{-1}(c)$ is a (2n-1)-dimensional manifold
 foliated by $(n-1)$-dimensional submanifolds $E_z\cap g^{-1}(c)$.
 Note that $\mathbb R\vec g(z)=\ker\sigma\bigr|_{T_zg^{-1}(c)}$,
 hence \linebreak $\Sigma_z^g\stackrel{def}{=}T_zg^{-1}(c)/\mathbb R\vec
 g(z)$ is a $2(n-1)$-dimensional symplectic space and \linebreak
 $\Delta_z^g\stackrel{def}{=}T_z\left(E_z\cap g^{-1}(c)\right)$ is
 a Lagrangian subspace in $L_z^g$, i.e. $\Delta^g_z\in
 L(\Sigma^g_z)$.

The submanifold $g^{-1}(c)$ is invariant for the flow $e^{t\vec
h}$. Moreover, $e^{t\vec h}_*\vec g=\vec g$. Hence $e^{t\vec h}_*$
induces a symplectic transformation
$e^{t\vec h}_*:\Sigma^g_z\to\Sigma^g_{e^{t\vec h}(z)}$. Set
$J_z^g(t)=e^{-t\vec h}_*\Delta^g_{e^{t\vec h}(z)}$. The curve
$t\mapsto J^g_z(t)$ in the Lagrange Grassmannian $L(\Sigma^g_z)$
is called a {\it reduced Jacobi curve} for the Hamiltonian field
$\vec h$ at $z\in N$.

The reduced Jacobi curve can be easily reconstructed from the
Jacobi curve
$J_z(t)=e_*^{-t\vec h}\left(T_{e^{t\vec h}(z)}E_{e^{t\vec
h}(z)}\right)\in L(T_zN)$ and vector $\vec g(z)$. An elementary
calculation shows that
$$
J^g_z(t)=J_z(t)\cap\vec g(z)^\angle+\mathbb R\vec g(z).
$$
Now we can temporary forget the symplectic manifold and
Hamiltonians and formulate everything in terms of the curves in
the Lagrange Grassmannian. So let $\Lambda(\cdot)$ be a smooth
curve in the Lagrange Grassmannian $L(\Sigma)$ and $\gamma$ a
one-dimensional subspace in $\Sigma$. We set
$\Lambda^\gamma(t)=\Lambda(t)\cap\gamma^\angle+\gamma$, a Lagrange
subspace in the symplectic space $\gamma^\angle/\gamma$. If
$\gamma\not\subset\Lambda(t)$, then $\Lambda^\gamma(\cdot)$ is
smooth and
$\underline{\dot\Lambda}^\gamma(t)=
\underline{\dot\Lambda}(t)\bigr|_{\Lambda(t)\cap\gamma^\angle}$ as
it easily follows from the definitions. In particular,
monotonicity of $\Lambda(\cdot)$ implies monotonicity of
$\Lambda^\gamma(\cdot)$; if $\Lambda(\cdot)$ is regular and
monotone, then $\Lambda^\gamma(\cdot)$ is also regular and
monotone. The curvatures and the Maslov indices of
$\Lambda(\cdot)$ and $\Lambda^\gamma(\cdot)$ are related in a more
complicated way. The following result is proved in \cite{acz}.

\begin{theorem} Let $\Lambda(t),\ t\in[t_0,t_1]$ be a smooth
monotone increasing curve in $L(\Sigma)$ and $\gamma$ a
one-dimensional subspace of $\Sigma$ such that
$\gamma\not\subset\Lambda(t),\ \forall t\in[t_0,t_1]$. Let
$\Pi\in L(\Sigma),\ \gamma\not\subset\Pi,\
\Lambda(t_0)\cap\Pi=\Lambda(t_1)\cap\Pi=0$. Then
\begin{itemize}\item
$\mu_{\Pi}(\Lambda(\cdot))\le\mu_{\Pi^\gamma}(\Lambda^\gamma(\cdot))
\le \mu_{\Pi}(\Lambda(\cdot))+1.$
\item If $\Lambda(\cdot)$ is regular, then
$r_{\Lambda^\gamma}(t)\ge
r_\Lambda(t)\bigr|_{\Lambda(t)\cap\gamma^\angle}$ and \\
$\mathrm{rank}\left(r_{\Lambda^\gamma}(t)-
r_\Lambda(t)\bigr|_{\Lambda(t)\cap\gamma^\angle}
\right)\le 1. \qquad $
\end{itemize}
\end{theorem}

The inequality $r_{\Lambda^\gamma}(t)\ge
r_\Lambda(t)\bigr|_{\Lambda(t)\cap\gamma^\angle}
$ turns into the equality if \linebreak
$\gamma\subset\Lambda^\circ(t),\ \forall t$. Then
$\gamma\subset\ker\underline{\dot\Lambda}^\circ(t)$. According to
Proposition II.9, to $\gamma$ there corresponds a one-dimensional
subspace in the kernel of $r_\Lambda(t)$; in particular, $r_\Lambda(t)$
is degenerate.

Return to the Jacobi curves $J_z(t)$ of a Hamiltonian field $\vec
h$. There always exists at least one first integral: the
Hamiltonian $h$ itself. In general, $\vec h(z)\notin J^\circ_z(0)$
and the reduction procedure has a nontrivial influence on the
curvature (see \cite{ac,acz} for explicit expressions). Still,
there is an important class of Hamiltonians and Lagrange
foliations for which the relation $\vec h(z)\in J^{\circ}_z(0)$
holds $\forall z$. These are homogeneous on fibers Hamiltonians on
cotangent bundles. In this case the generating homotheties of the
fibers Euler vector field belongs to the kernel of the curvature
form.

\section{Hyperbolicity}

\noindent{\bf Definition.} We say that a Hamiltonian function $h$
on the symplectic manifold $N$ is regular with respect to the
Lagrange foliation $E$ if the functions $h\bigr|_{E_z}$ have
nondegenerate second derivatives at $z,\ \forall z\in N$ (second
derivative is well-defined due to the canonical affine structure
on $E_z$). We say that $h$ is monotone with respect to $E$ if
$h\bigr|_{E_z}$ is a convex or concave function $\forall z\in N$.

\medskip Typical examples of regular monotone Hamiltonians on the
cotangent bundles are energy functions of natural mechanical
systems. Such a function is the sum of the kinetic energy whose
Hamiltonian system generates the Riemannian geodesic flow and a
``potential" that is a constant on the fibers function.
Proposition II.8 implies that Jacobi curves associated to the
regular monotone Hamiltonians are also regular and monotone. We'll
show that negativity of the curvature operators of such a
Hamiltonian implies the hyperbolic behavior of the Hamiltonian
flow. This is a natural extension of the classical result about
Riemannian geodesic flows.

Main tool is the structural equation derived in Section 13. First
we'll show that this equation is well coordinated with the
symplectic structure. Let $\Lambda(t),\ t\in\mathbb R,$ be a
regular curve in $L(\Sigma)$ and
$\Sigma=\Lambda(t)\oplus\Lambda^\circ(t)$ the correspondent
canonical splitting. Consider the structural equation
$$
\ddot e(t)+R_\Lambda(t)e(t)=0,\quad \mathrm{where}\ e(t)\in\Lambda(t),
\ \dot e(t)\in\Lambda^\circ(t), \eqno (23)
$$
(see Corollary II.1).

\begin{lemma} The mapping
$e(0)\oplus\dot e(0)\mapsto e(t)\oplus\dot e(t)$, where $e(\cdot)$
and $\dot e(\cdot)$ satisfies (23), is a symplectic transformation of $\Sigma$.
\end{lemma}
{\bf Proof.} We have to check that $\sigma(e_1(t),e_2(t)),\
\sigma(\dot e_1(t),\dot e_2(t)),\ \sigma(e_1(t),\dot e_2(t))$ do
not depend on $t$ as soon as $e_i(t), \dot e_i(t),\ i=1,2$,
satisfy (23). First two quantities vanish since $\Lambda(t)$ and
$\Lambda^\circ(t)$ are Lagrangian subspaces. The derivative of the
third quantity vanishes as well since $\ddot e_i(t)\in\Lambda(t). \quad
\square $

Let $h$ be a regular monotone Hamiltonian on the symplectic
manifold $N$ equipped with a Lagrange foliation $E$. As before, we
denote by $J_z(t)$ the Jacobi curves of $\vec h$ and by $J_z^h(t)$
the reduced to the level of $h$ Jacobi curves (see previous
Section). Let $R(z)=R_{J_z}(0)$ and $R^h(z)=R_{J_z^h}(0)$ be the
curvature operators of $J_z(\cdot)$ and $J_z^h(\cdot)$
correspondently. We say that the Hamiltonian field $\vec h$ has a
negative curvature at $z$ with respect to $E$ if all eigenvalues of
$R(z)$ are negative. We say that $\vec h$ has a negative reduced
curvature at $z$ if all eigenvalues of $R_z^h$ are negative.

\begin{prop} Let $z_0\in N,\ z_t=e^{t\vec h}(z)$. Assume that
that $\overline{\{z_t:t\in\mathbb R\}}$ is a compact subset of $N$ and that
$N$ is endowed with a Riemannian structure. If $\vec h$ has a
negative curvature at any $z\in\overline{\{z_t:t\in\mathbb R\}}$,
then there exists a constant $\alpha>0$ and a splitting
$T_{z_t}N=\Delta^+_{z_t}\oplus\Delta^-_{z_t}$, where
$\Delta^\pm_{z_t}$ are Lagrangian subspaces of $T_{z_t}N$ such
that
$e^{\tau\vec h}_*(\Delta^\pm_{z_t})=\Delta^\pm_{z_{t+\tau}}\
\forall\,t,\tau\in\mathbb R$ and
$$
\|e^{\pm\tau\vec
h}_*\zeta_\pm\|\ge e^{\alpha\tau}\|\zeta_{\pm}\| \quad
\forall\,\zeta\in\Delta^\pm_{z_t},\,\tau\ge 0,\,t\in\mathbb R.
\eqno (24)
$$
Similarly, if $\vec h$ has a
negative reduced curvature at any $z\in\overline{\{z_t:t\in\mathbb R\}}$,
then there exists a splitting
$T_{z_t}(h^{-1}(c)/\mathbb R h(z_t))=\hat\Delta^+_{z_t}\oplus\hat\Delta^-_{z_t}$, where
$c=h(z_0)$ and $\hat\Delta^\pm_{z_t}$ are Lagrangian subspaces of
$T_{z_t}(h^{-1}(c)/\mathbb R h(z_t))$ such
that
$e^{\tau\vec h}_*(\hat\Delta^\pm_{z_t})=\hat\Delta^\pm_{z_{t+\tau}}\
\forall\,t,\tau\in\mathbb R$ and $\|e^{\pm\tau\vec
h}_*\zeta_\pm\|\ge e^{\alpha\tau}\|\zeta_{\pm}\| \quad
\forall\,\zeta\in\hat\Delta^\pm_{z_t},\,\tau\ge 0,\,t\in\mathbb R.$
\end{prop}
{\bf Proof.} Obviously, the desired properties of
$\Delta^\pm_{z_t}$ and $\hat\Delta^\pm_{z_t}$ do not depend on the
choice of the Riemannian structure on $N$. We'll introduce a
special Riemannian structure determined by $h$. The Riemannian
structure is a smooth family of inner products
$\langle\cdot,\cdot\rangle_z$ on $T_zN$, $z\in N$. We have
$T_zN=J_z(0)\oplus J^\circ_z(0)$, where $J_z(0)=T_zE_z$. Replacing
$h$ with $-h$ if necessary we may assume that $h\bigr|_{E_z}$ is a
strongly convex function. First we define
$\langle\cdot,\cdot\rangle_z\bigr|_{J_z(0)}$ to be equal to the
second derivative of $h\bigr|_{E_z}$. Symplectic form $\sigma$
induces a nondegenerate pairing of $J_z(0)$ and $J^\circ_z(0)$. In
particular, for any $\zeta\in J_z(0)$ there exists a unique
$\zeta^\circ\in J_z^\circ(0)$ such that
$\sigma(\zeta^\circ,\cdot)\bigr|_{J_z(0)}=
\langle\zeta,\cdot\rangle_z\bigr|_{J_z(0)}$. There exists a unique
extension of the inner product $\langle\cdot,\cdot\rangle_z$ from
$J_z(0)$ to the whole $T_zN$ with the following properties:
\begin{itemize}\item $J_z^\circ(0)$ is orthogonal to $J_z(0)$ with
respect to $\langle\cdot,\cdot\rangle_z$;
\item $\langle\zeta_1,\zeta_2\rangle_z=
\langle\zeta^\circ_1,\zeta^\circ_2\rangle_z\
\forall\,\zeta_1,\zeta_2\in J_z(0)$.
\end{itemize}

We'll need the following classical fact from Hyperbolic Dynamics
(see, for instance, \cite[Sec.\,17.6]{kh}).

\begin{lemma} Let $A(t),\ t\in\mathbb R$, be a bounded family of
symmetric $n\times n$-matrices whose eigenvalues are all negative
and uniformly separated from 0. Let $\Gamma(t,\tau)$ be the
fundamental matrix of the $2n$-dimensional linear system
$\dot x=-y$, $\dot y=A(t)x$, where $x,y\in\mathbb R^n$, i.e.
$$
\frac\partial{\partial t}\Gamma(t,\tau)=\left(\begin{smallmatrix}0 &
-I\\A & 0 \end{smallmatrix}\right)\Gamma(t,\tau), \quad
\Gamma(\tau,\tau)=\left(\begin{smallmatrix}I & 0\\ 0 &
I\end{smallmatrix}\right). \eqno (25)
$$
Then there exist closed conic neighborhoods
$C^+_\Gamma,\,C^-_\Gamma$, where
$C^+_\Gamma\cap C^-_\Gamma=0$, of some $n$-dimensional subspaces of $\mathbb
R^{2n}$ and a constant $\alpha>0$
such that $$\Gamma(t,\tau)C^+_\Gamma\subset C^+_\Gamma, \quad
|\Gamma(t,\tau)\xi_+|\ge e^{\alpha(\tau-t)}|\xi_+|, \ \forall\,
\xi_+\in C^+_\Gamma,\,t\le\tau,$$ and $$\Gamma(t,\tau)C^-_\Gamma\subset
C^-_\Gamma, \quad
|\Gamma(t,\tau)\xi_-|\ge e^{\alpha(t-\tau)}|\xi_-|, \  \forall\,
\xi_-\in C^-_\Gamma,\,t\ge\tau.$$ The constant $\alpha$ depends only on upper
and lower
bounds of the eigenvalues of $A(t). \qquad \square$
\end{lemma}

\begin{corollary} Let $C^\pm_\Gamma$ be the cones described in
Lemma~II.14; then \linebreak $\Gamma(0,\pm
t)C^\pm_\Gamma\subset\Gamma(0;\pm\tau)C^\pm_\Gamma$ for any
$t\ge\tau\ge 0$ and the subsets \linebreak $K^\pm_\Gamma=\bigcap\limits_{t\ge
0}\Gamma(0,t)C^\pm_\Gamma$ are Lagrangian subspaces of $\mathbb
R^n\times\mathbb R^n$ equipped with the standard symplectic
structure.
\end{corollary}
{\bf Proof.} The relations $\Gamma(\tau,t)C^+_\Gamma\subset
C^+_\Gamma$ and $\Gamma(\tau,t)C^-_\Gamma\subset
C^-_\Gamma$ imply:
$$
\Gamma(0,\pm t)C^\pm_\Gamma=\Gamma(0,\pm\tau)\Gamma(\pm\tau,\pm t)C^\pm_\Gamma
\subset\Gamma(0,\pm\tau)C^\pm_\Gamma.
$$

In what follows we'll study $K^+_\Gamma$; the same arguments work
for $K^-_\Gamma$. Take vectors $\zeta,\zeta'\in K^+_\Gamma$; then
$\zeta=\Gamma(0,t)\zeta_t$ and $\zeta'=\Gamma(0,t)\zeta'_t$ for
any $t\ge 0$ and some $\zeta_t,\zeta'_t\in C^+_\Gamma$. Then,
according to Lemma II.14, $|\zeta_t|\le e^{-\alpha t}|\zeta|,\
|\zeta'_t|\le e^{-\alpha t}|\zeta'|$, i.e. $\zeta_t$ and
$\zeta'_t$ tend to 0 as $t\to+\infty$.
On the other hand,
$$
\sigma(\zeta,\zeta')=\sigma(\Gamma(0,t)\zeta_t,\Gamma(0,t)\zeta'_t)=
\sigma(\zeta_t,\zeta'_t) \quad \forall t\ge 0
$$
since $\Gamma(0,t)$ is a symplectic matrix. Hence
$\sigma(\zeta,\zeta')=0$.

We have shown that $K^+_\Gamma$ is an isotropic subset of $\mathbb
R^n\times\mathbb R^n$. On the other hand, $K^+_\Gamma$ contains an
$n$-dimensional subspace since $C^+_\Gamma$ contains one and
$\Gamma(0,t)$ are invertible linear transformations. Isotropic
$n$-dimensional subspace is equal to its skew-orthogonal
complement, therefore $K^+_\Gamma$ is a Lagrangian subspace.
$\quad \square$

Take now a regular monotone curve $\Lambda(t),\ t\in\mathbb R$ in
the Lagrange Grassmannian $L(\Sigma)$.
We may assume that $\Lambda(\cdot)$ is monotone
increasing, i.e. $\dot\Lambda(t)>0$. Recall that
$\underline{\dot\Lambda}(t)(e(t))=\sigma(e(t),\dot e(t))$, where
$e(\cdot)$ is an arbitrary smooth curve in $\Sigma$ such that
$e(\tau)\in\Lambda(\tau),\ \forall\tau$. Differentiation of the
identity \linebreak $\sigma(e_1(\tau),e_2(\tau))=0$ implies:
$\sigma(e_1(t),\dot e_2(t))=-\sigma(\dot e_1(t),e_2(t))=
\sigma(e_2(t),\dot e_1(t))$ if $e_i(\tau)\in\Lambda(\tau)$,
$\forall\tau$, $i=1,2$. Hence the Euclidean structure
$\langle\cdot,\cdot\rangle_{\dot\Lambda(t)}$ defined by the
quadratic form $\underline{\dot\Lambda}(t)$ reads:
$\langle
e_1(t),e_2(t)\rangle_{\dot\Lambda(t)}=\sigma(e_1(t),\dot e_2(t))$.

Take a basis $e_1(0),\ldots,e_n(0)$ of $\Lambda(0)$ such that the
form $\underline{\dot\Lambda}(t)$ has the unit matrix in this
basis, i.e. $\sigma(e_i(0),\dot e_j(0))=\delta_{ij}$. In fact,
vectors $\dot e_j(0)$ are defined modulo $\Lambda(0)$; we can
normalize them assuming that $\dot e_i(0)\in\Lambda^\circ(0),\
i=1,\ldots,n$. Then $e_1(0),\ldots,e_n(0),\dot e_1(0),\ldots,\dot
e_n(0)$ is a Darboux basis of $\Sigma$. Fix coordinates in
$\Sigma$ using this basis: $\Sigma=\mathbb R^n\times\mathbb R^n$,
where
$\left(\begin{smallmatrix}x\\y\end{smallmatrix}\right)\in\mathbb
R^n\times\mathbb R^n$ is identified with
$\sum\limits_{j=1}^n\left(x^je_j(0)+y^j\dot
e_j(0)\right)\,\in\Sigma,$ $x=(x^1,\ldots,x^n)^\top$,
$y=(y^1,\ldots,y^n)^\top.$

We claim that there exists a smooth family $A(t),\ t\in\mathbb R,$
of symmetric $n\times n$ matrices such that $A(t)$ has the same
eigenvalues as $R_\Lambda(t)$ and
$$
\Lambda(t)=\Gamma(0,t)
\left(\begin{smallmatrix}\mathbb R^n\\0\end{smallmatrix}\right),
\quad \Lambda^\circ(t)=\Gamma(0,t)
\left(\begin{smallmatrix}0\\ \mathbb R^n\end{smallmatrix}\right),
\quad \forall t\in\mathbb R
$$
in the fixed coordinates, where $\Gamma(t,\tau)$ satisfies (25).
Indeed, let $e_i(t),\ i=1,\ldots,n,$ be solutions to the
structural equations (23). Then
$$
\Lambda(t)=span\{e_1(t),\ldots,e_n(t)\},\quad
\Lambda^\circ(t)=span\{\dot e_1(t),\ldots,\dot e_n(t)\}.
$$
Moreover,
$\ddot e_i(t)=-\sum\limits_{i=1}^na_{ij}(t)e_j(t)$, where
$A(t)=\{a_{ij}(t)\}_{i,j=1}^n$ is the matrix of the operator
$R_\Lambda(t)$ in the `moving' basis $e_1(t),\ldots,e_n(t)$. Lemma
I.13 implies that $\langle
e_i(t),e_j(t)\rangle_{\dot\Lambda(t)}=\sigma(e_i(t),\dot
e_j(t))=\delta_{ij}$. In other words, the Euclidean structure
$\langle\cdot,\cdot\rangle_{\dot\Lambda(t)}$ has unit matrix in
the basis $e_1(t),\ldots,e_n(t)$. Operator $R_\Lambda(t)$ is
self-adjoint for the Euclidean structure
$\langle\cdot,\cdot\rangle_{\dot\Lambda(t)}$ (see Propositon
II.9). Hence matrix $A(t)$ is symmetric.

Let $e_i(t)=\left(\begin{smallmatrix}x_i(t)\\y_i(t)\end{smallmatrix}\right)\in\mathbb
R^n\times\mathbb R^n$ in the fixed coordinates. Make up $n\times
n$-matrices $X(t)=(x_1(t),\ldots,x_n(t))$,
$Y(t)=(y_1(t),\ldots,y_n(t))$ and a $2n\times 2n$-matrix
$\left(\begin{smallmatrix}X(t) & \dot X(t)\\Y(t) & \dot
Y(t)\end{smallmatrix}\right).$ We have
$$
\frac d{dt}\begin{pmatrix}X & \dot X\\Y & \dot
Y\end{pmatrix}(t)=\begin{pmatrix}X & \dot X\\Y & \dot
Y\end{pmatrix}(t)\begin{pmatrix}0 & -A(t)\\I &
0\end{pmatrix}, \quad \begin{pmatrix}X & \dot X\\Y & \dot
Y\end{pmatrix}(0)=\begin{pmatrix}I & 0\\0 &
I\end{pmatrix}.
$$
Hence $\left(\begin{smallmatrix}X & \dot X\\Y & \dot
Y\end{smallmatrix}\right)(t)=\Gamma(t,0)^{-1}=\Gamma(0,t)$.

Let now $\Lambda(\cdot)$ be the Jacobi curve,
$\Lambda(t)=J_{z_0}(t)$. Set $\xi_i(z_t)=e_*^{t\vec h}e_i(t)$,
$\eta_i(z_t)=e_*^{t\vec h}\dot e_i(t)$; then
$$
\xi_1(z_t),\ldots,\xi_n(z_t),\eta_1(z_t),\ldots,\eta_n(z_t) \eqno
(26)
$$
is a Darboux basis of $T_{z_t}N$, where
$J_{z_t}(0)=span\{\xi_1(z_t),\ldots,\xi_n(z_t)\}$,
$J^\circ_{z_t}(0)=span\{\eta_1(z_t),\ldots,\eta_n(z_t)\}$.
Moreover, the basis (26) is orthonormal for the inner product
$\langle\cdot,\cdot\rangle_{z_t}$ on $T_{z_t}N$.

The intrinsic nature of the structural equation implies the
translation invariance of the construction of the frame (26): if
we would start from $z_s$ instead of $z_0$ and put
$\Lambda(t)=J_{z_s}(t)$, $e_i(0)=\xi_i(z_s)$, $\dot
e_i(0)=\eta_i(z_s)$ for some $s\in\mathbb R$, then
we would obtain $e^{t\vec h}_*e_i(t)=\xi_i(z_{s+t})$,
$e^{t\vec h}_*\dot e_i(t)=\eta_i(z_{s+t})$.

The frame (26) gives us fixed orthonormal Darboux coordinates in
$T_{z_s}N$ for $\forall\,s\in\mathbb R$ and the correspondent
symplectic $2n\times 2n$-matrices $\Gamma_{z_s}(\tau,t)$. We have:
$\Gamma_{z_s}(\tau,t)==\Gamma_{z_0}(s+\tau,s+t)$; indeed,
$\Gamma_{z_s}(\tau,t)
\left(\begin{smallmatrix}x\\y\end{smallmatrix}\right)$ is the
coordinate presentation of the vector
$$
e_*^{(\tau-t)\vec
h}\left(\sum\limits_ix^i\xi^i(z_{s+t})+y^i\eta_i(z_{s+t})\right)
$$
in the basis $\xi_i(z_{s+\tau}),\,\eta_i(z_{s+\tau})$. In
particular,
$$
\left|\Gamma_{z_s}(0,t)\left(\begin{smallmatrix}x\\y\end{smallmatrix}\right)
\right|=\left\|e_*^{-t\vec
h}\left(\sum\limits_ix^i\xi^i(z_{s+t})+y^i\eta_i(z_{s+t})\right)\right\|_{z_s}. \eqno
(27)
$$
Recall that $
\xi_1(z_\tau),\ldots,\xi_n(z_\tau),\eta_1(z_\tau),\ldots,\eta_n(z_\tau)$
is an orthonormal frame for the scalar product
$\langle\cdot,\cdot\rangle_{z_\tau}$ and
$\|\zeta\|_{z_\tau}=\sqrt{\langle\zeta,\zeta\rangle}_{z_\tau}$.

We introduce the notation :
$$
\lfloor W\rfloor_{z_s}=\left\{\sum\limits_ix^i\xi^i(z_s)+y^i\eta_i(z_s):
\left(\begin{smallmatrix}x\\y\end{smallmatrix}\right)\in W\right\},
$$
for any $W\subset\mathbb R^n\times\mathbb R^n$. Let
$C^\pm_{\Gamma_{z_0}}$ be the cones from Lemma II.14. Then
$$
e_*^{-\tau\vec
h}\lfloor\Gamma_{z_s}(0,t)C^\pm_{\Gamma_{z_0}}\rfloor_{z_{s-\tau}}=
\lfloor\Gamma_{z_{s-\tau}}(0,t+\tau)C^\pm_{\Gamma_{z_0}}\rfloor_{z_{s-\tau}},
\quad\forall\, t,\tau,s.  \eqno (28)
$$
Now set $K^+_{\Gamma_{z_s}}=\bigcap\limits_{t\ge
0}C^+_{\Gamma_{z_0}}$, $K^-_{\Gamma_{z_s}}=\bigcap\limits_{t\le
0}C^-_{\Gamma_{z_0}}$ and $\Delta^\pm_{z_s}=\lfloor
K^\mp_{\Gamma_{z_s}}\rfloor_{z_s}$. Corollary II.6 implies that
$\Delta^\pm_{z_s}$ are Lagrangian subspaces of $T_{z_s}N$.
Moreover, it follows from (28) that
$e_*^{t\vec h}\Delta^\pm_{z_s}=\Delta^\pm_{z_{s+t}}$, while (28)
and (27) imply inequalities (24).

This finishes the proof of the part of Proposition II.10 which concerns
Jacobi curves $J_z(t)$. We leave to the reader a simple
adaptation of this proof to the case of reduced Jacobi curves
$J^h_z(t). \qquad\square$

\medskip\noindent {\bf Remark.} Constant $\alpha$ depends, of
course, on the Riemannian structure on $N$. In the case of the
special Riemannian structure defined at the beginning of the proof
of Proposition II.10 this constant depends only on the upper and
lower bounds for the eigenvalues of the curvature operators and
reduced curvature operators correspondently (see Lemma II.14 and
further arguments).

\medskip Let $e^{tX} \ ,t\in\mathbb R$ be the flow generated by the the vector
field $X$ on a manifold $M$. Recall that a compact invariant
subset $W\subset M$ of the flow $e^{tX}$ is called a hyperbolic
set if there exists a Riemannian structure in a neighborhood of
$W$, a positive constant $\alpha$, and a splitting $T_zM=E^+_z\oplus
E^-_z\oplus \mathbb RX(z),\ z\in W$, such that $X(z)\ne 0,
\ e^{tX}_*E^\pm_z=E^\pm_{e^{tX}(z)}$, and $\|e^{\pm
tX}_*\zeta^\pm\|\ge e^{\alpha t}\|\zeta^\pm\|,\ \forall t\ge 0,\
\zeta^\pm\in E^\pm_z$. Just the fact some invariant set is
hyperbolic implies a rather detailed information about asymptotic
behavior of the flow in a neighborhood of this set (see \cite{kh}
for the introduction to Hyperbolic Dynamics). The flow $e^{tX}$ is
called an Anosov flow if the entire manifold $M$ is a hyperbolic
set.

The following result is an immediate corollary of Proposition
II.10 and the above remark.

\begin{theorem} Let $h$ be a regular monotone Hamiltonian on $N$,
$c\in\mathbb R$, $W\subset h^{-1}(c)$ a compact invariant set of
the flow $e^{t\vec h},\ t\in\mathbb R$, and $d_zh\ne 0,\
\forall z\in W$. If $\vec h$ has a negative reduced curvature at
every point of $W$, then $W$ is a hyperbolic set of the flow $e^{t\vec
h}\bigr|_{h^{-1}(c)}. \qquad\square$
\end{theorem}

This theorem generalizes a classical result about geodesic flows
on compact Riemannian manifolds with negative sectional
curvatures. Indeed, if $N$ is the cotangent bundle of a Riemannian
a Riemannian manifold and $e^{t\vec h}$ is the geodesic flow, then
negativity of the reduced curvature of $\vec h$ means simply
negativity of the sectional Riemannian curvature. In this case,
the Hamiltonian $h$ is homogeneous on the fibers of the cotangent
bundle and the restrictions $e^{t\vec h}\bigr|_{h^{-1}(c)}$ are
equivalent for all $c>0$.

The situation changes if $h$ is the energy function of a general
natural mechanical system on the Riemannian manifold. In this
case, the flow and the reduced curvature depend on the energy
level. Still, negativity of the sectional curvature implies
negativity of the reduced curvature at $h^{-1}(c)$ for all
sufficiently big $c$. In particular, $e^{t\vec h}\bigr|_{h^{-1}(c)}$
is an Anosov flow for any sufficiently big $c$; see \cite{ac,acz}
for the explicit expression of the reduced curvature in this case.

Theorem II.3 concerns only the reduced curvature while the next
result deals with the (not reduced) curvature of $\vec h$.

\begin{theorem} Let $h$ be a regular monotone Hamiltonian and $W$
a compact invariant set of the flow $e^{t\vec h}$. If $\vec h$ has
a negative curvature at any point of $W$, then $W$ is a finite set and each
point of $W$ is a hyperbolic equilibrium of the field $\vec h$.
\end{theorem}
{\bf Proof.} Let $z\in W$; the trajectory $z_t=e^{t\vec h}(z),\
t\in\mathbb R$, satisfies conditions of Proposition II.10. Take
the correspondent splitting
$T_{z_t}N=\Delta^+_{z_t}\oplus\Delta^-_{z_t}$. In particular, $\vec
h(z_t)=\vec h^+(z_t)+\vec h^-(z_t)$, where $\vec
h^{\pm}(z_t)\in\Delta^\pm_{z_t}$.

We have $e^{\tau\vec h}_*\vec h(z_t)=\vec h(z_{t+\tau})$. Hence
$$
\|\vec h(z_{t+\tau})\|=\|e^{\tau\vec h}_*\vec h(z_t)\|\ge
\|e^{\tau\vec h}_*\vec h^+(z_t)\|-\|e^{\tau\vec h}_*\vec
h^-(z_t)\|
$$
$$\ge e^{\alpha\tau}\|\vec h^+(z_t)\|-e^{-\alpha\tau}\|\vec
h^-(z_t)\|, \quad \forall\tau\ge 0.
$$
Compactness of $\overline{\{z_t:t\in\mathbb R\}}$ implies that
$\vec h^+(z_t)$ is uniformly bounded; hence $\vec h^+(z_t)=0$.
Similarly, $\|\vec h(z_{t-\tau}\|\ge e^{\alpha\tau}\|\vec h^-(z_t)\|-
e^{-\alpha\tau}\|\vec h^+(z_t)\|$ that implies the equality $\vec
h^-(z_t)=0$. Finally, $\vec h(z_t)=0$. In other words, $z_t\equiv z$
is an equilibrium of $\vec h$ and
$T_zN=\Delta^+_z\oplus\Delta^-_z$ is the splitting of $T_zN$ into
the repelling and attracting invariant subspaces for the
linearization of the flow $e^{t\vec h}$ at $z$. Hence $z$ is a
hyperbolic equilibrium; in particular, $z$ is an isolated
equilibrium of $\vec h. \qquad \square$

We say that a subset of a finite dimensional manifold is bounded
if it has a compact closure.
\begin{corollary} Assume that $h$ is a regular monotone Hamiltonian
and $\vec h$ has everywhere negative curvature. Then any bounded
semi-trajectory of the system $\dot z=\vec h(z)$ converges to an
equilibrium with the exponential rate while another
semi-trajectory of the same trajectory must be unbounded.
$\quad\square$
\end{corollary}

Typical Hamiltonians which satisfy conditions of Corollary II.7 are
energy functions of natural mechanical systems in $\mathbb R^n$
with a strongly concave potential energy. Indeed, in this case, the
second derivative of the potential energy is equal to the matrix
of the curvature operator in the standard Cartesian coordinates (see Sec.~15).

\end{document}